\documentclass[11pt]{article}
\usepackage{amsmath,amssymb,amsthm, geometry,booktabs}
\usepackage{enumitem}

\usepackage{authblk}

\usepackage[normalem]{ulem}  	
\usepackage{changes}  		

\newcommand{\subalpha}{{\underline\alpha}}
\newcommand{\overs}{{\overline s}}
\newcommand{\overd}{{d_\infty}}

\DeclareMathOperator*{\esssup}{ess\,sup}

\newcommand{\dxmax}{{\Delta x}_{max}}
\newcommand{\dxmin}{{\Delta x}_{min}}

\newcommand{\graphics}[1]{}
\usepackage{subcaption}
\usepackage{caption}
\usepackage{float}

\usepackage{xcolor}
\usepackage{amsfonts}
\usepackage{amsmath,amssymb}
\usepackage{graphicx}
\usepackage{epstopdf}

\usepackage{algorithm}
\usepackage{algpseudocode}

\usepackage{caption}
\usepackage{subcaption}

\usepackage{stmaryrd}

\ifpdf
  \DeclareGraphicsExtensions{.eps,.pdf,.png,.jpg}
\else
  \DeclareGraphicsExtensions{.eps}
\fi

\newcommand{\R}{\mathbb{R}}
\newcommand{\N}{\mathbb{N}}
\newcommand{\Z}{\mathbb{Z}}

\newcommand{\KH}{[\cH]} 

\newcommand{\dt}{{\Delta t}}
\newcommand{\dx}{{\Delta x}}
\newcommand{\conv}{\rightarrow}
\newcommand{\converge}{\rightarrow}
\newcommand{\cB}{{\mathcal B}}

\newcommand{\cG}{{\mathcal G}}
\newcommand{\cF}{{\mathcal F}}
\newcommand{\cH}{{\mathcal H}}
\newcommand{\cJ}{{\mathcal J}}

\newcommand{\cR}{{\mathcal R}}
\newcommand{\cS}{{\mathcal S}}
\newcommand{\cV}{{\mathcal V}}
\newcommand{\ma}{\alpha}
\newcommand{\mb}{\beta}

\newcommand{\mk}{\kappa}
\newcommand{\ml}{\lambda}

\newcommand{\mO}{\Omega}
\newcommand{\hD}{{\nabla_G}}
\newcommand{\be}{\begin{eqnarray}}
\newcommand{\ee}{\end{eqnarray}}
\newcommand{\beno}{\begin{eqnarray*}}
\newcommand{\eeno}{\end{eqnarray*}}

\newcommand{\COMMENT}[1]{}   

\newtheorem{definition}{Definition}
\newtheorem{lem}{Lemma}
\newtheorem{corol}{Corollary}

\title{Solving Hamilton-Jacobi equations by residual minimization 
of monotone finite-difference discretizations
}

\author[1]{O.~Bokanowski\thanks{\texttt{olivier.bokanowski@u-paris.fr}}}
\author[2]{C.~Esteve-Yag\"ue\thanks{\texttt{c.esteve@ua.es}}}
\author[3]{R.~Tsai\thanks{\texttt{ytsai@math.utexas.edu}}}

\affil[1]{Laboratoire Jacques-Louis Lions, Universit\'e Paris Cit\'e, France}
\affil[2]{Departamento de Matem\'aticas, Universidad de Alicante, Spain}
\affil[3]{Department of Mathematics and Oden Institute, The University of Texas at Austin, USA}

\date{\today}

\theoremstyle{plain}
\newtheorem{theorem}{Theorem}

\begin{document}
\maketitle

\begin{abstract}
We introduce a method for solving Hamilton--Jacobi equations, both inviscid and viscous, by minimizing the squared residuals of monotone finite-difference discretizations on grids of varying resolution. The method is designed to leverage neural networks and modern GPUs to solve these equations in higher dimensions; consequently, the setting for our analysis is the minimization of the residual functionals via gradient-based optimization. We establish a well-posedness theory for this approach: any critical point of the finite-difference loss solves the monotone scheme together with the prescribed Dirichlet boundary conditions, and the error of an approximation is controlled by its residual. We then derive the rate of convergence of the gradient flow that minimizes the residual in several settings, noting that the rate may depend on the grid resolution and the domain dimension. Building on this foundation, we propose a multi-level training algorithm that exploits the faster convergence available on the coarser grids of the discretization. Combined with the convergence theorem of Barles and Souganidis for monotone and consistent schemes, our results guarantee convergence to the unique viscosity solution as the grid is refined. We illustrate the approach on eikonal equations, level-set problems, and a Hamilton--Jacobi--Isaacs equation arising from a stochastic differential game, in dimensions up to eight.
\end{abstract}


\newpage

\tableofcontents

\section{Introduction}
\label{sec:intro}

Hamilton--Jacobi (HJ) equations are ubiquitous in mathematical physics and control theory, governing systems ranging from level-set dynamics and front propagation to optimal control, differential games, and mathematical finance.
We consider fully nonlinear equations of the general form
\begin{equation}\label{eq:hj}
H(x, u, Du, D^2 u) = 0, \qquad x\in\Omega \subset \mathbb{R}^d,
\end{equation}
with boundary data $u|_\Gamma = g$, and the time-dependent version
\begin{equation}\label{eq:time-hj}
\partial_t u + H(x, u, Du, D^2 u) = 0, \qquad x\in\Omega \subset \mathbb{R}^d,~t>0,
\end{equation}
with initial condition $u(x,0)=u_0(x)$ and boundary data $u|_\Gamma = g$.
Here $Du$ and $D^2 u$ denote the gradient and the Hessian of $u$, respectively.
This formulation encompasses first-order (inviscid) Hamilton--Jacobi equations---including the eikonal equation $c(x)\|Du\|=1$ and level-set equations---as well as second-order elliptic or parabolic problems.
Among the latter are the Hamilton--Jacobi--Bellman (HJB) and Hamilton--Jacobi--Isaacs (HJI) equations arising in stochastic optimal control and differential games.
In these problems a diffusion matrix $a(x)=\frac{1}{2}\sigma(x)\sigma(x)^\top\geq 0$ encodes the noise structure of the underlying process, and $H$ involves terms of the form $-\mathrm{tr}(a(x)\,D^2 u)$.

The solution of \eqref{eq:hj}--\eqref{eq:time-hj} is generally non-smooth, particularly in the first-order (inviscid) case; even when second-order terms are present, solutions may lack classical regularity when the diffusion is degenerate.
The theory of \textit{viscosity solutions}, introduced by Crandall, Evans, and Lions~\cite{crandall1984some, CIL92}, provides the appropriate weak solution framework that handles the inherent non-smoothness and selects the physically relevant solution through a comparison principle.
Existence, uniqueness, and stability results for viscosity solutions are now well established for broad classes of degenerate elliptic and parabolic equations; see~\cite{CIL92, fleming2006controlled, bardi2008optimal} and the references therein.

\paragraph{The discretized problem.}
We consider finite-difference discretizations of \eqref{eq:hj} on a set of nodes $\{x_j\}_{j\in V}$, forming a grid. The nodes are split into interior nodes $\{x_j\}_{j\in \cJ}\subset \Omega$ and boundary nodes $\{x_j\}_{j\in\cB}\subset \Gamma$, so that $V = \cJ \sqcup\cB$.
The discrete problem takes the form
\begin{equation}
\label{eq: disc eq intro}
    \cH \left( x_j, u_j, \nabla_G u_j \right)=0,\quad j\in \cJ,
    \qquad \text{and} \qquad u_j = g(x_j), \quad j\in \cB,
\end{equation}
where the unknown $u = \{u_j\}_{j\in V}\in \R^V$ is a grid function, and
\begin{equation}
\label{finite-differences graph}
\nabla_G u_j :=  \left\{\dfrac{u_j-u_k}{\Delta x_{jk}}\right\}_{k\in\cV_j}
\end{equation}
is the vector of finite-differences approximating the directional derivatives of $u$ at $x_j$ in the directions $\{ (x_j - x_k)/\dx_{jk} \}_{k\in \cV_j}$, with $\cV_j$ denoting the set of neighbors of $x_j$ and $\Delta x_{jk}>0$ the distance between $x_j$ and $x_k$.
The function $\cH$ is the so-called \emph{numerical Hamiltonian}, and is assumed to be consistent with a first-order differential operator $H(x, u, D u)$.
The same framework applies to second-order PDEs such as \eqref{eq:hj}, provided the function $\cH(x_j, u, \nabla_G u)$ also depends on the discretization steps $\{ \dx_{jk} \}_{k\in \cV_j}$.
Note that, after discretization, both first-order terms (involving $Du$) and second-order terms (involving $D^2 u$) in the continuous equation~\eqref{eq:hj} are absorbed into the dependence of $\cH$ on $\nabla_G u_j$. For instance, in a uniform Cartesian grid, the standard central-difference approximation of the Laplacian can be written as $\Delta_G u_j = \sum_{k\in \cV_j} \frac{u_j -u_k}{\dx^2}$.
The numerical Hamiltonian $\cH$ thus encodes the full structure of $H(x, u, Du, D^2 u)$, including any diffusive terms.
We focus on the class of \emph{monotone} discretizations: we require that $\cH(x_j, u, p)$ is non-decreasing in $u$ and in each component $p_k$ of $p$ (see hypotheses {\bf (H2)}--{\bf(H4)} below for the precise conditions).

\paragraph{Convergence of monotone schemes.}
A fundamental result is the convergence theorem of Barles and Souganidis~\cite{barles1991convergence}: if a numerical scheme is \emph{monotone}, \emph{stable}, and \emph{consistent}, then its solutions converge uniformly to the unique viscosity solution as the mesh is refined.
Monotonicity requires $\cH$ to be non-decreasing in~$u_j$ and $\nabla_G u_j$; stability requires uniform bounds for the solution to the discrete problem; and consistency requires
$$
  \cH(x, \phi(x), \nabla_G \phi(x)) \longrightarrow H(x, \phi(x), D\phi(x), D^2\phi(x))
  \qquad \text{as } \Delta x \to 0,
$$
for all smooth $\phi$.
If these three conditions are met, and the continuous equation satisfies a comparison principle, then the solution $\bar u$ of the discrete equation $\cH (x, u, \nabla_G u) = 0$ converges to the viscosity solution of \eqref{eq:hj} as $\Delta x \to 0$.

This convergence theorem provides the theoretical foundation for the present work.
Our goal is to solve the discrete equation \eqref{eq: disc eq intro}  by minimizing a least-squares residual, and Barles--Souganidis framework guarantees that---provided $\cH$ is monotone and consistent---the solutions so obtained converge to the correct viscosity solution as the discretization is refined.
Therefore, this paper concentrates on the analytical questions that arise from solving $\cH (x,u, \nabla_G u ) = 0$ via residual minimization, rather than on the convergence of $\cH$ to the continuous operator, which is handled by the Barles--Souganidis theory.

\paragraph{Residual minimization and the loss functional.}
Our main motivation is to develop numerical algorithms for PDE problems of the form \eqref{eq:hj} and \eqref{eq:time-hj} in high-dimensional domains.
Rather than solving the nonlinear system \eqref{eq: disc eq intro} by classical iterative methods (which are infeasible in high-dimensional settings),
we formulate it as an optimization problem.
We define the residual operator $\cR:\R^V\to\R^V$ by $\cR(u)_j = \cH(x_j, u_j, \nabla_G u_j)$ for $j\in\cJ$ and $\cR(u)_j = \lambda_b(u_j - g_j)$ for $j\in\cB$, and consider the least-squares loss functional
\begin{equation}
\label{eq:loss func def intro}
    L(u) := \| \cR (u)\|_{\ell^2(V)}^2 = \sum_{j\in\cJ} \left( \cH (x_j, u_j, \nabla_G u_j )\right)^2 + \lambda_b^2 \sum_{j\in\cB} \big(u_j - g(x_j)\big)^2.
\end{equation}
The discrete solution is recovered as the minimizer of $L$, and in practice, we approximate it using a neural network $\Phi(x;\theta)$ trained by stochastic gradient descent on $L(\Phi(x;\theta))$.
This approach is motivated by the need for mesh-free, GPU-friendly solvers that can scale to high dimensions, circumventing the exponential complexity of grid-based discretizations.

A common alternative is to minimize the residual of the \emph{continuous} PDE directly, as in Physics-Informed Neural Networks (PINNs) \cite{raissi2019physics}.
However, for nonlinear and nonconvex equations such as \eqref{eq:hj}, minimizing the continuous residual does not guarantee convergence to the viscosity solution, since multiple weak solutions---that solve the PDE at almost every point in $\Omega$---may produce identical zero-residual values.
Our approach instead minimizes the residual of a \emph{monotone numerical scheme}, thereby inheriting the convergence guarantees of the Barles--Souganidis framework once the discrete equation is solved exactly or to sufficient accuracy.

\paragraph{Key analytical questions.}
Solving the discrete equation $\cH(x_j, u_j, \nabla_G u_j) = 0$ by minimizing the least-squares functional $L(u) = \|\cR(u)\|_{\ell^2 (V)}^2$ raises several fundamental analytical issues:
\begin{enumerate}
\item \emph{Critical points vs.\ solutions.} Does every critical point of $L(\cdot)$ solve $\cH = 0$? Since $L(\cdot)$ is not convex in general, one must ensure that gradient-based optimization cannot become trapped at a spurious critical point that does not correspond to a solution of the discrete scheme.

\item \emph{Stability and error estimates.} If $L(\cdot)$ is minimized only approximately---as is inevitable in practice with stochastic gradient descent---how does the approximation error $\|u - \bar u\|_{\ell^\infty(V)}$ relate to the magnitude of the residual $\|\cR(u)\|_{\ell^\infty(V)}$? Is it possible to control the $\ell^2$-error, $\| u- \bar u\|_{\ell^2(V)}$, by the $\ell^2$-norm of the residual, which is indeed $\sqrt{L(u)}$?

\item \emph{Conditioning and convergence of gradient descent.} How does the conditioning of the gradient system $\nabla L(u) = 2D\cR(u)^\top \cR(u)$ depend on the discretization parameters? In particular, how do the small eigenvalues of the Jacobian $D\cR(u)$ scale with $\Delta x$, and what are the implications for the convergence rate of gradient descent?
\end{enumerate}
We address all three questions in this paper.

\paragraph{Main contributions.}

The framework studied here was introduced in \cite{esteve2025finite} for Lax--Friedrichs-type discretizations on Cartesian grids, where it was shown that, if the diffusion dominates the nonlinear part of $\cH$, then every critical point of $L$ solves the discrete equation---answering question~1 for that case. The present work removes both restrictions, treating general monotone schemes on unstructured meshes, and addresses all three questions.
\begin{itemize}
\item \emph{Well-posedness} (question~1; Theorems~\ref{thm:exist-uniqueness-scheme-gen} and \ref{thm: euler map is a contraction}). Under {\bf (H1)}--{\bf (H2)} and either a properness condition {\bf (H3)} or a new uniform-ellipticity condition {\bf (H4)}, the residual operator $\cR:\R^V\to\R^V$ is a homeomorphism, via a Lipschitz extension of Hadamard's theorem due to Pourciau~\cite{pourciau-82}. In particular, every critical point of $L(\cdot)$ solves $\cR(u)=0$, i.e., there are no spurious minima.
\item \emph{A posteriori error estimates} (question~2; Theorems~\ref{thm:stability-H3}, \ref{thm:dim-robust error estimate proper case}, \ref{thm:stability-H4}, and \ref{thm:  viscous HJ}). The error $\|u-\bar u\|_{\ell^\infty(V)}$ is controlled by $\|\cR(u)\|_{\ell^\infty(V)}$ with explicit, computable constants.
Under {\bf (H3)}, and when the discretization parameter $\dx$ is sufficiently big, we obtain a dimension-robust $\ell^2$ estimate, with a constant that depends on local stencils, but not on the total number of grid points. A similar dimension-robust $\ell^2$ estimate is obtained for schemes which are consistent with second-order semilinear elliptic equations.
\item \emph{Conditioning and convergence} (question~3; Theorems~\ref{thm:dim-robust error estimate proper case}, \ref{thm:  viscous HJ}). Under the same conditions as when dimension-robust $\ell^2$-error estimates apply (i.e. $\dx$ sufficiently big), we obtain a Polyak--\L{}ojasiewicz inequality. This ensures that the gradient flow that minimizes $L(\cdot)$ converges linearly to the global minimizer. We also provide explicit upper bounds for the condition number of the residual's differential, which serve to estimate the convergence rate of the associated gradient flow.
Our estimates show that the condition number is smaller for larger values of $\dx$ and deteriorates when $\dx$ is increased.
This deterioration on finer meshes motivates the multi-level training strategy of Section~\ref{sec:numerical algos}.
\end{itemize}

We also extend the analysis to time-dependent problems (Theorems~\ref{thm:exist-uniqueness-scheme-time-implicit} and \ref{thm:error-estimates-scheme-time-implicit}), discretizing the time derivative by a single implicit or explicit Euler step, which recasts the evolution as the minimization of one space-time residual, for which {\bf (H1)}--{\bf (H2)} suffice. The implicit discretization is distinguished in the $\ell^\infty$ setting: the error is controlled by the residual, with no accumulation in time, for every time step, whereas the explicit one requires a CFL-type condition. In the $\ell^2$ setting (the norm defining the loss $L(\cdot)$) dimension-robust error estimates require the CFL-type condition for both discretizations, and, under it, the two residual minimization problems are essentially equivalent, with comparable training efficiency (Section~\ref{sec:time}). Although the implicit Euler discretization is unconditionally stable, without the CFL-type restriction, the associated $\ell^2$-error estimate and the rate of convergence during training suffer from the curse of dimensionality.

The framework applies to Lax--Friedrichs and upwind schemes, and extends to second-order elliptic and parabolic PDEs on graphs with monotone graph Laplacians. Combined with the Barles--Souganidis theorem, the discrete solutions obtained by residual minimization converge to the viscosity solution of \eqref{eq:hj} as $\dx\to0$. The numerical experiments in Section~\ref{sec:numerics} demonstrate the scalability of the approach: eikonal equations in high dimensions, time-dependent level-set problems with obstacles, and a Hamilton--Jacobi--Isaacs equation with second-order diffusion from a stochastic differential game.

\paragraph{Related work.}
Classical numerical methods for solving Hamilton--Jacobi equations, such as Fast Marching Methods \cite{tsitsiklis2002efficient,sethian1999fast, sethian2001ordered} and Fast Sweeping Methods \cite{zhao2005fast, tsai2003fast, kao2005fast}, rely on discretizing the PDE on a grid and are effective in low dimensions but fundamentally infeasible for high-dimensional problems due to exponential grid complexity.
Grid-free algorithms based on classical representation formulas (Hopf--Lax, Lax--Oleinik) and fast optimization have been developed in \cite{darbon2016algorithms, chow2019algorithm, chen2024lax}; these compute the solution at query points with complexity that does not grow exponentially.
Recently, \cite{park2025neural} introduced a neural implicit solution formula that generalizes the Hopf and Lax formulas via deep learning.

\paragraph{Organization.}
Section~\ref{sec: main results stationary case} sets up the framework of monotone schemes on geometric graphs and proves well-posedness (Theorem~\ref{thm:exist-uniqueness-scheme-gen}) under {\bf (H3)} and {\bf (H4)}.
Section~\ref{sec: stability} derives $\ell^\infty$ {\it a posteriori} error estimates, a dimension-robust $\ell^2$ estimate with a Polyak--\L{}ojasiewicz inequality, and upper bounds for the condition number of the linearized problem.
Section~\ref{sec:time} extends the framework to time-dependent problems and compares implicit and explicit time stepping.
Section~\ref{sec:monotone-hamiltonians} verifies that standard schemes satisfy our hypotheses and discusses the Kruzhkov transform as a technique to enforce the properness condition in {\bf (H3)}.
Section~\ref{sec:numerical algos} presents the multi-level algorithm with spectral analysis.
Section~\ref{sec:numerics} provides computational examples.

\section{Monotone schemes and their residual functionals}
\label{sec: main results stationary case}

We use the notation introduced in Section~\ref{sec:intro}: the nodes $\{x_j\}_{j\in V}$, with $V$ being a finite index set, are the vertices of a directed graph $G=(V,E)$, where $(j,k)\in E$ whenever $k\in \cV_j$; we assume $j\notin \cV_j$ for all $j$.
The interior and boundary index sets are
$$
\cJ := \{ j\in V : \  \cV_j \neq \emptyset\}
\quad \text{and} \quad
\cB := \{ j\in V : \  \cV_j = \emptyset\},
$$
both assumed non-empty.
It is typical to assume $\{x_j\}_{j\in \cJ}\subset \Omega$ and $\{ x_j\}_{j\in \cB}\subset \partial\Omega$, though our analysis does not require it.
The space of grid functions $u:V\conv \R$ is denoted $\R^V$.
The discrete problem and the finite-difference gradient $\nabla_G u_j$ are as in~\eqref{eq: disc eq intro}--\eqref{finite-differences graph};
we restate the problem here for ease of reference:
\begin{equation}\label{eq:scheme-steady-1}
\begin{cases}
  \cH(x_j,u_j, \nabla_Gu_j)=0, & j\in \cJ \\
 u_j=g_j, & j\in \cB
\end{cases}
\end{equation}
for given boundary values $\{ g_j \}_{j\in \cB}$, where the numerical Hamiltonian at each node is a function $\cH(x_j,\cdot,\cdot):\R\times\R^{|\cV_j|}\converge\R$ which depends on $u_j$ and $\nabla_Gu_j$.

\subsection{Hypotheses on the numerical Hamiltonian}

We define the following conditions on $\cH$:
\begin{itemize}
 \item[\bf (H1)]
  for all $x_j$, the function $(u,p)\mapsto\cH(x_j,u,p)$ is Lipschitz continuous, i.e. there exist two constants $\KH_u, \KH_p>0$ (independent of $x_j$), such that
  $$
  | \cH (x_j, u, p) - \cH (x_j, v, p)| \leq \KH_u |u-v|, \qquad
  \forall u,v\in \R \quad \text{and} \quad
  p\in \R^{|\cV_j|}.
  $$
  $$
  | \cH (x_j, u, p) - \cH (x_j, u, q)| \leq \KH_p \,  \| p - q\|_\infty, \qquad
  \forall u\in \R \quad \text{and} \quad
  p,q\in \R^{|\cV_j|}.
  $$
  We denote $\KH := \max (\KH_u,\, \KH_p)$.
 \item[\bf (H2)]
  $\cH$ is monotone in the following sense: $\forall j\in \cJ$,
  $\forall u,v\in \R$, and $\forall p,q \in \R^{|\cV_j|}$
  $$
  \mbox{ $p\leq q$ (component wise) $\Rightarrow$ $\cH(x_j, u, p)\leq \cH(x_j, u, q)$. }
  $$
  and
  $$
  \mbox{ $u\leq v$ $\Rightarrow$ $\cH(x_j, u, p)\leq \cH(x_j, v, p)$.}
  $$
\end{itemize}

To ensure well-posedness, we furthermore assume one of the following two conditions---a properness condition {\bf (H3)}, or a weaker uniform-ellipticity condition {\bf (H4)} that also covers the degenerate case $\ml=0$ in which properness fails:
\begin{itemize}
 \item[\bf (H3)]
  $\cH$ is ``proper", that is,
  there exists a constant $\ml>0$ such that $\forall j\in \cJ$, $\forall u, v \in \R$, and $\forall p \in \cR^{|\cV_j|}$,
  $$
    v\geq  u \quad \Longrightarrow\ \quad
    \cH(x_j,v,p)-\cH(x_j,u,p) \geq \ml (v-u).
  $$

 \item[\bf (H4)] $\cH$ is ``uniformly elliptic'', that is,
  \begin{enumerate}[label=(\roman*)]
    \item there exists a constant $\subalpha > 0$  such that
    $\forall j\in \cJ$, $\forall u \in \R^V$:
    \be\label{eq:equiv1}
      \mbox{for all $k\in \cV_j$,}
      \qquad \frac{\partial \cH}{\partial p_k}(x_j,u,p)\geq \subalpha, \quad \mbox{for a.e. $p\in \R^{|\cV_j|}$.}
    \ee
    \item for any $j\in\cJ$, there exists $n\in \N$, a sequence
    $(j_1, \ldots, j_n)\in \cJ^n$ and $j_{n+1}\in\cB$
    such that $j_1 = j$ and $j_{k+1}\in \cV_{j_k}$ for all $k=1,2,\ldots, n$.
    \item for any $j,k\in \cJ$, it holds that $j\in \cV_k$ if and only if $k\in \cV_j$.
  \end{enumerate}
\end{itemize}

The properness condition {\bf (H3)} was introduced in \cite{oberman2006convergent}
to obtain the existence of a fixed point of the Euler map $u- \rho\cR(u)$ for $\rho$ sufficiently small;
it is also related to well-posedness in the continuous setting (see \cite{CIL92}).

Condition {\bf (H4)} introduces a uniform ellipticity property for the numerical Hamiltonian to handle cases where $\ml=0$ may occur. We emphasize that this condition characterizes the discrete scheme rather than the continuous equation; a uniformly elliptic scheme can be consistent with either a first-order or a second-order fully nonlinear PDE (such as those arising from a diffusion term). Under this framework, condition {\bf (H4)}$(i)$ acts as a ``proper''-like condition on the difference quotients, and condition {\bf (H4)}$(ii)$ requires that every interior node can be connected to the boundary through a finite path of neighbors.

\subsection{Residual function and well--posedness}

We define the residual operator $\cR: \R^V \longrightarrow \R^V$ by
\begin{equation}
    \label{residual function def}
    \cR(u)_j :=
    \begin{cases}
         \cH(x_j,u_j, \nabla_Gu_j), & \text{if} \quad j\in \cJ,  \\
         \lambda_b (u_j - g_j) , & \text{if} \quad  j\in \cB,
    \end{cases}
\end{equation}
for some arbitrary constant $\lambda_b>0$.
Note that $u\in \R^V$ is a solution of \eqref{eq:scheme-steady-1} if and only if $\cR (u) = 0$.
We shall prove that under {\bf (H1)}--{\bf (H2)} and ({\bf (H3)} or {\bf (H4)}), $\cR$ is a homeomorphism of $\R^V$, so \eqref{eq:scheme-steady-1} has a unique solution $u=\cR^{-1}(0)$. As a by-product, this solution is the unique critical point of the least-squares functional
\begin{equation}
\label{Lp functional steady state eq}
L(u) : = \frac{1}{|\cJ|} \sum_{j\in \cJ} \left(\cH(x_j,u_j, \nabla_Gu_j) \right)^2
 + \frac{\mu_b}{|\cB|}  \sum_{j\in \cB} \left( u_j - g_j \right)^2,
\end{equation}
with any $\mu_b>0$.
This motivates the approximation of the solution of \eqref{eq:scheme-steady-1} via neural networks trained by stochastic gradient descent on $L(u)$.
The constant $\mu_b$ and the normalizing factors $|\cJ|$ and $|\cB|$ do not play
a fundamental role in our analysis in this section, but, {as we will see, they become relevant in the error estimates and rates of convergence during training that we obtain in the next section.}
 
Our strategy is to show that $\cR$ is a global homeomorphism via Hadamard's theorem~\cite{hadamard1906transformations}: if $\cR: \R^V \to \R^V$ is continuously differentiable and there exists a constant $C>0$ such that $\|(D\cR(u))^{-1}\| \leq C$, $\forall u\in \R^V$, then $\cR$ is a diffeomorphism of $\R^V$ onto itself.
Throughout, we use the $\ell^\infty$-norm $\|x\|_{\ell^\infty(V)} := \sup_{j\in V} |x_j|$ and its subordinate matrix norm.

To cover Lipschitz (not necessarily $C^1$) residuals, we use the generalized differential introduced by Pourciau~\cite{pourciau-82} (see also Clarke~\cite{clarke}).

\begin{definition}[\textit{Generalized differential}]
    \label{def: generalized differential}
Let $\cR: \R^{V}\to \R^{V}$ be a.e.\ differentiable. The generalized differential $\partial \cR(u)$ is the convex hull of all limits of $D\cR(u_n)$ over sequences $(u_n)_n$ of differentiability points converging to~$u$.
Here, $D\cR(u_n)$ denotes the classical Jacobian matrix of $\cR(u_n)$.
\end{definition}

The following generalization of the mean value theorem will be used throughout the paper.

\begin{lem}[generalized mean value,~{\cite[Proposition 2.6.5]{clarke}}]
\label{lem:mvt}
Let $\cR: \R^{V}\to \R^{V}$ be a Lipschitz map.
For any $u,v\in\R^V$ there exists a matrix $A\in co\big(\partial\cR([u,v])\big)$, the convex hull of $\{Z\in\partial\cR(q):q\in[u,v]\}$, such that $\cR(u)-\cR(v)=A\,(u-v)$.
\end{lem}

We start by proving that, under {\bf(H1)}, {\bf(H2)}, and ({\bf(H3)} or ${\bf(H4)}$),
there is a constant $C>0$, independent of $u$, such that $\| A^{-1} \|_{\ell^\infty (V)} \leq C$ for all $A\in \partial \cR(u)$.
The first case (when {\bf (H3)} holds) is dealt with in Lemma~\ref{lem:DR-invertible lambda},
while the second case (when {\bf (H4)} holds) will be treated in Lemmas~\ref{lem:barrier} and \ref{lem:DR-invertible alpha}.
\begin{lem}\label{lem:DR-invertible lambda}
    Assume that the Hamiltonian function $\cH$ satisfies {\bf(H1),(H2),(H3)}.
    Then, for any $u\in \R^V$, any matrix $A\in \partial \cR(u)$
    is $\mu$-diagonally dominant with 
    $\mu = \min \left(  \lambda, \lambda_b \right).$
    Hence, any $A\in \partial\cR(u)$ is invertible with $\| A^{-1}\|_{\ell^\infty (V)} \leq \frac{1}{\mu}$. 
\end{lem}

\begin{proof}
We recall that, for $\mu>0$, a square matrix $A\in \R^{V\times V}$ is said to be $\mu$-diagonally dominant if $|A_{jj}|\geq \mu + \sum_{k\neq j} |A_{jk}|$ for all $1\le j\le V$. 

  If $A$ is $\mu$-diagonally dominant with $\mu > 0$, then $A$ is invertible with $\| A^{-1}\|_\infty \leq \frac{1}{\mu}$. Indeed, for any $x\in \R^V$, we pick $i$ with $|x_i|=\|x\|_{\ell^\infty (V)}$, and obtain
  $$
  |(Ax)_i| = \Big| A_{ii}x_i+\sum_{k\neq i}A_{ik}x_k\Big|
   \geq |A_{ii}|\,|x_i|-\sum_{k\neq i}|A_{ik}|\,|x_k|
   \geq \Big(|A_{ii}|-\sum_{k\neq i}|A_{ik}|\Big)\|x\|_{\ell^\infty (V)}
   \geq \mu\,\|x\|_{\ell^\infty (V)},
  $$
  so $\|x\|_{\ell^\infty (V)}\leq\frac{1}{\mu}\|Ax\|_{\ell^\infty (V)}$. In particular $Ax=0$ forces $x=0$, so $A$ is invertible. For any $y\in \R^V$, applying this inequality to $x = A^{-1}y$, we obtain $\| A^{-1} y \|_{\ell^\infty (V)} \leq \frac{1}{\mu} \| y\|_{\ell^\infty (V)}$, implying that $\| A^{-1}\|_{\ell^\infty (V)} \leq \frac{1}{\mu}$.

  Recall that $\hD u_j:=\{\frac{u_j-u_k}{\dx_{kj}}\}_{k\in \cV_j}$. For any $j\in \cJ$ and $k\in \cV_j$, let us denote by $p_k$ the variable of $\cH(x_j, u_j, \nabla_G u_j)$ 
  associated to the value $p_k=\frac{u_j-u_k}{\dx_{jk}}$.
  If $\cR(\cdot)$ is differentiable at $u$, then $\partial\cR(u)$ consists of the classical Jacobian matrix $D\cR(u)$.
  Then, for any $j\in \cJ$, we can compute the diagonal term in the $j$-th row of the Jacobian as
  \begin{equation}\label{eq:dH(uj)/duj}
    \frac{\partial \cR(u)_j}{\partial u_j} 
    =  \frac{\partial \cH}{\partial u}(x_j,u_j,\hD u_j) 
	+ \sum_{k\in \cV_j} \frac{1}{\dx_{jk}} \frac{\partial \cH}{\partial p_k}(x_j,u_j,\hD u_j), 
  \end{equation}
   whereas the off-diagonal terms read as
   \begin{equation}\label{eq:dH(uj)/duk}
     \frac{\partial \cR(u)_j}{\partial u_k}  = 
	-\frac{1}{\dx_{jk}} \frac{\partial \cH}{\partial p_k}(x_j,u_j,\hD u_j), \qquad \text{for all $k\in \cV_j$,}     
   \end{equation}
   and $\frac{\partial \cR(u)_j}{\partial u_k} =  0$ for all $k\notin \cV_j\cup \{j\}$.
   Because of the monotonicity assumptions \textbf{(H2),(H3)}, it holds that
    \begin{align*}
    \left| \frac{\partial \cR(u)_j}{\partial u_j} \right| & \geq  \lambda + \sum_{k\in \cV_j} \frac{1}{\dx_{jk}} \frac{\partial \cH}{\partial p_k}(x_j,u_j,\hD u_j) = \lambda + \sum_{k\neq j} \left| \frac{\partial \cR(u)_j}{\partial u_k} \right|, \qquad \forall j\in \cJ.
    \end{align*}
    Here we used, in view of {\bf (H2)}, that $\frac{\partial \cH}{\partial p_k}(x_j,u_j,\hD u_j) \geq 0$ for all $k\in \cV_j$.
    Hence, all the rows in the Jacobian $D\cR(u)$ are $\mu$-diagonally dominant with $\mu = \lambda$.

    For any $j\in \cB$, the corresponding row is trivially $\mu$-diagonally dominant with $\mu = \lambda_b$, since 
    $\frac{\partial \cR(u)_j}{\partial u_j} = \lambda_b$ and
    $\frac{\partial \cR(u)_j}{\partial u_k} = 0$ $\forall k\neq j.$
    Hence, $D\cR (u)$ is $\mu$-diagonally dominant with $\mu = \min (\lambda, \lambda_b)$, and therefore, $D\cR (u)$ is invertible with $\| (D\cR (u))^{-1}\|_{\ell^\infty (V)} \leq \frac{1}{\mu}$.

    For the non-differentiable case, we note that the convex combination of $\mu$-diagonally dominant matrices
    is also $\mu$-diagonally dominant. Moreover, the limit of a sequence of $\mu$-diagonally dominant matrices 
    is also $\mu$-diagonally dominant, concluding the proof for the case when $\cR(\cdot)$ is only Lipschitz.
\end{proof}

Now, let us consider the case in which {\bf (H4)} holds. 
This case is more delicate since the Jacobian matrix $D\cR(u)$ 
is no longer diagonally dominant with strict inequality.

The first step to obtain a uniform upper bound for the inverse Jacobian
is the construction of a barrier function $\psi\in \R^V$ such that $D\cR(u) \psi \geq 1$ for all $u\in \R^V$. The barrier $\psi$ must be independent of $u$.

\begin{definition}
[{\bf Notation related to the graph, the mesh, and the numerical Hamiltonian}]
\label{def: graph distance constants}
We define the graph-distance to the boundary as follows:
$$
  d_j := \min \{ n \ : \ \exists (j_1, j_2, \ldots, j_n, j_{n+1})\in \cJ^n \times \cB \quad \text{s.t.} \quad j_1=j, \quad j_{k+1}\in \cV_{j_k} \ \forall k\},
$$
for $j\in \cJ$, and with $d_j := 0$ for $j\in \cB$.

We furthermore introduce the quantities
    $$ \overs := \esssup_{j,u,p}
      \left(\sum_{k\in \cV_j} \frac{\partial\cH}{\partial p_k}(x_j,u,p)\right), 
      \quad \overd:=\max_{j\in \cJ} d_j
    $$
    as well as 
    $$
     \dxmin := \min_{j,\,k} \dx_{jk}, \qquad \dxmax:= \max_{j,\,k} \dx_{jk}.
    $$
\end{definition}
Since $p\conv \cH(x_j,u,p)$ is Lipschitz continuous, with Lipschitz contant $\KH_p$, we can upper bound $\overs$ by
$\overs \leq \KH_p\ \max_j |\cV_j|$, where $\max_j |\cV_j|$ is the maximum degree of the graph. The quantity $\overd$ is the so-called in-radius of the graph and can be estimated using the geometry of the domain $\Omega$. For uniform Cartesian meshes with constant discretization step $\dx$ in a $d$-dimensional domain, we have $\overs \leq 2 d \KH_p $ and $\overd = O( 1/\dx)$.

\begin{lem}\label{lem:barrier}
    Assume that the Hamiltonian function $\cH$
    satisfies {\bf(H1), (H2)} and {\bf (H4)}.
    Then, there exists a grid function $\psi\in \R^V$ such that, for any $u\in \R^V$, it holds that
    $$
    A \psi \geq \mathbf{1} \quad 
    \text{for any} \quad A \in \partial \cR(u),
    $$
    where $\partial \cR(u)$ denotes the generalized differential (see Definition~\ref{def: generalized differential}).
    Moreover, the function $\psi$ can be chosen independently of $u$ and satisfying
    \be\label{eq:psi-bound}
      \|\psi \|_{\ell^\infty (V)} \leq  \dfrac{\dxmax}{\subalpha}
       \left(2 + \frac{\overs}{\subalpha}\frac{\dxmax}{\dxmin}\right)^\overd + \dfrac{1}{\lambda_b}.
    \ee
\end{lem}

\begin{proof}
If we assume that $\cR(\cdot)$ is differentiable at $u$, by the same computations as in \eqref{eq:dH(uj)/duj} and \eqref{eq:dH(uj)/duk}, the Jacobian matrix $D\cR(u)$ defines the linear operator in $\R^V$ such that, for any $w\in \R^V$, the $j$-th component of the vector $D\cR(u) w$ is given by
    \begin{equation}
\label{jacobian general}
 [D\cR(u)w]_j = \begin{cases}
 \beta_j w_j + \displaystyle\sum_{k\in \cV_j} \alpha_{jk} (w_j - w_k), & \text{if} \quad j\in \cJ\\
 \lambda_b w_j, & \text{if} \quad j\in \cB
  \end{cases}
\end{equation}
where 
  \begin{equation}
  \label{eq:alpha_and_beta}
    \alpha_{jk}:= \dfrac{1}{\dx_{jk}} \dfrac{\partial \cH}{\partial p_k} (x_j,u_j, \nabla_G u_j)
    \quad \mbox{and} \quad
    \beta_j := \dfrac{\partial \cH}{\partial u} (x,u_j, \nabla_G u_j).
  \end{equation}
  Due to {\bf (H1)}, {\bf (H2)} and 
  {\bf (H4)}, we have, 
  $\forall j\in \cJ$ and $k\in \cV_j$:
\begin{equation}
\label{estimates beta alpha lemma}
  \beta_j \geq 0, \quad  \alpha_{jk}\geq \dfrac{\subalpha}{\dx_{max}} 
  \quad \text{and} \quad
  \sum_{k\in \cV_j} \alpha_{jk}\leq\dfrac{\overs}{\dx_{min}}.
\end{equation}

If $\cR(\cdot)$ is not differentiable at $u$, then any matrix $A\in \partial \cR(u)$ can be written as $A = \lim_{n\to \infty} A_n$, where $A_n$ is a convex combination of Jacobian matrices $D\cR(u_n)$. Hence, every matrix $A_n$ defines an operator of the form \eqref{jacobian general}, with coefficients $\beta_j$ and $\alpha_{jk}$ satisfying \eqref{estimates beta alpha lemma}.
We shall construct a barrier $\psi\in \R^V$ such that $A_n\psi \geq \mathbf{1}$ for all $n$, so taking the limit, we obtain that $A\psi \geq \mathbf{1}$.

Let
\begin{equation}
\label{supersol psi def}
\psi_j := C(1- e^{-\gamma d_j}) + \dfrac{1}{\lambda_b}, \qquad \forall j\in \cJ \cup \cB,
\end{equation}
for some constants $C, \gamma>0$ to be chosen. Since $\psi_j\geq 0$ and $\mb_j\geq0$, it suffices to verify
\begin{equation}
\label{supersol psi}
\begin{cases}
  \displaystyle\sum_{k\in \cV_j} \alpha_{jk} (\psi_j - \psi_k) \geq  1, & \text{if} \quad j\in \cJ \\
  \lambda_b \psi_j\geq 1 & \text{if} \quad j\in \cB.
\end{cases}
\end{equation}

First, notice that for any $j\in \cB$, $\lambda_b \psi_j = 1$.
Then, for $j\in \cJ$, by the definition of $d_j$, and using {\bf(H4)}(iii) we have
$$
   d_k \leq d_j + 1, \qquad \forall k\in \cV_j.
$$
Since, by {\bf(H4)}(ii), every node is connected to the boundary, some neighbor $k^*\in \cV_j$ lies one step closer to the boundary, i.e.\ $d_{k^*} = d_j - 1$,
and it holds
$$
   \psi_j - \psi_{k^*} = C \left(1-e^{-\gamma d_j} - \left(1-e^{-\gamma(d_j-1)}\right) \right) 
     = C e^{-\gamma d_j} \left( e^\gamma - 1 \right).
$$
On the other hand, for $k\in \cV_j$,
$$
  \psi_j - \psi_k
   \geq C \left( 1 - e^{-\gamma d_j} - \left( 1 - e^{-\gamma (d_j + 1)}\right) \right) 
       = C e^{-\gamma d_j} \left( e^{-\gamma} -1 \right)
   \geq -C e^{-\gamma d_j}.
$$

Assembling everything and using \eqref{estimates beta alpha lemma}, we obtain
\begin{eqnarray*}
   \sum_{k\in \cV_j} \alpha_{jk} (\psi_j - \psi_k) 
   &  =   & \alpha_{jk^*} (\psi_j - \psi_{k^*}) + \sum_{k\in \cV_j\setminus \{k^*\}} \alpha_{jk} (\psi_j - \psi_k) \\ 
   & \geq & \alpha_{jk^\ast} C e^{-\gamma d_j} \left( e^\gamma - 1 \right)  
     - \sum_{k\in \cV_j \setminus \{k^\ast\}} \alpha_{jk} C e^{-\gamma d_j} \\ 
   & \geq &   C e^{-\gamma d_j} \left( \frac{\subalpha}{\dxmax} \left( e^\gamma - 1 \right)  - \frac{\overs}{\dxmin}\right). 
\end{eqnarray*}
Now we choose $\gamma>0$ such that $e^\gamma = 2+\frac{\overs}{\subalpha}\frac{\dxmax}{\dxmin}$, 
and $C := e^{\gamma \overd} \frac{\dxmax}{\subalpha}$. 
Therefore, 
$\sum_{k\in \cV_j} \alpha_{jk} (\psi_j - \psi_k)  \geq C e^{-\gamma \overd} \frac{\subalpha}{\dxmax} \geq 1,$
which leads to \eqref{supersol psi}.
By construction we also have the estimate $\|\psi\|_{\ell^\infty (V)} \leq C + \frac{1}{\lambda_b}$, proving \eqref{eq:psi-bound}.
This concludes the proof of Lemma~\ref{lem:barrier}.
\end{proof}

We can now prove an analogue of Lemma \ref{lem:DR-invertible lambda} for the case when {\bf (H3)} is not satisfied
(that is, we may have $\lambda =0$), but {\bf (H4)} holds.

\begin{lem}\label{lem:DR-invertible alpha}
    Assume that the Hamiltonian function $\cH$ satisfies {\bf(H1), (H2)} and {\bf (H4)}.
    Then, for any $u\in \R^V$, all the matrices $A\in \partial\cR(u)$ 
    are invertible and satisfy
    $$
    \| A^{-1}\|_{\ell^\infty (V)} \leq  
         \frac{\dxmax}{\subalpha}
         \left(2 + \frac{\overs}{\subalpha}\frac{\dxmax}{\dxmin}\right)^\overd + \dfrac{1}{\lambda_b}.
    $$
    Moreover, all the elements in the inverse $A^{-1}$ are non-negative.
\end{lem}

\begin{proof}
Let $A\in \partial \cR(u)$. We can write $A$ as the limit $A = \lim_{n\to \infty} A_n$, where $A_n$ is the convex combination of Jacobian matrices $D\cR(u_n)$. As in the proof of Lemma~\ref{lem:barrier}, for every $n$ the matrix $A_n$ defines a linear operator in $\R^V$ of the form \eqref{jacobian general}, with coefficients $\beta_j$ and $\alpha_{jk}$ satisfying \eqref{estimates beta alpha lemma}.
We shall prove that $A_n$ is invertible for every $n$, with $A^{-1}_n$ having all the elements non-negative, and
$$
  \| A_n^{-1}\|_{\ell^\infty (V)}  \leq \|\psi\|_{\ell^\infty (V)},
$$
where $\psi$ is the barrier function from Lemma~\ref{lem:barrier}.
Since the upper bound is independent of $n$, we deduce that the limit $A$ is invertible with the same upper bound for its inverse. Since all the elements of $A_n^{-1}$ are non-negative, the same property holds for the limit $A^{-1}$.

We prove the upper bound in two steps.

\textit{Step 1: Comparison principle.} First we verify that the operator associated to the matrix $A_n$ satisfies the discrete comparison principle, i.e.
$$
\forall w_1, w_2\in \R^V, \qquad A_n w_1 \geq A_n w_2 \quad  \text{implies} \qquad w_1 \geq w_2.
$$
This is a classical result for Laplacian-type matrices, which we recall here for the sake of completeness.

Let $w_1, w_2\in \R^V$ be such that $A_n w_1 \geq A_n w_2$, and define $\varphi := w_1 - w_2$.
Let $j^\ast\in V$ be such that $\varphi_{j^\ast} = \min_j \varphi_j$, and suppose for contradiction that $\min_j \varphi_j < 0$.
Recall that the linear operator associated to $A_n$ is of the form \eqref{jacobian general}, and note that $(A_n \varphi)_j\geq 0$ for all $j\in \cJ \cup \cB$.

For any $j\in \cB$, we have $\varphi_j = \frac{1}{\lambda_b}  [A_n (w_1 - w_2)]_j \geq 0$, hence the minimum of $\varphi$ cannot be attained on the boundary $\cB$, i.e. $j^\ast\in \cJ$. Using that $\beta_j \geq 0$, $\alpha_{jk}>0$, $\varphi_{j^\ast} <0$, and that $\varphi_{j^\ast} - \varphi_k\leq 0$ for all $k\in \cV_{j^\ast}$, we obtain
$$
0 \leq \left( A_n \varphi \right)_{j^\ast} = \beta_{j^\ast} \varphi_{j^\ast} + \sum_{k\in \cV_{j^\ast}} \alpha_{j^\ast k} (\varphi_{j^\ast} - \varphi_k) \leq  0.
$$
  Now, since $\alpha_{jk}>0$ for all $k\in \cV_j$, the above inequalities imply that $\varphi_k = \varphi_{j^\ast}$ for all $k\in \cV_j$, i.e., the minimum of $\varphi$ is attained also at all the neighbors of $j^\ast$. We can apply the same argument to all the neighbors of $j^\ast$. By a recursive argument, due to the connectivity property {\bf (H4)}$(ii)$, the minimum of $\varphi$ is also attained at a boundary point $j\in \cB$, which contradicts the fact that  $\min_j \varphi_j <0$.

The maximum principle implies that $A_n$ is invertible and that its inverse $A^{-1}_n$ has all the elements greater or equal than $0$.\\

\textit{Step 2: Use of the barrier function.}
By Lemma \ref{lem:barrier}, there exists a grid function $\psi\in \R^V$ such that $A_n \psi \geq \bf{1}$ for all $n$.
By the maximum principle and the fact that all the elements of $A_n^{-1}$ are non-negative, we have
$$
\psi \geq A_n^{-1} \mathbf{1} \geq 0, \quad \text{which implies} \quad \| \psi\|_{\ell^\infty (V)} \geq \| A_n^{-1} \mathbf{1}\|_{\ell^\infty (V)}.
$$
Then, since we deal with the supremum norm, using again that $A_n^{-1}$ has all the entries greater than or equal to $0$, we obtain 
\begin{eqnarray*}
\| A_n^{-1}\|_{\ell^\infty (V)} & = & \max_{j} \sum_k \left|(A_n^{-1})_{jk} \right|
   =  \| A_n^{-1} \mathbf{1}\|_{\ell^\infty (V)} \leq \| \psi \|_{\ell^\infty (V)},
\end{eqnarray*}
which concludes the proof using the upper bound estimate for $\psi$ from Lemma~\ref{lem:barrier}.
\end{proof}

We are now in a position to prove the main result of this section.
Here, $\partial L(u)$ denotes the generalized differential of $L(\cdot)$ at $u$, and is defined analogously to the differential of $\cR(u)$ in Definition~\ref{def: generalized differential}.

\begin{theorem}[Well-posedness of residual minimization] 
\label{thm:exist-uniqueness-scheme-gen}
  Assume (H1), (H2) and ((H3) or (H4)). Then,
  \begin{enumerate}[label=(\roman*)]
  \item the solution of \eqref{eq:scheme-steady-1} exists and is unique.
  \item the solution of \eqref{eq:scheme-steady-1} is the unique critical point of the functional 
  $L(u)$ defined in \eqref{Lp functional steady state eq}:
  $$
     0\in \partial L(u) 
     \qquad \text{if and only if} \qquad
     \cR (u) = 0.
  $$
  \end{enumerate}
\end{theorem}
\begin{proof}
$(i)$
  Using Lemma~\ref{lem:DR-invertible lambda} when {\bf(H3)} holds, 
  and Lemma~\ref{lem:DR-invertible alpha} when {\bf(H4)} holds, 
  we obtain that, for any $u\in \R^{V}$ and any matrix $A$ in the generalized differential $\partial \cR (u)\subset \R^{V \times V}$ (see Definition~\ref{def: generalized differential}), the inverse of $A$ exists and satisfies $\|A^{-1}\|_{\ell^\infty (V)} \leq C$, with a constant $C>0$ independent of $u$, and given explicitly in Lemmas~\ref{lem:DR-invertible lambda} and \ref{lem:DR-invertible alpha} respectively. By a generalization of Hadamard's theorem by Pourciau \cite{pourciau-82}, we deduce that the map $u\conv \cR(u)$ is a homeomorphism of $\R^{V}$ onto itself. 
   In particular, there is a unique $u\in \R^V$ such that $\cR(u)=0$.

$(ii)$
    If $\cR(\bar u)=0$, then it is clear that $\bar u$ is a global minimum of $L(\cdot)$  and therefore
    $0 \in \partial L(\bar u)$. 
    We now assume that we have a critical point, in the sense that $0\in \partial L(\bar u)$, and aim to prove
    $\cR(\bar u)=0$.
    Without loss of generality, we can assume here that $\ml_b=1$ and that 
    $$
      L(u) := \|\cR(u)\|_{\ell^2(V)}^2.
    $$
    
    If we assume that $\cR(\cdot)$ is differentiable at $u\in \R^V$, then $\partial L(u)$ contains only the classical gradient of $L(\cdot)$ at $u$, which is given by
    $$
    \nabla L(u) = 2\sum_{j\in V} \nabla [\cR(u)_j] \cR(u)_j = 2 D\cR(u^\ast)^\top  \cR(u),
    $$
    where $(D \cR(u))^\top$ denotes the matrix transpose of $D \cR(u)$.
    Note that $\nabla [\cR(u)_j]$ is the $j$-th row of the Jacobian matrix $D \cR(u)$.
    
    In general, in view of Definition~\ref{def: generalized differential}, for any $u\in \R^V$, any vector in the generalized differential $w\in \partial L(u)$ can be written as
    $$
    w = 2A^\top \cR(u), \qquad \text{for some} \quad A\in \partial\cR(u).
    $$
    Since $A$ is invertible, and therefore so is $A^\top$, we obtain $\cR(u) = \frac{1}{2} (A^\top)^{-1} w$. Hence, if $0\in \partial L(\bar u)$, then $\cR(\bar u) = \frac{1}{2} (A^\top)^{-1} 0 = 0$.
\end{proof}

\subsection{Alternative proof using a fixed point argument}

We give an alternative existence and uniqueness proof via the Euler map $u- \rho\cR(u)$.
Under {\bf(H3)}, this was established in Oberman~\cite{oberman2006convergent}; here we extend the argument to {\bf(H4)} using a weighted norm constructed from the barrier of Lemma~\ref{lem:barrier}. Beyond re-proving well-posedness, the contraction property yields a convergent fixed-point iteration $u_{n+1}=u_n-\rho\,\cR(u_n)$ for the discrete solution.

\begin{theorem}
\label{thm: euler map is a contraction}
Let $\cR(\cdot)$ be a residual function of the form \eqref{residual function def}, with $\cH$ satisfying {\bf (H1)}--{\bf (H2)}, and {\bf (H3)} or {\bf (H4)}, then there exists a norm in $\R^V$ such that the Euler map  $\Phi:\R^V\conv \R^V$ defined by
$$
  \Phi(u):= u- \rho \cR(u)
$$ 
is a contraction for $\rho>0$ small enough.
As a consequence, the equation $\cR(u) = 0$ admits a unique solution in $\R^V$.
\end{theorem}

\begin{proof}
The case of (H3) was already proved in Oberman \cite{oberman2006convergent}. 

Let us consider the case when only (H4) holds.
By Lemma~\ref{lem:barrier}, there exists a barrier $w \in \R^V$ such that, for all $u\in \R^V$ and $A\in \partial \cR(u)$, we have $A w\geq \mathbf{1}$. Since, by virtue of Lemma~\ref{lem:DR-invertible alpha}, $A$ is invertible and all the elements of $A^{-1}$ are non-negative, we have $w \geq A^{-1}\mathbf{1}$. Additionally, $w_j>0$ for all $j\in V$.

Let us define the norm
$$
\|x\|_w:= \max_{j\in V} \frac{|x_j|}{w_j}, \qquad \text{for any vector} \  x\in \R^{V}.
$$
The associated weighted matrix norm satisfies
$$
\|M\|_w:=\max_{x\neq 0} \frac{\|M x\|_w}{ \|x\|_w} \leq  \max_{j} \left( \frac{1}{w_j}\sum_i |M_{ji}| w_{i} \right), \qquad \text{for any matrix}\ M\in \R^{V\times V}.
$$
We aim to find a $\rho>0$ such that, for any $u\in \R^V$, and any $A\in \partial \cR(u)$, 
$\| I - \rho A\|_w\leq  \kappa <1$ for some $\kappa>0$ independent of $u$.

Consider the case when  $\cR (\cdot)$ is differentiable at $u$.
Notice that, if we take $\rho>0$ such that $\rho \max_j D\cR(u)_{jj}\leq 1$, we have
    \beno
      \sum_k |(I-\rho D\cR(u))_{jk}| w_k
	& = & (1- \rho D\cR(u)_{jj}) w_j + \rho \sum_{k\neq j} |D\cR(u)_{jk}| w_k \\
	& = & w_j - \rho  \left(D\cR(u)_{jj}w_j -  \sum_{k\neq j} |D\cR(u)_{jk}| w_k \right) \\ 
    & = & w_j - \rho  \sum_k D\cR(u)_{jk} w_k = w_j -\rho (D\cR(u) w)_j \\
	& \leq & w_j -\rho.
     \eeno
    Here we used that the off-diagonal terms of the Jacobian matrix, $D\cR(u)_{jk}$ with $j\neq k$, are all non-positive due to \eqref{eq:dH(uj)/duk} and the monotonicity condition {\bf (H2)}; and that $(D\cR(u)w)_j\geq 1$ for all $j$.
  Hence,
   \beno
     \|I-\rho D\cR(u)\|_w 
       &    \leq & \max_j \frac{1}{w_j} \sum_k | (I-\rho D\cR(u))_{jk} | w_k \\
       & \leq & \max_j \frac{w_j -\rho}{w_j} = 1  - \frac{\rho}{\max_j w_j} = 1 - \frac{\rho}{\|w\|_{\ell^\infty(V)}}.
   \eeno
   We only need to select $\rho>0$ such that $\rho \max_j D\cR(u)_{jj}\leq 1$ for all $u$.
   Let $Q=\max_u \max_j D\cR(u)_{jj}$ (by the assumptions on $\cR$, we have $Q\leq \max\!\big(\lambda_b,\; \KH \, (1 + N/\dxmin)\big) < \infty$, where $N=\max_j |\cV_j|$).
   Hence, we choose $\rho$ such that $0< \rho < \frac{1}{Q}$,
   which gives $\|I-\rho D\cR(u)\|_w \leq \kappa$ for some $\kappa = 1-\frac{\rho}{\|w\|_{\ell^\infty(V)}}$.
   Therefore $u- \rho \cR(u)$ is a contraction for the norm $\|\cdot \|_w$.

   When $\cR(\cdot)$ is not differentiable, the proof can be adapted by using the generalized differential.
   Let $u,v$ be given in $\R^V$. By Lemma~\ref{lem:mvt}, $\cR(u)-\cR(v)=A(u-v)$ for some $A\in co(\partial\cR([u,v]))$. Every element of $\partial\cR([u,v])$ is a limit of Jacobians $A_n$ with $\|I-\rho A_n\|_w \leq \kappa$. Since $\{M:\|I-\rho M\|_w\leq \kappa\}$ is closed and convex, the same bound $\|I-\rho A\|_w \leq \kappa$ holds for $A$.
   Therefore $\| u-\rho \cR(u) - (v-\rho\cR(v)) \|_w= \| (I- \rho A)(u-v)\|_w\leq \mk \|u-v\|_w$, 
   which concludes the proof.
\end{proof}

\section{Error estimates and Polyak--\L{}ojasiewicz inequalities}
\label{sec: stability}

\subsection{Proper discretizations}
We prove stability estimates for the solution of \eqref{eq:scheme-steady-1} in terms of the residual $\cR(u)$.
These serve as {\it a posteriori} error estimates, since the right-hand sides of \eqref{stability estimate Linfty-H3} and \eqref{stability estimate Linfty-H4} are computable (e.g., via stochastic sampling on the grid).

For the case when $\cH$ is proper (i.e., {\bf (H3)} holds), we have the following error estimate.
\begin{theorem}\label{thm:stability-H3}
  Let $\cR(\cdot)$ be a residual function of the form \eqref{residual function def}, with $\cH$ satisfying {\bf (H1)}, {\bf (H2)} and {\bf (H3)}, and let $\bar u$ be the unique solution to $\cR(\bar u)=0$.
Then, for any $u\in \R^V$, we have 
   \begin{equation}
     \label{stability estimate Linfty-H3}
     \|u - \bar u\|_{\ell^\infty(V)} \leq \max\bigg(
     \frac{1}{\ml} \|\cR(u)\|_{\ell^\infty(\cJ)},\
     \| u - g \|_{\ell^\infty(\cB)}
     \bigg).   
   \end{equation}
\end{theorem}
\begin{proof} 
The proof is obtained by classical maximum principle arguments and is left to the reader.
\end{proof}

As a consequence of Lemma \ref{lem:DR-invertible lambda}, and the equivalence inequality $\| A\|_{2} \leq \sqrt{N}\,  \| A\|_\infty$ for all matrices $A\in \R^{N\times N}$, we can obtain the following $\ell^2$-error estimate, which also shows that the loss function $L(u) := \| \cR(u)\|_{\ell^2(V)}^2$ is coercive in $\R^V$.

\begin{corol}
\label{cor: L is coercive}
Let $\cR(\cdot)$ be a residual function of the form \eqref{residual function def}, with $\cH$ satisfying {\bf (H1)}, {\bf (H2)} and {\bf (H3)}, and let $\bar u$ be the unique solution to $\cR(\bar u)=0$.
Then, for any $u\in \R^V$, we have 
$$
\| u-\bar{u}\|_{\ell^2(V)}^2 \leq \frac{|V|}{\mu^2} \| \cR (u)\|_{\ell^2(V)}^2 = \frac{|V|}{\mu^2} L(u),
$$
where $\mu = \min (\lambda, \lambda_b)$.
\end{corol}

\begin{proof}
    By Lemma~\ref{lem:mvt} and $\cR(\bar u)=0$, $\cR(u)=A(u-\bar u)$ for some $A\in co(\partial\cR([u,\bar u]))$.

Due to Lemma \ref{lem:DR-invertible lambda}, we have $\| A^{-1} \|_{\ell^\infty(V)} \leq \frac{1}{\mu}$, and using the standard equivalency between $\ell^2$ and $\ell^\infty$ norms, we obtain $\| A^{-1} \|_{\ell^2(V)} \leq \frac{\sqrt{|V|}}{\mu}$, with $\mu = \min (\lambda, \lambda_b)$.

Then, we conclude the proof by left-multiplying the identity $\cR(u) = A  (u - \bar u)$ by $A^{-1}$, to obtain
$$
\| u-\bar u\|_{\ell^2(V)}^2 = \| A^{-1} \cR (u) \|_{\ell^2(V)}^2 \leq \frac{|V|}{\mu^2} \| \cR(u)\|_{\ell^2(V)}^2.
$$
\end{proof}

The drawback of the $\ell^2$-error estimate from Corollary~\ref{cor: L is coercive} is that the coercivity constant depends on the number of nodes $|V|$: on a uniform Cartesian grid in $\R^d$ it scales as $O(\dx^{-d})$, which is especially problematic on the fine grids needed in high dimensions. The loss $L(u)$ can be very small while $u$ is still far from $\bar u$. We can say that, without further assumptions, the $\ell^2$-error estimate suffers from the curse of dimensionality.

In the following, we prove that if $\lambda \dx$ and $\lambda_b \dx$ are large enough, then the constant for the $\ell^2$ estimate depends on local stencils, but does not depend on the total number of grid points. In other words, if the damping $\lambda$ is sufficiently large relative to the grid spacing, the error is controlled by the residual in $\ell^2$, and the associated constant does not suffer from the curse of dimensionality. This is established in Theorem~\ref{thm:dim-robust error estimate proper case} below.

This result is a consequence of the following elementary lemma in linear algebra.

\begin{lem}
\label{lem:L1Linfty-interpolation}
For any invertible matrix $A\in\R^{N\times N}$ with nonnegative inverse $A^{-1}\geq 0$ (entrywise), if $A \mathbf{1} \geq c_1 \mathbf{1}$ and $A^\top\mathbf{1}\geq c_2\mathbf{1}$ for some $c_1, c_2 >0$, then 
$$
\|A^{-1}\|_\infty \leq \frac{1}{c_1} \quad \text{and} \quad \|A^{-1}\|_1 \leq \frac{1}{c_2}.$$ 
Hence, we have $\|A^{-1}\|_2\le\sqrt{\|A^{-1}\|_1\,\|A^{-1}\|_\infty}\leq \sqrt{\frac{1}{c_1 c_2}}$.
\end{lem}

\begin{proof}
Since $A$ is invertible with $A^{-1}\geq 0$, the inequality $A\mathbf{1} \geq c_1 \mathbf{1}$ implies $A^{-1}\mathbf{1} \leq \frac{1}{c_1}\mathbf{1}$, which in turn implies $\| A^{-1}\mathbf{1} \|_\infty \leq \frac{1}{c_1}$. Hence, using $A^{-1}\geq 0$ to write each absolute row sum as a row sum,
$$
\| A^{-1}\|_\infty = \max_{k = 1, \ldots, N} \sum_{j=1}^N (A^{-1})_{kj} = \| A^{-1} \mathbf{1}\|_\infty \leq \frac{1}{c_1}.
$$
Similarly, we obtain $\| (A^\top)^{-1}\|_\infty \leq \frac{1}{c_2}$, and therefore we have
$$
\| A^{-1}\|_1 = \| (A^{-1})^\top\|_\infty = \| (A^\top)^{-1} \|_\infty \leq \frac{1}{c_2}.
$$
And the conclusion follows.
\end{proof}

We can now prove the following dimension-robust $\ell^2$-error estimate.

\begin{theorem}
\label{thm:dim-robust error estimate proper case}
  Let $\cR(\cdot)$ be a residual function of the form
  \eqref{residual function def}, with $\cH$ satisfying {\bf (H1)}, {\bf (H2)}
  and {\bf (H3)}, and let
  \begin{equation}
  \label{def V'j}
    \cV_j' := \{\, i \in \cJ \ :\ j \in \cV_i \,\}
  \end{equation}
  be the set of interior nodes whose stencil contains $x_j$. Suppose that
  \begin{equation}\label{eq:dimrobust-conditions}
    \lambda > \frac{\KH_p}{\dx_{min}}\,  \max_{j\in\cJ} \left( \max(|\cV_j| , \, |\cV_j'|)\right)
    \qquad\text{and}\qquad
    \lambda_b > \frac{\KH_p}{\dx_{min}}\,\max_{j\in\cB}|\cV_j'|,
  \end{equation}
  and set
  \[
    c_1 := \min(\lambda,\lambda_b), \qquad
    c_2 := \min\!\Big(
      \lambda - \tfrac{\KH_p}{\dx_{min}}\max_{j\in\cJ}|\cV_j'|,\;
      \lambda_b - \tfrac{\KH_p}{\dx_{min}}\max_{j\in\cB}|\cV_j'|
    \Big) > 0 .
  \]
  Let $\bar u$ be the unique solution of $\cR(\bar u)=0$. Then the loss
  $L(u):=\|\cR(u)\|_{\ell^2(V)}^2$ satisfies
  \[
    \| u-\bar u\|_{\ell^2(V)}^2 \leq \frac{1}{c_1\,c_2}\, L(u),
    \qquad \forall u\in\R^V,
  \]
  together with the generalized Polyak--\L{}ojasiewicz inequality
  \[
    \| w \|_{\ell^2(V)}^2 \geq 4\,c_1\,c_2\, L(u),
    \qquad \forall u\in\R^V,\ \forall w\in\partial L(u),
  \]
  where $\partial L(u)$ denotes the generalized gradient of $L$ in the sense of Clarke.

  Moreover, for any $u\in \R^V$ and $A\in \partial \cR(u)$, the condition number of $A$ can be estimated as
  \begin{equation}
  \label{eq:dim-robust condition number proper case}
  \kappa (A) := \| A\|_{\ell^2(V)} \| A^{-1}\|_{\ell^2(V)} \leq (c_1 c_2)^{-1/2} \max \left\{ \KH_u + 2\lambda, \, 2\lambda_b \right\}.
  \end{equation}
\end{theorem}

\begin{proof}
The argument is that of Corollary~\ref{cor: L is coercive}, using
Lemma~\ref{lem:L1Linfty-interpolation}. For any $A\in\partial\cR(u)$,
diagonal dominance gives $A\mathbf 1 \geq c_1\mathbf 1$ with
$c_1=\min(\lambda,\lambda_b)$. For $A^\top\mathbf 1$ we use
\eqref{eq:dH(uj)/duj}--\eqref{eq:dH(uj)/duk} separately on interior and
boundary rows: for $j\in\cJ$,
\begin{eqnarray*}
  (A^\top\mathbf 1)_j & = & \frac{\partial\cH}{\partial u} (x_j, u_j, \nabla_G u_j) + \sum_{k\in \cV_j} \frac{1}{\dx_{jk}} \frac{\partial \cH}{\partial p_k} (x_j, u_j, \nabla_G u_j) - \sum_{i\in \cV_j'} \frac{1}{\dx_{ij}} \frac{\partial \cH}{\partial p_j} (x_i, u_i, \nabla_G u_i)  \\
  &\geq & \lambda - \tfrac{\KH_p}{\dx_{min}}\max_{i\in\cJ}|\cV_i'|,
\end{eqnarray*}
while for $j\in\cB$ the column-$j$ off-diagonal entries come only from the
$|\cV_j'|$ interior neighbours of $j$, so
\[
  (A^\top\mathbf 1)_j = \lambda_b
    - \sum_{i\in\cV_j'}\tfrac{1}{\dx_{ij}}\tfrac{\partial\cH}{\partial p_j}(x_i,u_i,\hD u_i)
  \geq \lambda_b - \tfrac{\KH_p}{\dx_{min}}\max_{i\in\cB}|\cV_i'|.
\]
Hence $A^\top\mathbf 1\geq c_2\mathbf 1$. Since $A$ has positive diagonal and nonpositive off-diagonal entries and is diagonally dominant (by \eqref{eq:dH(uj)/duj}--\eqref{eq:dH(uj)/duk} and {\bf (H2)}--{\bf (H3)}), it is an M-matrix, so $A^{-1}\geq 0$ (entrywise), and
Lemma~\ref{lem:L1Linfty-interpolation} yields $\|A^{-1}\|_{\ell^2(V)}\leq (c_1 c_2)^{-1/2}$.
The error estimate then follows as in Corollary~\ref{cor: L is coercive}. For the Polyak--\L{}ojasiewicz inequality, every $w\in\partial L(u)$ has the form $w = 2A^\top\cR(u)$ with $A\in\partial\cR(u)$, so $\|w\|_{\ell^2(V)} \geq 2\,\sigma_{\min}(A)\,\|\cR(u)\|_{\ell^2(V)} \geq 2\sqrt{c_1 c_2}\,\|\cR(u)\|_{\ell^2(V)}$, whence $\|w\|_{\ell^2(V)}^2 \geq 4 c_1 c_2\, L(u)$.

To estimate the condition number of $A$ as in \eqref{eq:dim-robust condition number proper case}, it only remains to upper bound $\| A\|_{\ell^2(V)}$. Using $\| A\|_{\ell^2(V)} \leq \sqrt{\| A\|_{\ell^\infty(V)} \| A\|_{\ell^1(V)}}$, it amounts to control the sums of the absolute values of the elements of $A$, row- and column-wise. Using \eqref{eq:dH(uj)/duj}--\eqref{eq:dH(uj)/duk} and \eqref{eq:dimrobust-conditions}, the row sums can be controlled by
$$
\sum_{k\in V} |A_{jk}| \leq |A_{jj}| + \sum_{k\in \cV_j} |A_{jk}| \leq \KH_u + 2 \frac{\KH_p}{ \dx_{\min} } \max_{j\in \cJ} |\cV_j| \leq \KH_u + 2\lambda, 
$$
for all $j\in \cJ$, and $\sum_{k\in V} |A_{jk}| = \lambda_b$ for all $j\in \cB$. Likewise, the column sums can be controlled by
$$
\sum_{k\in V} |A_{kj}| \leq |A_{jj}| + \sum_{k\in \cV_j'} |A_{kj}| \leq \KH_u + \frac{\KH_p}{ \dx_{\min} } \left( \max_{j\in \cJ}|\cV_j| + \max_{j\in \cJ}|\cV_j'|\right) \leq \KH_u + 2\lambda,
$$
for all $j\in \cJ$, and by
$$
\sum_{k\in V} |A_{kj}| \leq |A_{jj}| + \sum_{k\in \cV_j'} |A_{kj}| \leq \lambda_b + \frac{\KH_p}{ \dx_{\min} } \max_{j\in \cB} |\cV_j'| \leq 2 \lambda_b, \qquad \text{for all $j\in \cB$.}
$$
\end{proof}

\paragraph{Stencil cardinalities and dimension-robustness.}
Recall the two stencils at a node $x_j$. The out-stencil $\cV_j$ holds the
neighbours used to form $\nabla_G u_j$. The in-stencil
$\cV_j'=\{\,i\in\cJ : j\in\cV_i\,\}$ holds the interior nodes whose scheme
involves $u_j$. Take a uniform Cartesian grid of a $d$-dimensional box
$\Omega=\prod_{k=1}^d(a_k,b_k)$ with constant step $\dx$, whose faces lie on
grid hyperplanes. An interior node has all $2d$ axis-neighbours, so
$\max_{j\in\cJ}|\cV_j'|=2d$. A boundary node has (at most) a single interior
axis-neighbour, in the inward normal direction. Its $2(d-1)$ tangential
neighbours lie on the same face, and its outward neighbour lies outside
$\Omega$. Hence $\max_{j\in\cB}|\cV_j'|=1$, with edges and corners giving
fewer. Both counts are exact, not merely generic.

Conditions~\eqref{eq:dimrobust-conditions} then read $\lambda > 2d\,\KH_p/\dx$ and $\lambda_b > \KH_p/\dx$, so the dimension enters only through the factor $2d$ in the condition on $\lambda$, while the one on $\lambda_b$ is dimension-independent. Consequently $c_1$ and $c_2$---and hence the error and PL constants---depend on the dimension only through the local stencil size $2d$, never through the global node count $|V|=O(\dx^{-d})$.

A boundary not aligned with the grid is resolved only up to a staircase. Then
$|\cV_j'|$ depends on the local geometry, rather than being uniformly $1$.
Such domains fall under the general framework of
Section~\ref{sec: main results stationary case}, where $\cV_j$ and the steps
$\dx_{jk}$ are arbitrary (for instance, stencils truncated at the boundary).
The well-posedness and error estimates still hold. One only replaces
$|\cV_j|$ and $|\cV_j'|$ by the maximal out- and in-degrees of the graph.
These degrees remain $O(d)$ for any local stencil, so dimension-robustness
persists. Only the clean value $|\cV_j'|=1$ is specific to the box. In the
neural-network implementation, the boundary condition is imposed by penalizing
sampled points. Curved domains such as the spherical shell of
Section~\ref{sec:isaacs} therefore need no grid-conforming boundary.

\subsection{Elliptic discretizations}

If {\bf (H4)} holds, then we have the following error estimate.

\begin{theorem}
\label{thm:stability-H4}
Let $\cR(\cdot)$ be a residual function of the form \eqref{residual function def}, with $\cH$ satisfying {\bf (H1)},{\bf (H2)} and {\bf (H4)}.
Let $\bar u$ be such that $\cR(\bar u)=0$ (i.e. solution of the scheme).
Then, for any $u\in \R^V$, we have
\begin{equation}
\label{stability estimate Linfty-H4}
\|u - \bar u\|_{\ell^\infty(V)} \leq C \|\cR(u)\|_{\ell^\infty(V)},
\end{equation}
where $C$ is the explicit constant of Lemma~\ref{lem:DR-invertible alpha} (see also Definition~\ref{def: graph distance constants} for $\overs, \subalpha, \dxmax, \dxmin$).
\end{theorem}

\begin{proof}
By Lemma~\ref{lem:mvt} and $\cR(\bar u)=0$, we have that $\cR(u)=A(u-\bar u)$ for some $A\in co(\partial\cR([u,\bar u]))$.
By the uniform estimate of Lemma~\ref{lem:DR-invertible alpha}, $\|A^{-1}\|_{\ell^\infty(V)} \leq C$ for every $A\in \partial\cR(v)$, and $v\in [u, \bar u]$. Any convex combination of such matrices is again of the form \eqref{jacobian general} with coefficients obeying \eqref{estimates beta alpha lemma}, so the same bound holds for our $A$. 

Hence, we have
$$
   \| u- \bar u\|_{\ell^\infty(V)} = \| A^{-1} A(u - \bar u)\|_{\ell^\infty(V)} \leq C \|A(u - \bar u)\|_{\ell^\infty(V)} = C \| \cR (u) \|_{\ell^\infty(V)}.
$$
\end{proof}

\paragraph{Discretizations on Cartesian grid.}
Let us consider that the graph $G$ is a uniform Cartesian grid of a $d$-dimensional bounded domain $\Omega\subset \R^d$. This means that the index set $V$ is a bounded subset of the $d$-dimensional lattice $\Z^d$, i.e. 
\begin{equation}
\label{Cartesian grid def}
\{ x_j\}_{j\in V} = \{ j\dx \ : \ j\in V\}, \quad \text{for some} \ V\subset \Z^d \ \text{bounded,}
\end{equation}
and that any interior point $x_j$ with $j\in \cJ$ has exactly $2d$ neighbors as follows:
\begin{equation}
\label{Cartesian neighbors}
\cV_j = \{ j \pm e_k \ : \quad k=1, \ldots, d \},
\end{equation}
where $\{ e_k\}_{k=1}^d$ is the canonical orthonormal basis of $\R^d$. We note that the distance between neighboring points is constant, i.e. $\dx_{jk} = \dx$ for all $k\in \cV_j$.

In this setting, the constant $\overs$ in Definition~\ref{def: graph distance constants} can be estimated as
\begin{equation}
  \label{over s Cartesian grid}
  \overs = \esssup_{j,u,p} \left(  \sum_{k\in \cV_j}
    \frac{\partial\cH}{\partial p_k} (x_j, u_j, p)  \right)
     \leq \KH_p \max_{j\in \cJ} | \cV_j| \leq 2d \KH_p, 
\end{equation}
(where $\KH_p$ is the Lipschitz constant of $\cH$ with respect to the $\nabla_G u$ parameters),
and the maximum of the graph distance function $d_\infty$ can be estimated by
\begin{equation}
\label{in-radius Cartesian grid}
d_\infty \leq \frac{R}{\dx}, \qquad \text{where} \ 
R:= \max_{x\in \Omega} \min_{y\in \partial \Omega} \| x - y\|_\infty.
\end{equation}
Hence, the constant $C$ in the error estimate from Theorem \ref{thm:stability-H4} scales as $O(\exp(\dx^{-1}))$. Next we show that, when the finite-difference residual is consistent with a semilinear elliptic PDE, then one can obtain an $\ell^\infty$-error bound which is independent of $\dx$.

\paragraph{Viscous Hamilton-Jacobi equation.}
We shall now prove that, when the scheme is consistent with a semilinear second-order elliptic PDE,
$$
- \sum_{\ell=1}^d \nu_\ell u_{x_{\ell \ell}} + H(x, u, \nabla u) = 0, \quad x\in \Omega \qquad \text{with} \qquad u(x) = g(x), \quad x\in \partial\Omega,
$$
then the constant in the error bound is independent of  $\dx$. Moreover, if the sum of the diffusion coefficients $\sum_{\ell=1}^d \nu_\ell$ and $\ml_b$ are sufficiently large, 
we obtain a dimension-robust $\ell^2$-estimate as well as the associated Polyak--\L{}ojasiewicz inequality.

Let us consider the following residual function on a Cartesian grid $\{ x_j\}_{j\in V}$ of $\overline{\Omega}\subset\R^d$:
\begin{equation}
\label{residual viscous HJ eq}
\cR(u)_j :=
\begin{cases}
 \displaystyle\sum_{\ell=1}^d \nu_\ell \dfrac{2u_j - u_{j-e_\ell}-u_{j+e_\ell}}{\dx^2}  + \cH (x_j, u_j, \nabla_G u_j), & \text{if} \ j\in \cJ, \\
 \lambda_b (u_j- g_j), & \text{if} \ j\in \cB,
\end{cases}
\end{equation}
with $\nu_\ell > 0$ for all $\ell=1,\dots,d$, and a numerical Hamiltonian $\cH$ satisfying {\bf (H1)} and {\bf (H2)}.
Notice that, on a Cartesian grid, the numerical gradient is given by 
$$\nabla_G u_j := (D_\ell^\pm u_j)_{\ell=1}^d, \quad \text{where} \quad D^\pm_\ell u_j = \frac{u_j-u_{j\pm e_\ell}}{\dx} \quad  \text{for all}  \quad  \ell=1,\dots,d,$$
and the discrete Laplacian term is a linear combination of $\nabla_G u_j$:
\begin{equation}
\label{disc-lap as diff gradients}
  \sum_{\ell=1}^d \nu_\ell \frac{2u_j - u_{j-e_\ell}-u_{j+e_\ell}}{\dx^2} 
   = \sum_{\ell = 1}^d  \frac{\nu_{\ell}}{\dx} \left( D^-_\ell u_j + D^+_\ell u_j \right).
\end{equation}

We have the following result.

\begin{theorem}
\label{thm:  viscous HJ}
Let $\Omega \subset \R^d$ be a bounded domain with diameter $R:= \max_{x,y\in \Omega} \|x-y\|_\infty$.
We consider a uniform Cartesian grid of $\overline{\Omega}$ with discretization step $\dx>0$.
Let $\cR(u)$ be a residual function of the form \eqref{residual viscous HJ eq}, 
with $\cH$ satisfying {\bf (H1)} and {\bf (H2)}. Then,
\begin{enumerate}[label=(\roman*)]
    \item
    There exists a unique $\bar u\in \R^V$ such that $\cR(\bar u) = 0$, and this solution is the unique critical point of the functional $L(u):= \| \cR (u)\|_{\ell^2(V)}^2$. 
    \item Assume there exists $\ell^\ast\in \{ 1, \ldots, d\}$ such that
    \be \label{eq:dx-condition-small}
      \dx \leq \frac{\sqrt{3}\nu_{\ell^\ast}}{2\KH_p} e^{-\frac{\KH_p}{\nu_{\ell^\ast}} R}.
    \ee
    Then
    \be \label{eq:u-ubar-estimate}
    \| u - \bar u \|_{\ell^\infty(V)} \leq   
         \left( \frac{\nu_{\ell^\ast}}{2 \KH_p^2}  e^{\frac{2 \KH_p}{\nu_{\ell^\ast}}R} + \dfrac{1}{\lambda_b} \right) \| \cR(u)\|_{\ell^\infty(V)},\qquad \text{for any $u\in \R^V$.}
    \ee
    \item 
    Furthermore, if $\dx$ and $\ml_b$ are large enough, namely,
    \begin{equation}
    \label{eq:l2 cond visc HJ}
    \dx \sum_{\ell=1}^d\nu_\ell  > \KH_p R^2 d\quad
    \text{and} \quad \lambda_b \geq \frac{2d\,(\nu_{\max}+\KH_p\dx)}{\dx^2}, 
    \end{equation}
    where $\nu_{\max}:=\max\limits_{\ell=1, \ldots, d} \nu_\ell$,
    then we have the $\ell^2$-error estimate
    $$
    \| u - \bar u\|_{\ell^2 (V)}^2 \leq \left( \frac{2}{m^2} + \frac{1}{\lambda_b^2} \right) L(u), \qquad \text{for all}\  u\in \R^V,
    $$
    together with the generalized Polyak--\L{}ojasiewicz inequality
    \[
    \| w \|_{\ell^2(V)}^2 \geq \left( \frac{1}{2m^2} + \frac{1}{4\lambda_b^2} \right)^{-1} \, L(u),
    \qquad \forall u\in\R^V,\ \forall w\in\partial L(u),
    \]
  with $m := \frac{4}{R^2}\sum_{\ell=1}^d\nu_\ell - \frac{4 \KH_p\, d}{\dx}$,
  and  where $\partial L(u)$ denotes the generalized gradient of $L$ in the sense of Clarke.
  Moreover, for any $u\in \R^V$ and $A\in \partial \cR(u)$, the condition number of $A$ can be estimated as
  \begin{equation*}
  \kappa (A) := \| A\|_{\ell^2(V)} \| A^{-1}\|_{\ell^2(V)} \leq \left( \frac{2}{m^2} + \frac{1}{\lambda_b^2} \right)^{1/2} \max \left\{ 4\left( \frac{1}{\dx^2} + \frac{1}{R^2} \right)\sum_{\ell=1}^d \nu_\ell + \KH_u, \, 2\lambda_b \right\}.
  \end{equation*}
\end{enumerate}
\end{theorem}

\begin{proof}
   $(i)$ In view of the definition of the residual function $\cR(u)$ in \eqref{residual viscous HJ eq} and the identity \eqref{disc-lap as diff gradients}, the residual at the interior points can be written as
    $$
    \cR(u)_j = \cH' (x_j, u_j, D^-u_j, D^+ u_j) :=
    \sum_{\ell=1}^d \frac{\nu_\ell}{\dx} \left( D^-_\ell u_j + D^+_\ell u_j \right) + \cH (x_j, u_j, \nabla_G u_j), \quad \forall j\in \cJ.
    $$
    which satisfies {\bf (H1)}, {\bf (H2)} and {\bf (H4)} with $\subalpha = \min_\ell \frac{\nu_\ell}{\dx}$. Note that we assumed that $\cH$ satisfies {\bf (H1)}, {\bf (H2)}, so it is already non-decreasing with respect to $D^- u_j$ and $D^+ u_j$. We can then apply directly Theorem \ref{thm:exist-uniqueness-scheme-gen}.
    
    $(ii)$ The proof follows the same idea as the proof of Theorem \ref{thm:stability-H4}, where now we prove the estimate
    \be\label{eq:A-inv-main-estimate} 
    \| A^{-1}\|_{\ell^\infty(V)} \leq   
         \frac{\nu_{\ell^\ast}}{2 \KH_p^2} \left( e^{\frac{2 \KH_p}{\nu_{\ell^\ast}}R} - 1 \right) - \frac{R}{\KH_p} + \dfrac{1}{\lambda_b}, \qquad \forall u\in \R^V \ \text{and} \ A\in \partial \cR(u),
    \ee 
    which we can simplify (by dropping the negative terms) to obtain the rough bound \eqref{eq:u-ubar-estimate}.
    As in Lemma~\ref{lem:DR-invertible alpha}, we shall construct a barrier function $\psi\in \R^V$ such that for all $u\in \R^V$
    and $A \in \partial\cR(u)$, it holds that $(A \psi)_j \geq 1$ for all $j\in V$. Then, $\|A^{-1}\|_{\ell^\infty(V)}\leq \| \psi\|_{\ell^\infty (V)}$ gives, precisely, the right-hand-side of \eqref{eq:A-inv-main-estimate}.

    We only consider the case when $\cR(\cdot)$ is differentiable at $u$. 
    The generalization to non-differentiable points can be done by taking the limit as in the proof of Lemma \ref{lem:DR-invertible alpha}.
    In the same way in which we obtained \eqref{jacobian general} (using the computations \eqref{eq:dH(uj)/duj} and \eqref{eq:dH(uj)/duk}), we obtain that the Jacobian matrix $A = D\cR(u)$ defines a linear operator in $\R^V$ such that, for any $w\in \R^V$, the $j$-th component of $A w$ is given by
    \begin{equation}
    \label{A viscous HJ}
    (Aw)_j = \begin{cases}
 \displaystyle\sum_{\ell = 1}^d \dfrac{\nu_\ell }{\dx^2} (2 w_j - w_{j-e_\ell} - w_{j+e_\ell}) + \displaystyle\sum_{k\in \cV_j} \alpha_{jk} (w_j - w_k) + \beta_j w_j , & \text{if} \quad j\in \cJ\\
 \lambda_b w_j, & \text{if} \quad j\in \cB
  \end{cases}
\end{equation}
where $0\leq \alpha_{jk} \leq \dfrac{\KH_p}{\dx}$ and $\beta_j \geq 0$ for all $j\in \cJ$ and $k\in \cV_j$.

We construct a barrier $\{\psi_j\}_{j\in V}\in \R^V$ which depends only on the $\ell^\ast$-th variable.
We can assume, without loss of generality, that $\ell^\ast = 1$, and also that the domain is contained 
in the band $\{ (x_1, x_2, \ldots , x_d) \in \R^d\, :\ 0\leq x_1 \leq R\}$, and then we consider
$$ 
  \psi_j := 2 \phi(y_n) + \frac{1}{\lambda_b}, \qquad \text{for any}\quad  j = (n, j_2, j_3, \ldots, j_d)\in V,
$$
with $y_n:=n\dx$, $n\in \{0 ,\ldots, N\}$ and $N = \left\lfloor \frac{R}{\dx}\right\rfloor$, and 
the function $\phi:[0,R]\conv \R$ is defined as in Lemma \ref{lem:continuous barrier} below, with constants $\eta$ and $F$ to be fixed later on.

Since $\psi_j \geq \frac{1}{\lambda_b}$ for all $j\in V$, we have in particular $(A\psi)_j = \lambda_b \psi_j \geq 1$ for all $j\in \cB$. 
Due to \eqref{Cartesian grid def}, \eqref{Cartesian neighbors}, and the fact that $\psi_j$ only depends on the first variable, for any $j\in \cJ$ we have
\begin{align*}
 (A\psi)_j 
 &  =  \underbrace{\dfrac{\nu_1}{\dx^2} (2 \psi_j - \psi_{j-e_1} - \psi_{j+e_1}) }_{(\ma)}
     + \underbrace{\alpha_{j, j+e_1} (\psi_j - \psi_{j+e_1}) + \alpha_{j,j-e_1}(\psi_j - \psi_{j-e_1})}_{(\mb)}+  \underbrace{\beta_j \psi_j}_{(\gamma)}.
\end{align*}
For the $(\ma)$ term, 
by using Taylor expansions and the fact that $\phi^{(4)}\leq 0$, we obtain 
\beno  
   (\ma) & = &   - 2 \frac{\nu_1}{\dx^2} (\psi_{j-e_1} - 2 \psi_j + \psi_{j+e_1})
       \ = \  - 2 \frac{\nu_1}{\dx^2} (\phi(y_{n-1}) - 2\phi(y_{n}) + \phi(y_{n+1})) \\
      & = & - 2\nu_1 (\phi''(y_n) + \frac{\dx^2}{12} \phi^{(4)}(\xi_n)) \ \geq \  - 2\nu_1 \phi''(y_n)
\eeno
(for some $\xi_n\in [0,R]$).

For the $(\mb)$ term, using the fact that $\phi'\geq 0$,
\beno 
   (\mb)
     & \geq & - 2 \KH_p \Big|\frac{\phi(y_n)-\phi(y_{n-1})}{\dx}\Big| - 2\KH_p  \Big|\frac{\phi(y_n)-\phi(y_{n+1})}{\dx} \Big| \\
     & \geq & - 2\KH_p \Big(\frac{\phi(y_n)-\phi(y_{n-1})}{\dx}  + \frac{\phi(y_{n+1})-\phi(y_{n})}{\dx}  \Big) 
        \ =\ - 4\KH_p \frac{\phi(y_{n+1})-\phi(y_{n-1})}{2\dx} \\
     & \geq & - 4 \KH_p \Big(\phi'(y_n)  + \frac{\dx^2}{6} \phi^{(3)}(\bar\xi_n)\Big), \quad \mbox{for some $\bar \xi_n\in[0,R]$,}
       \\
     & \geq & - 4 \KH_p  |\phi'(y_n)|  -  \frac{2\KH_p\dx^2}{3} \frac{F}{\eta^2} e^{FR/\eta}.
\eeno 

Finally we have $(\gamma)=\beta_j\psi_j\geq 0$, since $\beta_j\geq 0$ and $\psi_j\geq 0$ by construction.

Now chosing $\eta:=\nu_1$ and $F:=2\KH_p$, by the previous estimates and the construction of the barrier $\phi$,
\beno
   (A \psi)_j 
     & \geq &  2 \Big(- \eta \phi''(y_n) - F |\phi'(y_n)|\Big) -   \frac{2\KH_p\dx^2}{3} \frac{F}{\eta^2} e^{FR/\eta}\\
     & \geq &  2 -  \frac{4\dx^2\KH_p^2}{3 \eta^2} e^{2\KH_p R/\eta}.
\eeno
Hence by choosing $\dx$ such that $\frac{4\dx^2\KH_p^2}{3 \eta^2} e^{2\KH_p R/\eta}\leq 1$, where $\eta:=\nu_1$, we obtain
that $(A\psi)_j\geq 1$ for all $j$ and therefore $\|A^{-1}\|_{\ell^\infty(V)} \leq \|\psi\|_{\ell^\infty(V)} \leq 2 \|\phi\|_{\ell^\infty(V)} + \frac{1}{\ml_b}$,
where $\|\phi\|_{\ell^\infty(V)}$ is bounded as in \eqref{eq:phi_bound}.
This concludes the proof of statement (ii) of the theorem.

$(iii)$ For both the $\ell^2$-error estimate and the Polyak--\L{}ojasiewicz inequality, it suffices to bound $\|A^{-1}\|_{\ell^2(V)}$ for any $A\in \partial\cR(u)$ and $u\in \R^V$. We obtain this bound by treating the interior and boundary blocks separately. Ordering $V=\cJ\cup\cB$, the operator \eqref{A viscous HJ} is block upper-triangular,
$$
  A = \begin{pmatrix} A_{\cJ\cJ} & A_{\cJ\cB} \\[2pt] 0 & \lambda_b I_{\cB} \end{pmatrix},
$$
since the boundary rows act as $w_j\mapsto\lambda_b w_j$. From \eqref{A viscous HJ}, the interior block is
$$
A_{\cJ\cJ} = \sum_{\ell = 1}^d \frac{\nu_\ell}{\dx^2} L_\ell  + A[\alpha]_{\cJ\cJ} + \operatorname{diag}(\beta_j)_{j\in\cJ},
$$
where $L_\ell$ is the symmetric Dirichlet graph-Laplacian on the interior nodes with diffusion in the $\ell$-th direction only, $\beta_j\geq 0$ by {\bf (H1)}, and $A[\alpha]: \R^V \to \R^\cJ$ is the operator
$$
(A[\alpha] w)_j = \sum_{k\in \cV_j} \alpha_{jk} (w_j - w_k), \qquad j\in\cJ .
$$

\emph{Interior block.} For any $x\in\R^{\cJ}$, using that each $L_\ell$ is symmetric, $\operatorname{diag}(\beta_j)\succeq 0$, and that $A[\alpha]_{\cJ\cJ}$ is a principal submatrix of $A[\alpha]$ (hence $\sigma_{\max}(A[\alpha]_{\cJ\cJ})\leq\sigma_{\max}(A[\alpha])$),
$$
x^\top A_{\cJ\cJ}\, x \;\geq\; \sum_{\ell = 1}^d \frac{\nu_\ell}{\dx^2}\,\lambda_{\min}(L_\ell)\,\|x\|_{\ell^2(\cJ)}^2 \;-\; \sigma_{\max}(A[\alpha])\,\|x\|_{\ell^2(\cJ)}^2 ,
$$
where $\lambda_{\min}(L_\ell)$ denotes the smallest eigenvalue of $L_\ell$, and $\sigma_{\max}(A[\alpha])$ denotes the largest singular value of $A[\alpha]$.
After a translation, $\Omega\subset[0,R]^d$, so the grid satisfies $\{x_j\}_{j\in V}\subset\{\dx j : j\in\{0,\dots,N\}^d\}$ with $N\leq R/\dx$, and the smallest eigenvalue of the one-dimensional Dirichlet Laplacian gives
$$
\lambda_{\min}(L_\ell) = 4\sin^2\!\Big(\frac{\pi}{2N}\Big) \geq \frac{4}{N^2} \geq \frac{4\dx^2}{R^2}, \qquad \ell=1,\dots,d.
$$
For $\sigma_{\max}(A[\alpha])$ we use $\sigma_{\max}(A[\alpha])\leq\sqrt{\|A[\alpha]\|_{\ell^1(V)}\,\|A[\alpha]\|_{\ell^\infty(V)}}$, the geometric mean of the maximal absolute column- and row-sums. Since $0\leq\alpha_{jk}\leq\KH_p/\dx$ by {\bf (H1)} and $|\cV_j|,|\cV_j'|\leq 2d$ on a Cartesian grid, 
the absolute row-sum of row $j$ is
$2\sum_{k\in\cV_j}\alpha_{jk}\leq 4\KH_p d/\dx$ 
and the same bound for the  absolute column-sum of column $j$.  
Hence $\sigma_{\max}(A[\alpha])\leq 4\KH_p d/\dx$
and therefore
\begin{equation}
  \label{AJJ l2 estimate}
  x^\top A_{\cJ\cJ}\, x \;\geq\; \left(\frac{4}{R^2}\sum_{\ell=1}^d\nu_\ell - \frac{4\KH_p d}{\dx}\right)\|x\|_{\ell^2(\cJ)}^2 = m\,\|x\|_{\ell^2(\cJ)}^2 ,
\end{equation}
so that $\sigma_{\min}(A_{\cJ\cJ})\geq m$ and $\|A_{\cJ\cJ}^{-1}\|_{\ell^2(\cJ)}\leq 1/m$.

\emph{Boundary split.} Concerning the block $A_{\cJ\cB} \in \R^{|\cJ| \times |\cB|}$, each entry of $A_{\cJ\cB}$ has modulus at most $\nu_{\max}/\dx^2+\KH_p/\dx$, and each interior node has at most $2d$ boundary neighbors while each boundary node has at most $2d$ interior neighbors. Therefore, using the assumed lower bound on $\lambda_b$, we have 
\begin{equation}
\label{AJB l2 estimate}
\|A_{\cJ\cB}\|_{\ell^2(\cB;\cJ)}\leq\sqrt{\|A_{\cJ\cB}\|_{\ell^1(\cB;\cJ)}\|A_{\cJ\cB}\|_{\ell^\infty(\cB;\cJ)}}\leq \dfrac{2d\,(\nu_{\max}+\KH_p\dx)}{\dx^2} \leq \lambda_b.
\end{equation}

It is easy to check that the matrix $A^{-1}$ is also upper triangular, and has the form
$$
A^{-1} = \begin{pmatrix} A_{\cJ\cJ}^{-1} & -\frac{1}{\lambda_b}A_{\cJ\cJ}^{-1}A_{\cJ\cB} \\[2pt] 0 & \frac{1}{\lambda_b} I_{\cB} \end{pmatrix}.
$$
For any $v = (v_\cJ, v_\cB)\in \R^V$, we write $A^{-1}v = ((A^{-1}v)_\cJ, (A^{-1}v)_\cB)$, where
$$
(A^{-1}v)_\cJ = A_{\cJ\cJ}^{-1} v_\cJ - \frac{1}{\lambda_b} A_{\cJ\cJ}^{-1}A_{\cJ\cB}v_\cB
\quad  \text{and} \quad
(A^{-1}v)_\cB = \frac{1}{\lambda_b}  v_\cB.
$$
Now, using \eqref{AJJ l2 estimate} and \eqref{AJB l2 estimate}, we can estimate
\begin{align*}
\| (A^{-1} v)_\cJ\|_{\ell^2(\cJ)} & \leq \| A_{\cJ\cJ}^{-1} v_\cJ\|_{\ell^2(\cJ)} + \frac{1}{\lambda_b} \| A_{\cJ\cJ}^{-1} A_{\cJ\cB} v_\cB\|_{\ell^2(\cJ)} \\
& \leq \frac{1}{m}\left( \| v_\cJ\|_{\ell^2(\cJ)} + \|v_\cB\|_{\ell^2(\cB)} \right) \leq \frac{\sqrt{2}}{m} \| v\|_{\ell^2(V)},
\end{align*}
and then,
\begin{equation}
\label{A-1 L2 bound elliptic}
\| A^{-1} v\|_{\ell^2(V)}^2 =  \| (A^{-1} v)_\cJ\|_{\ell^2(\cJ)}^2 + \| (A^{-1} v)_\cB\|_{\ell^2(\cB)}^2 \leq \left(  \frac{2}{m^2} + \frac{1}{\lambda_b^2} \right) \| v\|_{\ell^2(V)}^2.
\end{equation}
The same upper bound is true for $(A^\top)^{-1}$.

\emph{Error estimate.} As in the proof of Theorem~\ref{thm:stability-H4}, for any $u\in \R^V$, the mean-value theorem (Lemma \ref{lem:mvt}) and $\cR(\bar u) = 0$ yields $\cR(u)=A\,(u-\bar u)$ for some $A\in co(\partial\cR([u,\bar u]))$.
Therefore, we can write $u-\bar u = A^{-1} \cR(u)$, and using \eqref{A-1 L2 bound elliptic} we obtain
$$
\| u - \bar u\|_{\ell^2(V)}^2 = \| A^{-1}\cR(u)\|_{\ell^2(V)}^2 \leq \left( \frac{2}{m^2} + \frac{1}{\lambda_b^2}\right) L(u). 
$$

\emph{Polyak--\L{}ojasiewicz inequality.}
As in the proof of Theorem \ref{thm:dim-robust error estimate proper case}, for any $u\in \R^V$, every $w\in \partial L(u)$ has the form $w = 2 A^\top \cR(u)$ for some $A\in \partial\cR(u)$. Hence, \eqref{A-1 L2 bound elliptic} applied to $v = A^\top\cR(u) = \frac{w}{2}$ gives 
$$
L(u) = \| \cR(u)\|_{\ell^2(V)}^2 = \| (A^\top)^{-1} A^\top\cR(u)\|_{\ell(V)}^2 \leq \frac{1}{4} \left( \frac{2}{m^2} + \frac{1}{\lambda_b^2}\right) \| w\|_{\ell^2(V)}^2.
$$

\emph{Condition number.} For any $u\in \R^V$ and $A\in \partial\cR(u)$, since we have the estimate \eqref{A-1 L2 bound elliptic} for $\| A^{-1}\|_{\ell^2(V)}$, it only remains to upper bound $\| A\|_{\ell^2(V)}^2$. Hence, we only need to control the sum of the absolute values of the elements of $A$, row- and column-wise. 

In view of \eqref{A viscous HJ} and \eqref{eq:l2 cond visc HJ}, since $0\leq\alpha_{jk}\leq\KH_p/\dx$ by {\bf (H1)} and $|\cV_j|,|\cV_j'|\leq 2d$ on a Cartesian grid, 
the sum of the absolute values of the $j$-th row of $A$ satisfies
\begin{eqnarray*}
\sum_{k\in V} |A_{jk}| &= & \frac{4}{\dx^2} \sum_{\ell = 1}^d \nu_\ell + 2 \sum_{k\in \cV_j} \alpha_{jk} + \beta_j 
\leq   \frac{4}{\dx^2} \sum_{\ell = 1}^d \nu_\ell + 4d \frac{\KH_p}{\dx} + \KH_u \\ 
&\leq & 4 \left( \frac{1}{\dx^2} + \frac{1}{R^2} \right) \sum_{\ell=1}^d \nu_\ell  + \KH_u,   \qquad \text{if $j\in \cJ$,}
\end{eqnarray*}
and $\sum_{k\in V} |A_{jk}| = \lambda_b$ if $j\in \cJ$.
Likewise, the sum of the absolute values of the $j$-th column of $A$ can be controlled by
\begin{eqnarray*}
\sum_{k\in V} |A_{kj}| &= & \frac{4}{\dx^2} \sum_{\ell = 1}^d \nu_\ell + \sum_{k\in \cV_j} \alpha_{jk} +\sum_{k\in \cV_j'} \alpha_{kj} + \beta_j 
\leq   \frac{4}{\dx^2} \sum_{\ell = 1}^d \nu_\ell + 4d \frac{\KH_p}{\dx} + \KH_u \\ 
&\leq & 4 \left( \frac{1}{\dx^2} + \frac{1}{R^2} \right) \sum_{\ell=1}^d \nu_\ell  + \KH_u,   \qquad \text{if $j\in \cJ$,}
\end{eqnarray*}
and $\sum_{k\in V} |A_{kj}| = \lambda_b +\sum_{k\in \cV_j'} \alpha_{kj}\leq \lambda_b + \frac{2d \KH_p}{\dx} \leq 2\lambda_b$ if $j\in \cJ$.
\end{proof}

The following elementary Lemma was used in the proof of Theorem~\ref{thm:  viscous HJ}(ii). The proof can be done by simple computations.

\begin{lem}[one dimensional continuous barrier]
\label{lem:continuous barrier}
Let $\eta>0$ and  $F>0$ be two constants. Then there is a unique $\phi\in C^\infty([0,R])$ such that $\phi(0)=0$, $\phi'(R)=0$, $\phi\ge 0$, and
\begin{equation}\label{eq:phi_ode}
  -\eta \phi''(s)-F|\phi'(s)|=1, \qquad\text{for all }s\in(0,R).
\end{equation}
The function is defined by
\begin{equation}\label{eq:phi_def}
  \phi(s):= \frac{\eta}{F^2}e^{FR/\eta}\bigl(1-e^{-Fs/\eta}\bigr) -\frac{s}{F}, \qquad s\in[0,R].
\end{equation}
Moreover, it satisfies $\phi'\geq 0$ on $[0,R]$, and
\begin{equation}\label{eq:phi_bound}
  \|\phi\|_{\ell^\infty(0,R)} = \phi(R) = \frac{\eta}{F^2}\bigl(e^{FR/\eta}-1\bigr) - \frac{R}{F} \le \frac{\eta}{F^2}e^{FR/\eta}. 
\end{equation}
\end{lem}

\section{Time dependent problems}
\label{sec:time}

The time-dependent evolution equation
\begin{equation}
    \label{eq:evolution_prob}
    \partial_t u (t,x) + H(x,u,\nabla u) = 0, \quad (t,x) \in (0,T)\times \Omega \subset \R^+ \times \R^d,
\end{equation}
can formally be viewed as a special case of the stationary framework \eqref{eq:hj}
by treating $(t,x)$ as a single augmented variable.
However, the causality in the time variable allows a separate treatment.

We discretize the time derivative in \eqref{eq:evolution_prob} by a single implicit or explicit Euler step. In both cases the fully discrete problem becomes the minimization of one space-time residual functional $L(u)=\|\cR(u)\|_{\ell^2}^2$ over all nodes $(x_j,t_n)$ at once, solved by gradient flow. We study the a posteriori error of the minimizer and the convergence of the training.

The main goal of this section is to extend the theory to time-dependent problems and a comparison of the two time steppers. Since the loss is the squared $\ell^2$-norm of the residual, the relevant error estimates are those in $\ell^2$. These estimates are dimension-robust only under a CFL-type restriction on the time step, for both discretizations, and under this restriction they coincide, up to the same Gr\"onwall factor. From the $\ell^2$ perspective there is essentially no difference between minimizing the residual of the implicit or of the explicit discretization. The implicit discretization with a monotone numerical Hamiltonian is distinguished only in the $\ell^\infty$ setting since the error is controlled by the $\ell^\infty$-norm of the residual, for every $\dt$, with no accumulation in time. On the contrary, the explicit scheme enjoys this stability only under the CFL-type condition. The training efficiency of the two setups is likewise comparable, provided the implicit time step also obeys the CFL-type restriction; without it, the implicit training still converges, but its rate is no longer dimension-robust. Since dimension-robust training is the prerequisite for the method in high dimension, the analysis does not favor one stepper over the other; the choice rests on the temporal truncation error one accepts.

\subsection{Implicit discretizations}

As in Section~\ref{sec: main results stationary case}, let $\{ x_j \}_{j\in V}\subset \bar\Omega$ be the nodes on a directed graph $(V, E)$.
Let $\{ t_n\}_{n=0,\ldots,N}$ be a uniform discretization of the time interval $[0,T]$, i.e. $t_n = n\dt$ with $\dt = T/N$.
We consider the fully discrete initial-boundary value problem:
\begin{equation}\label{eq:scheme-time}
  \begin{cases}
  \cF_j^n:=
    \frac{u^{n}_j - u^{n-1}_j}{\dt}
  + \cH \big(x_j,u^{n}_j,\nabla_G u^{n}_j\big)=0,
  & \text{for} \quad (j,n) \in \cJ\times \{ 1, \ldots, N\}
  \\
  u_j^n = g_j^n, & \text{for} \quad (j,n) \in \cB\times \{ 1, \ldots, N\} \\
  u_j^0 = g_j^0, & \text{for} \quad  j\in \cJ.
  \end{cases}
\end{equation}
where the unknown $\{ u_j^n\}_{j,n}$ is a grid function in the space-time product space  
$$
\cG_T := \{ (x_j , t_n)\ : \ (j,n)\in V \times \{0,\ldots, N\}  \},
$$
and where $\cH(x_j, u_j^n, \nabla_G u_j^n)$ is a Hamiltonian function as in section~\ref{sec: main results stationary case}.
The space of grid functions in $\cG_T$ will be denoted by $\R^{V\times (N+1)}$. Throughout this section, when there is no risk of confusion, we will denote the $\ell^2$-norm in $\R^{V\times (N+1)}$ by $\| \cdot\|_{2}$.

If we consider $N=1$ and the initial values $u_j^0 = g_j^0$ being fixed, the problem \eqref{eq:scheme-time} would read as
\begin{equation}\label{eq:scheme-time 1 step}
\begin{cases}
\cF_j := \frac{1}{\dt} (u_j - g^{0}_j)
+ \cH \big(x_j,u_j,\nabla_G u_j\big)=0,
& \text{for} \quad j \in \cJ
\\
u_j = g_j^1, & \text{for} \quad j \in \cB.
\end{cases}
\end{equation}
In this context, all the theory of sections \ref{sec: main results stationary case} 
and \ref{sec: stability} applies,
provided $\cH$ satisfies {\bf (H1)} and {\bf (H2)}. 
Note that $\cF$ in \eqref{eq:scheme-time 1 step}  satisfies {\bf (H3)} naturally with $\ml=1/\dt$, due to the implicit discretization  of the time--derivative.

Arguing recursively, we can prove existence and uniqueness for the initial-boundary value problem \eqref{eq:scheme-time}. Moreover, applying Theorem~\ref{thm:exist-uniqueness-scheme-gen} recursively over the time discretization, we can prove that the unique solution is the unique critical point of the least-squares functional associated to the residual function $\cR: \R^{V\times (N+1)} \to \R^{V\times (N+1)}$ given by
\begin{equation}
\label{residual time-dep implicit}
\cR(u)_j^n := \begin{cases}
    \frac{u_j^n - u_j^{n-1}}{\dt} + \cH (x_j, u_j^n, \nabla_G u_j^n), & \text{if} \ (j,n)\in \cJ \times \{ 1, \ldots, N\} \\
    \lambda_b (u_j^n - g_j^n), & \text{if} \ (j,n)\in \cB \times \{ 1, \ldots, N\} \\
    \lambda_0 (u_j^0 - g_j^0), & \text{if} \ (j,n)\in \cJ \times \{0\}.
\end{cases}
\end{equation}
for some constants $\lambda_b, \lambda_0>0$.

\paragraph{\bf Well-posedness.}
\begin{theorem}
\label{thm:exist-uniqueness-scheme-time-implicit}
    Let $\cR(\cdot)$ be given by \eqref{residual time-dep implicit} and assume {\bf (H1)} and {\bf (H2)}. Then,
    \begin{enumerate}[label=(\roman*)]
        \item there exists a unique $\bar u\in \R^{V\times (N+1)}$ such that $\cR(\bar u) = 0$,
        \item this solution is the unique critical point of the functional $L(u) := \| \cR(u)\|^2_2$.
    \end{enumerate}
\end{theorem}

\begin{proof}
    For statement $(i)$, we proceed recursively over the time discretization. 
    Fix $u_j^0 = g_j^0$ for all $j\in \cJ$. Theorem~\ref{thm:exist-uniqueness-scheme-gen} applied to \eqref{eq:scheme-time 1 step} gives existence and uniqueness of $\{u_j^1\}_{j\in V}$.
    Assuming $\{u_j^n\}_{j\in V}$ exists and is unique, a second application of Theorem~\ref{thm:exist-uniqueness-scheme-gen} yields existence and uniqueness of $\{u_j^{n+1}\}_{j\in V}$. 
    For the one-step problem \eqref{eq:scheme-time 1 step}, condition {\bf (H3)} 
    holds with $\lambda=1/\dt$.
    By induction, $\{ u_j^n\}_{j\in V}$ exists and is unique for all $1\leq n \leq N$.  

    For statement $(ii)$, we first note that if $u\in  \R^{V\times (N+1)}$ is the solution of $\cR(u) = 0$, then it is the global minimizer of $L(u)$ and hence a critical point.
    We only need to prove that any critical point of $L(u)$ satisfies $\cR(u) = 0$.
    The functional $L(u) = \| \cR(u)\|^2_2$ can be written as
    $$
    L(u) = \sum_{n=1}^N L_n (u) + \lambda_0^2 \sum_{j\in \cJ} (u_j^0 - g_j^0)^2,
    $$
    where
    $$
    L_n (u) := \sum_{j\in \cJ} \left( \dfrac{u_j^{n} - u_j^{n-1}}{\dt} + \cH (x_j , u_j^n, \nabla_G u_j^n) \right)^2 
    + \lambda_b^2 \sum_{j\in\cB} (u_j^n - g_j^n)^2 .
    $$
    Let $u \in \R^{V\times (N+1)}$ be a critical poinut of $L(u)$, i.e. $0\in \partial L(u)$. Note that elements $w\in \partial L(u)$ are grid functions in $\R^{V\times (N+1)}$. We denote by $w^{n} := \{ w_j^n \}_{j\in V}$ the restriction of $w$ to the time $t_n$. If $w\in \partial L(u)$, then the restriction $w^{n}$ is in the generalized differential of $L(u)$ considering only the variables $\{ u_j^n\}_{j\in V}$. We shall denote this by $w^{n}\in \partial_{u^n} L(u)$.
    
    To prove that $\cR(u)^{n} = \{ \cR(u)_j^n \}_{j\in V} = 0$ for all $0\leq n\leq N$, we argue recursively again, this time backward in time.
    Note that the variables $u^{N} = \{ u_j^N\}_{j\in V}$ only appear in the term $L_N (u)$. Hence, $0\in \partial L(u)$ implies that $0 \in \partial_{u^N} L_N(u)$. Since $0\in \partial_{u^N} L_N(u)$ considers only the differential with respect to the variables $\{ u_j^N\}_{j\in V}$, we deduce from Theorem \ref{thm:exist-uniqueness-scheme-gen} that $\cR(u)^{N} = 0$, independently of the variables $\{ u_j^{N-1}\}_{j\in V}$.
    
    Now, for any $2 \leq n \leq N$, let us assume that $\cR(u)^{n} = 0$. We see that the variables $\{ u_j^{n-1}\}_{j\in V}$ appear only in the terms $L_{n}(u)$ and $L_{n-1}(u)$. Hence, we have that 
    $$
    0 \in \partial_{u^{n-1}} (L_{n} (u) + L_{n-1}(u)).
    $$
    By the same computations as in the proof of Theorem~\ref{thm:exist-uniqueness-scheme-gen}, we have
    \begin{align*}
    & \partial_{u^{n-1}} L_n(u) = 2  \left(\partial_{u_{n-1}} \cR(u)^{n}\right)^{\top} \cR(u)^{n} \\
    & \partial_{u^{n-1}} L_{n-1} (u) = 2  \left(\partial_{u_{n-1}} \cR(u)^{n-1} \right)^{\top} \cR(u)^{n-1}
    \end{align*}
    where $(\partial_{u^{n-1}} \cR(u)^{n})^{\top}$ denotes the set of transpose matrices of the generalized differential of $\cR(u)^{(n)}$ with respect to the variables $u^{(n-1)} = \{ u_j^{n-1}\}_{j\in V}$.
    Since we assumed $\cR(u)^{n} = 0$, we have that $0\in \partial_{u^{n-1}} L_{n-1} (u)$. Applying Theorem~\ref{thm:exist-uniqueness-scheme-gen} to the term $L_{n-1} (u)$, we deduce that $\cR(u)^{n-1} = 0$. And by a recurrence argument, we obtain $\cR(u)^{n} = 0$ for all $1\leq n \leq N$.

    Finally, since $u$ being a critical point implies that $\cR(u)^{n} = 0$ for all $1\leq n \leq N$, we deduce that $L(u) = \lambda_0^2\sum_{j\in \cJ} (u_j^0 - g_j^0)^2$, and taking derivatives with respect to the variables $u^{0} = \{ u_j^{0}\}_{j\in \cJ}$, we conclude that $\cR(u)^{0} = 0$.
\end{proof}

\paragraph{\bf Error estimates and the CFL condition.}
We obtain error estimates by applying Theorem~\ref{thm:stability-H3} recursively to the one-step problem \eqref{eq:scheme-time 1 step}. Since the loss is the squared $\ell^2$-norm of the residual, we want $\ell^2$-error estimates. On a finite grid one could pass through the $\ell^2$--$\ell^\infty$ norm equivalence, but its constant scales as $O(\dx^{-d})$, which is poor in high dimension. To get a dimension-robust $\ell^2$ estimate, we apply Theorem~\ref{thm:dim-robust error estimate proper case} recursively instead, which requires $\lambda = 1/\dt$ to be large relative to a quantity of order $O(1/\dx)$. Implicit Euler is unconditionally stable, so this is not a stability restriction; it is a CFL-type condition that makes the $\ell^2$ constant free of the curse of dimensionality.

\begin{theorem}
\label{thm:error-estimates-scheme-time-implicit}
    Let $\cR(\cdot)$ be given by \eqref{residual time-dep implicit} and assume {\bf (H1)} and {\bf (H2)}. Then,
    \begin{enumerate}[label=(\roman*)]
        \item for any $u\in \R^{V \times (N+1)}$ and $0\leq n \leq N$,  it holds
        \be
        \label{eq:time-dep Linf estimate}
           \| u^{n} - \bar u^{n}\|_{\ell^\infty(V)} 
           & \leq & 
           \max\bigg(\| u^0 - g^0\|_{\ell^\infty(\cJ)}, \max_{1\leq k\leq n} \|u^k - g^k\|_{\ell^\infty(\cB)}\bigg)
           \nonumber \\
           & &  \hspace{2cm}
           + \sum_{k=1}^n \dt \| \cR(u)^k \|_{\ell^\infty(\cJ)}.
        \ee
        \item if $\lambda_b > \frac{1}{\dt}$ and the following condition holds:
        $$
         \frac{\dt}{\dx_{min}} < \frac{1}{ \KH_p  \max_j | \cV'_j|}
        $$
        then, setting
        \begin{equation}\label{eq:gamma-CFL}
          \gamma := \frac{1}{1 -  \frac{\dt}{\dxmin}\, \KH_p \max_j |\cV_j'|},
        \end{equation}
        the one-step recursive estimate
        \begin{equation}\label{eq:one-step-l2}
           \| u^{n} - \bar u^{n}\|_{\ell^2(V)}
           \;\leq\;
           \sqrt{\gamma} \Big( \| u^{n-1} - \bar u^{n-1}\|_{\ell^2(\cJ)} +  \dt \, \| \cR(u)^n \|_{\ell^2(V)}\Big)
        \end{equation}
        holds for every $1\leq n\leq N$.
    \end{enumerate}
\end{theorem}
The factor $\sqrt{\gamma}$ in \eqref{eq:one-step-l2} exceeds $1$, unlike the factor $1$ in the $\ell^\infty$ bound \eqref{eq:time-dep Linf estimate}. 
The reason is that (as in Theorem \ref{thm:dim-robust error estimate proper case}), for any $u\in \R^V$ and any $A\in \partial\cR(u)$,
the dimension-robust $\ell^2$ estimate controls $A^{-1}$ through both the row and column sums of $A$, and the column-sum constant $c_2$
is smaller than the row-sum constant $c_1 = 1/\dt$ by the spatial coupling $\frac{\KH_p}{\dx_{min}}\max_j|\cV_j'|$.
\begin{proof}    
    For statement $(i)$, we compare the solution $\bar u =\{\bar u^n_{j}\}_{j,n}\in \R^{V\times (N+1)}$ of $\cR(\bar u) = 0$ and any grid function
    $u = \{ u^n_{j}\}_{j,n}\in \R^{V\times (N+1)}$. For all $n\in \{1, \ldots, N\}$, we can write
    \beno
       u_j^n + \dt \cH (x_j, u^n_j, \nabla_G u^n_j) -  \bar u_j^{n-1} =  u^{n-1}_j - \bar u^{n-1}_j + \dt\cR(u)^n_j, \quad \forall j\in \cJ
    \eeno
    (the term $-\bar u^{n-1}_j$ is added on both sides of the equation for convenience)
    and
    \beno
       \bar u_j^n + \dt \cH (x_j, \bar u^n_j, \nabla_G \bar u^n_j) -  \bar u_j^{n-1} = 0, \quad \forall j\in \cJ,
    \eeno
    (with  boundary values $u_j^n$ and $g_j^{n}$, for $j\in \cB$, respectively).
    Then (H3) holds with  $\ml\geq 1$ (for the operator $u\mapsto u + \dt\cH(x,u,\hD u)-\bar u^{n-1}$).
    Therefore, by Theorem~\ref{thm:stability-H3}, we have
    \beno
       \|u^n-\bar u^n\|_{\ell^\infty(V)} 
       & \leq & 
       \max\bigg( \frac{1}{\ml} \|u^{n-1}-\bar u^{n-1}+ \dt \cR(u)^n\|_{\ell^\infty(\cJ)},\
         \|u^n- g^n\|_{\ell^\infty(\cB)} \bigg) \\
       & \leq & 
       \max\bigg( \|u^{n-1}-\bar u^{n-1}\|_{\ell^\infty(\cJ)}, \|u^n-g^n\|_{\ell^\infty(\cB)} \bigg) 
        +   \dt \|\cR(u)^n\|_{\ell^\infty(\cJ)}.
    \eeno
    The desired bound is then obtained by using a recursion argument.

    $(ii)$ For the dimension-robust $\ell^2$--estimate we follow the same scheme as in $(i)$, but invoke Theorem~\ref{thm:dim-robust error estimate proper case} in place of Theorem~\ref{thm:stability-H3}. Fix $1\leq n\leq N$ and introduce the frozen-step residual $\cS^n:\R^V\to\R^V$,
    \beno
       \cS^n(v)_j :=
       \begin{cases}
         \frac{1}{\dt}\big(v_j - \bar u_j^{n-1}\big) + \cH(x_j, v_j, \nabla_G v_j), & j\in\cJ,\\[3pt]
         \lambda_b\,(v_j - g_j^n), & j\in\cB.
       \end{cases}
    \eeno
    This is a residual of the form \eqref{residual function def} for the Hamiltonian $\tilde\cH^n(x_j,v_j,p) := \frac{1}{\dt}(v_j-\bar u_j^{n-1}) + \cH(x_j,v_j,p)$. Adding the affine term $\frac{1}{\dt}(v_j - \bar u_j^{n-1})$ leaves {\bf (H1)} and {\bf (H2)} intact, with the same Lipschitz constant $\KH_p$ in the gradient variables and the same in-stencils $\cV_j'$, and makes {\bf (H3)} hold with $\lambda = 1/\dt$. By construction $\bar u^n$ solves $\cS^n(\bar u^n)=0$, and for any $u\in\R^{V\times(N+1)}$,
    \beno
       \cS^n(u^n) = \cR(u)^n + \frac{1}{\dt}\,\mathcal{E}^{n-1},
       \qquad
       \mathcal{E}^{n-1}_j :=
       \begin{cases} u_j^{n-1} - \bar u_j^{n-1}, & j\in\cJ,\\ 0, & j\in\cB,\end{cases}
    \eeno
    since on $\cJ$ one has $\frac{1}{\dt}(u_j^n-\bar u_j^{n-1}) + \cH(x_j,u_j^n,\nabla_G u_j^n) = \cR(u)_j^n + \frac{1}{\dt}(u_j^{n-1}-\bar u_j^{n-1})$, while on $\cB$ both sides equal $\lambda_b(u_j^n-g_j^n)=\cR(u)_j^n$. In particular $\|\mathcal{E}^{n-1}\|_{\ell^2(V)} = \|u^{n-1}-\bar u^{n-1}\|_{\ell^2(\cJ)}$.

    The hypothesis $\lambda_b > 1/\dt$ together with the CFL condition $\frac{\dx_{min}}{\dt} > \KH_p\max_j|\cV_j'|$ guarantee \eqref{eq:dimrobust-conditions} for $\cS^n$ with $\lambda = 1/\dt$, since
    \beno
       \lambda = \frac{1}{\dt} > \frac{\KH_p}{\dx_{min}}\max_{j\in\cJ}|\cV_j'|,
       \qquad
       \lambda_b > \frac{1}{\dt} > \frac{\KH_p}{\dx_{min}}\max_{j\in\cB}|\cV_j'|.
    \eeno
    Hence Theorem~\ref{thm:dim-robust error estimate proper case} applies to $\cS^n$, with
    \beno
       c_1 = \min\Big(\tfrac{1}{\dt},\lambda_b\Big) = \frac{1}{\dt},
    \eeno
    and
    \beno
       c_2 \;\geq\; \frac{1}{\dt} - \frac{\KH_p}{\dx_{min}}\max_j|\cV_j'|
            \;=\; \frac{1}{\dt}\cdot\frac{\dx_{min}-\dt\,\KH_p\max_j|\cV_j'|}{\dx_{min}}
            \;=\; \frac{1}{\dt\,\gamma},
    \eeno
    so that $c_1 c_2 \geq (\dt^2\gamma)^{-1}$. As in the proof of Corollary~\ref{cor: L is coercive}, the generalized mean-value theorem yields $A\in co\big(\partial\cS^n([\bar u^n, u^n])\big)$ with $\cS^n(u^n) = A\,(u^n-\bar u^n)$ and $\|A^{-1}\|_2 \leq (c_1 c_2)^{-1/2} \leq \dt\sqrt{\gamma}$. Therefore
    \beno
       \|u^n-\bar u^n\|_{\ell^2(V)}
       = \|A^{-1}\cS^n(u^n)\|_{\ell^2(V)}
       \leq \dt\sqrt{\gamma}\,\Big(\|\cR(u)^n\|_{\ell^2(V)} + \tfrac{1}{\dt}\|\mathcal{E}^{n-1}\|_{\ell^2(V)}\Big),
    \eeno
    which, using $\|\mathcal{E}^{n-1}\|_{\ell^2(V)} = \|u^{n-1}-\bar u^{n-1}\|_{\ell^2(\cJ)}$, is exactly the one-step estimate~\eqref{eq:one-step-l2}.
 \end{proof}

The two estimates in Theorem \ref{thm:error-estimates-scheme-time-implicit} differ in how small the residual must be for a target accuracy. Assume the initial and boundary data are matched exactly, so only the interior residual contributes. For the $\ell^\infty$ estimate $(i)$, the residual accumulates linearly in time,
$$
  \|u^n - \bar u^n\|_{\ell^\infty(V)} \;\le\; \sum_{k=1}^n \dt\,\|\cR(u)^k\|_{\ell^\infty(\cJ)} \;\le\; t_n \max_{1\le k\le n}\|\cR(u)^k\|_{\ell^\infty(\cJ)},
$$
Then $\|u^n - \bar u^n\|_{\ell^\infty(V)} < \epsilon$ holds once $\max_{1\le k\le n}\|\cR(u)^k\|_{\ell^\infty(\cJ)} < \epsilon/t_n$. Since $t_n\le T$, the mesh-independent bound $\|\cR(u)^k\|_{\ell^\infty(\cJ)} < \epsilon/T$ suffices. For the $\ell^2$ estimate $(ii)$, iterating \eqref{eq:one-step-l2} gives
$$
  \|u^n - \bar u^n\|_{\ell^2(V)} \;\le\; \dt\sum_{k=1}^n \gamma^{(n-k+1)/2}\,\|\cR(u)^k\|_{\ell^2(V)} \;\le\; t_n\,\gamma^{n/2}\max_{1\le k\le n}\|\cR(u)^k\|_{\ell^2(V)},
$$
so $\|u^n-\bar u^n\|_{\ell^2(V)} < \epsilon$ now requires $\max_{1\le k\le n}\|\cR(u)^k\|_{\ell^2(V)} < \epsilon/(t_n\,\gamma^{n/2})$, a smaller bound. The factor $\gamma^{n/2}>1$ has no counterpart in $(i)$. It is largest at the final time, where $\gamma^{N/2}\to \exp\!\big(\tfrac{T\,\KH_p\max_j|\cV_j'|}{2\,\dx_{min}}\big)$ as $\dt\to0$, and it grows under spatial refinement. So the $\ell^\infty$ error is controlled by a mesh-independent residual tolerance, while the $\ell^2$ error needs the residual, and hence the loss $L(u)=\|\cR(u)\|_2^2$, below a threshold smaller by the factor $\gamma^{n/2}$.

\paragraph{\bf Training convergence.}
The discrete solution minimizes the joint loss $L(u)=\|\cR(u)\|_2^2$ over $\R^{V\times(N+1)}$ by gradient flow. The space-time Jacobian is block lower-bidiagonal in time, with diagonally dominant, invertible diagonal blocks $\tfrac1\dt I + (\text{spatial part})$, so $L(\cdot)$ has a unique global minimizer $\bar u$ and gradient flow converges linearly through the Polyak--\L{}ojasiewicz inequality of Theorem~\ref{thm:dim-robust error estimate proper case}. The conditioning that sets the rate is computed for the explicit residual in Section~\ref{sec:explicit}; it has the same form here, with per-step amplification $(1-\beta)^{-1/2}$, where $\beta$ is the CFL number.

\subsection{Explicit discretizations}
\label{sec:explicit}
We now consider explicit time-marching schemes of the form
\begin{equation}\label{eq:scheme-time explicit}
\frac{u^{n}_j - u^{n-1}_j}{\dt_n}
+ \cH \big(x_j,u^{n-1}_j,\nabla_G u^{n-1}_j\big)=0,
\qquad \text{for} \quad (j,n) \in \cJ\times \{ 1, \ldots, N\},
\end{equation}
with the corresponding initial and boundary conditions as in \eqref{eq:scheme-time}.
Proving the existence and uniqueness of a solution for the initial-boundary value problem associated with \eqref{eq:scheme-time explicit} is straightforward.

Let us define the associated residual function
$$
\cR (u)_j^n:= \begin{cases}
    \dfrac{u_j^n - u_j^{n-1}}{\dt_n} + \cH (x_j, u_j^{n-1}, \nabla_G u_j^{n-1}), & \text{if} \ (j,n)\in \cJ \times \{1, \ldots, N\} \\
    \lambda_b (u_j^n - g_j^n), & \text{if} \ (j,n) \in \cB \times \{ 0, \ldots, N\} \\
    \lambda_0 (u_j^0 - g_j^0), & \text{if} \ (j,n)\in \cJ \times \{0\}.
\end{cases}
$$
\paragraph{\bf Well-posedness.}
A similar analysis as for implicit schemes proves that, if $u\in \R^{V\times (N+1)}$ is a critical point of the least-squares functional $L(u):= \| \cR(u)\|_2^2$, then the associated residual satisfies $w_j^n : = \cR(u)_j^n = 0$, for all $(j,n)$. Indeed, we observe that the terms $u_j^N$ only appear in $\cR(u)_j^N$, whence $\partial_{u_j^N} L(u) = 0$ implies $w_j^N = 0$ for all $j\in V$. Then, from the optimality condition $\partial_{u_j^n}L(u)=0$, one can quickly derive an explicit recurrence relation for the residual  $w^n_j= \mathcal{R}(u)^n_j$. For $j\in\cJ$ and $n\in \{1, \ldots,N-1\}$, the optimality condition reads as
\begin{align*}
     w^n_j =& \left(\frac{\dt_n}{\dt_{n+1}} - \dt_n \frac{\partial \mathcal{H}}{\partial u}(x_j, u^n_j, \nabla_G u^n_j)  -\sum_{k\in\cV_j} \frac{\dt_n}{\dx_{jk}} \frac{\partial \cH}{\partial p_k}(x_j, u^n_j, \nabla_G u^n_j) \right) w^{n+1}_j\\
      & +\sum_{k\in\cV_j'} \frac{\dt_n}{\dx_{kj}} \frac{\partial \cH}{\partial p_j} (x_k, u^n_k, \nabla_G u^n_k)w^{n+1}_k,
\end{align*}
where the last sum runs over the in-stencil $\cV_j' = \{ k\in\cJ : j\in\cV_k\}$, i.e. over the nodes whose stencils contain $x_j$. For $j\in\cB$, and for $n=0$, analogous relations express $\lambda_b w^n_j$ and $\lambda_0 w^0_j$ as linear combinations of $\{ w^{n+1}_k \}_{k\in \cV_j'}$.
Therefore, since the update for $w^n_j$ requires only the values $w^{n+1}_k$ at the next time level, and $w^N_j=0$ for all $j\in V$, we see that $w_j^n\equiv 0$; i.e., the residual at the critical point must be $0$.
In particular, as for the implicit scheme, the unique critical point of the loss is the discrete solution $\bar u$ for every $\dt$.

\paragraph{\bf Error estimates.}
Under the CFL condition the explicit Euler map is a contraction in $\ell^\infty$, so the estimates of Theorem~\ref{thm:error-estimates-scheme-time-implicit} carry over: the $\ell^\infty$ error is controlled by the residual with no accumulation in time, and the $\ell^2$ error up to the same Gr\"onwall factor $\gamma^{n/2}$. Unlike the implicit case, both require the CFL condition.

\paragraph{\bf Training convergence.}
The discrete solution minimizes the joint loss $L(u)=\|\cR(u)\|_2^2$ over $\R^{V\times(N+1)}$ by gradient flow. The space-time Jacobian $D\cR$ is block lower-bidiagonal in time and invertible for every $\dt$, so $L(\cdot)$ has a unique global minimizer $\bar u$, and gradient flow converges linearly through the Polyak--\L{}ojasiewicz inequality of Theorem~\ref{thm:dim-robust error estimate proper case} (see \cite{karimi2016linear}). The rate is set by the conditioning of $D\cR^\top D\cR$, which we now compute. We order the unknowns by time level, writing $u = (u^0, \ldots, u^N)\in\R^{V\times(N+1)}$, and write linear maps on $\R^{V\times(N+1)}$ in the corresponding block form: for a matrix $A$, the block $A_{n,m}\in\R^{V\times V}$ maps time level $m$ to time level $n$, so that $(D\cR)_{n,m} = \partial\cR(u)^n/\partial u^m$. To keep all rows of $\cR$ on the same scale, we take a uniform step $\dt_n\equiv\dt$ and $\lambda_b = \lambda_0 = 1/\dt$; other choices of the weights modify only the boundary and initial rows of the diagonal blocks below, by bounded factors. Define the one-step Euler maps $F^n:\R^V\to\R^V$ by
$$
F^n(v)_j := \begin{cases}
   v_j - \dt\,\cH(x_j,v_j,\nabla_G v_j), & j\in\cJ,\\
   g_j^n, & j\in\cB,
\end{cases}
$$
so that $\cR(u)^0 = \tfrac1\dt(u^0 - g^0)$ and $\cR(u)^n = \tfrac1\dt\big(u^n - F^n(u^{n-1})\big)$ for $1\leq n\leq N$. Differentiating, the diagonal blocks of $D\cR(u)$ are $\tfrac1\dt I$, the only other nonzero blocks are the subdiagonal ones, $(D\cR)_{n,n-1} = -\tfrac1\dt M_{n-1}$, where
$$
M_n := DF^{n+1}(u^n)\in\R^{V\times V}
$$
agrees with $I - \dt\,D\cH(u^n)$ on the rows indexed by $j\in\cJ$ and has zero rows for $j\in\cB$ (the boundary rows of $F^{n+1}$ are constant in $v$). The Jacobian therefore factors as
\begin{equation}\label{eq:DR-explicit-factor}
   D\cR = \tfrac1\dt\,(I-B),
   \qquad
   B_{n,m} := \begin{cases}
      M_{n-1}, & m = n-1,\\
      0, & \text{otherwise},
   \end{cases}
\end{equation}
with $B$ strictly lower block-bidiagonal in time. Being strictly block lower-triangular, $B$ is nilpotent, $B^{N+1}=0$, so $I-B$, and with it $D\cR$, is invertible for every $\dt$, with inverse given by the finite Neumann series $(I-B)^{-1} = \sum_{k=0}^N B^k$. Since $B$ has nonzero blocks only on the first block subdiagonal, $B^k$ has nonzero blocks only on the $k$-th one: $(B^0)_{n,m} = \delta_{n,m}\,I$ and, for $k\geq 1$,
$$
(B^k)_{n,m} = B_{n,n-1}\,B_{n-1,n-2}\cdots B_{m+1,m} = M_{n-1}M_{n-2}\cdots M_m,
\qquad \text{if} \ n-m=k,
$$
and $(B^k)_{n,m} = 0$ otherwise. Summing over $k$ contributes exactly one block to each position $n\geq m$:
\begin{equation}\label{eq:DR-explicit-inverse}
   (D\cR)^{-1} = \dt\sum_{k=0}^{N} B^k,
   \qquad
   \big((D\cR)^{-1}\big)_{n,m} =
   \begin{cases}
      \dt\,\Phi_{n,m}, & 0\leq m\leq n\leq N,\\
      0, & m>n,
   \end{cases}
\end{equation}
where
$$
\Phi_{n,m} := M_{n-1}M_{n-2}\cdots M_m \quad (n>m), \qquad \Phi_{n,n} := I.
$$
The matrix $\Phi_{n,m}$ is the discrete propagator from $t_m$ to $t_n$: it maps a perturbation of the state at time $t_m$ to the resulting perturbation at time $t_n$ under the linearized explicit scheme, and satisfies the semigroup property $\Phi_{n,m} = \Phi_{n,l}\,\Phi_{l,m}$ for $m\leq l\leq n$. This gives two facts. First, $D\cR^\top D\cR = \tfrac1{\dt^2}(I - B - B^\top + B^\top B)$ is symmetric and block-tridiagonal in time, that is, elliptic in time. Already for $\cH\equiv 0$ (so $M_n=I$) it is the one-dimensional discrete Laplacian in time, whose condition number scales like $N^2$. Second, the spatial coupling enters the rate only through $\|\Phi_{n,m}\|_2$, hence through $\|M_n\|_2$. By the monotonicity {\bf (H2)}, $M_n$ is entrywise nonnegative under the CFL condition, with row sums $\|M_n\|_\infty \le 1$ and column sums $\|M_n\|_1 \le 1+\beta$, where $\beta:=\dt\,\KH_p\max_j|\cV_j'|/\dx_{min}$. Hence
\begin{equation}\label{eq:Mn-bound}
   \|M_n\|_2 \le \sqrt{\|M_n\|_1\,\|M_n\|_\infty} \le (1+\beta)^{1/2},
\end{equation}
and by submultiplicativity, $\|\Phi_{n,m}\|_2 \leq (1+\beta)^{(n-m)/2}$.
This per-step amplification involves only the local in-stencil $\max_j|\cV_j'| = O(d)$, never $|V|=O(\dx^{-d})$, so it is dimension-robust. The condition number is therefore a time-elliptic factor of order $N^2$ times a dimension-robust propagator factor. The implicit residual \eqref{residual time-dep implicit} has the same block-bidiagonal form, with diagonal blocks $D_n$ that agree with $\tfrac1\dt I + D\cH(u^n)$ on the interior rows, and its per-step propagator $\tfrac1\dt D_n^{-1}$ obeys the analogous bound $\|\tfrac1\dt D_n^{-1}\|_2\leq \sqrt\gamma = (1-\beta)^{-1/2}$, since $\|D_n^{-1}\|_2\le\dt\sqrt\gamma$ by Theorem~\ref{thm:dim-robust error estimate proper case}.
Although $(1+\beta)^{1/2} < (1-\beta)^{-1/2}$ for $0<\beta<1$, this does not mean that the explicit scheme is more stable. Both quantities are worst-case bounds, derived from the same $\ell^1$--$\ell^\infty$ interpolation. 
In both cases the row sums are bounded by $1$, and the spatial coupling shifts the column sums by $\beta$. For the explicit scheme, the propagator $M_n$ is bounded directly, giving the factor $1+\beta$. The implicit propagator is instead the inverse of $\dt\,D_n$, so the factor $1-\beta$ enters through its reciprocal. 
The bound obtained for the implicit scheme is therefore not sharp: it degenerates as $\beta\to 1$, even though the scheme itself does not. Moreover, accumulated over $n-m$ steps, the two bounds yield the same leading-order Gr\"onwall factor $\exp\!\big(\tfrac{T\,\KH_p\max_j|\cV_j'|}{2\,\dx_{min}}\big)$ as $\dt\to 0$. Indeed, for $\beta\geq 1$ the implicit scheme remains stable---only this particular $\ell^2$ bound becomes unavailable---while the explicit scheme is unstable.

\section{Examples of monotone numerical Hamiltonians}
\label{sec:monotone-hamiltonians}

\newcommand{\LF}{\scriptscriptstyle LF}
\newcommand{\UU}{\scriptscriptstyle U}
\newcommand{\Hu}{{\bf (H1)}}
\newcommand{\Hd}{{\bf (H2)}}
\newcommand{\Ht}{{\bf (H3)}}
\newcommand{\Hq}{{\bf (H4)}}

We discuss different standard classes of numerical Hamiltonians that 
satisfy the hypotheses {\bf (H1)}--{\bf (H4)}, in particular the monotonicity condition {\bf (H2)}, 
and give also second-order non-monotone variants.
In this section, we consider mainly Cartesian grids $\{x_i\}_{j\in V}$ of $\R^d$, where $V\subset \Z^d$ is a multi-index set. In other words, we assume \eqref{Cartesian grid def} and \eqref{Cartesian neighbors}.

We denote by $(e_k)_{1\leq k\leq d}$ the canonical basis, and $(\dx_k)_{1\leq k\leq d}$ the mesh steps in each coordinate direction.
In the entire section, we shall use the notation
\begin{equation}\label{eq:one-sided-differences}
    D^{\pm} u_i = \{D_k^\pm u_i\}_{1\leq k\leq d}\in \R^{2d}, \quad \text{where} \quad D_k^\pm u_i := \frac{u_{i} - u_{i \pm e_k}}{\dx_k} \qquad \text{for} \quad  k=1,\ldots,d.
\end{equation}

\paragraph{\bf Lax--Friedrichs scheme (LF).}
Let us consider a first-order PDE of the form 
\be\label{eq:H---LFscheme}
    H(x,u(x),\nabla u(x))=0,
\ee 
where $H(x,u,p)$ is a Lipschitz continuous function, non-decreasing with respect to $u$.

The Lax--Friedrichs numerical Hamiltonian is defined by
\begin{equation}
\label{LxF scheme def}
\cH_{\LF}(x_i, u_i, D^{\pm} u_i) := H\!\left(x_i, u_i, \frac{D^{-} u_i - D^{+} u_i}{2}\right) 
       + \sum_{k=1}^d c_k(x_i) \frac{D_k^{-} u_i + D_k^{+} u_i}{2},
\end{equation}
where  the function $c_k: \Omega \to \R^+$ is chosen such that 
$c_k(x) \geq \subalpha + \left\|\frac{\partial H}{\partial p_k}(x,u,\cdot)\right\|_\infty$ for all $k$,
for some $\subalpha>0$. This ensures that $\cH_{\LF}$ satisfies {\bf (H1)}, {\bf (H2)} and {\bf (H4)}.
The consistency of the scheme with the PDE \eqref{eq:H---LFscheme}
amounts to having $\cH_{\LF}(x,u,p^\pm)= H(x,u,p)$ whenever $p^- = - p^+ = p$.
Other variants of the LF scheme ensuring monotonicity are also possible (see \cite{osher1991high}).

\paragraph{\bf Upwind scheme (U).}
Let us consider a PDE of the form
$$ H(x,u,\nabla u) := \min_{a\in A}
   \max_{b\in B} \bigl(f(x,a,b) \cdot \nabla u  + r(x,a,b) u  - \ell(x,a,b)\bigr), 
$$
where $f: \Omega \times A\times B \to \R^d$ and $r,\ell:\Omega\times A\times B\to \R$ are Lipschitz functions {(with $A,B$ two non-empty compact sets), so that $H(x,u,p)$ is well-defined for all $(x, u, p)\in \Omega\times \R\times \R^d$.}
We also assume that $r(x,a,b)\geq 0$ for all $(x,a,b)$, which implies that the function $H(x,u,p)$ is non-decreasing with respect to $u$. 

For this case, we may consider the following upwind scheme:
$$
\begin{aligned}
  \cH_{\UU}(x_i,u_i,D^\pm u_i)
  :=
  \min_{a\in A} \max_{b\in B} 
  \Big(
  \sum_{k=1}^d & \bigl[f_k(x_i,a,b)_{+} D^-_ku_i + f_k(x_i,a,b)_{-} D^+_k u_i \bigr] \\
  & + r(x_i,a,b) u_i
  - \ell(x_i,a,b)
  \Big),
  \label{eq:H---Upwind-b}
  \end{aligned}
$$
where $f_k(x,a,b)_+ = \max \left( 0,\,  f_k(x,a,b) \right)$ and $f_k(x,a,b)_- = -\min \left( 0,\,  f_k(x,a,b) \right)$.
Note that $\cH_{\UU}$ satisfies \Hu, \Hd, and {\bf (H3)} provided $r(x, a,b)\geq \lambda$ for some $\lambda>0$.

For the classical eikonal equation, where $f(x,a,b) = c(x) b$ and $B =  \mathbb{S}^{d-1}$,  $r(x,a,b)\equiv 0$ 
and $\ell(x,a,b)\equiv 1$, we have
\begin{equation}
\label{eikonal positive}
  H(x,u,\nabla u)= |c(x)|\,  \|\nabla u\|_2 - 1.
\end{equation}
The associated upwind Hamiltonian in this case is
\begin{equation}
\label{upwind scheme eikonal}
  \cH_{\UU}(x_i,u_i,D^\pm u_i) =
  |c(x_i)| \left( \sum_{k=1}^d \left[\max( D^-_k u_i, D^+_k u_i,0)\right]^2 \right)^{\frac{1}{2}} - 1.
\end{equation}
This Hamiltonian satisfies \Hu\ and \Hd. However, it is neither proper nor uniformly elliptic.\\

If instead we consider $f(x,a,b) = c(x) a$ and $A =  \mathbb{S}^{d-1}$ and $\ell (x,a,b)\equiv -1$, we obtain the Hamiltonian
\begin{equation}
\label{eikonal negative}
H(x,u,\nabla u)= - |c(x)|\,  \|\nabla u\|_2 + 1.
\end{equation}
Note that, although the PDEs $H(x,u,\nabla u) = 0$, with $H$ given by \eqref{eikonal positive} and \eqref{eikonal negative}, appear to be the same, they are not in the sense of viscosity solutions. For instance, if $c(x) = 1$, and we consider zero boundary condition, the viscosity solution of \eqref{eikonal positive} is the distance to the boundary, whereas the viscosity solution of \eqref{eikonal negative} is the negative distance to the boundary. 

In the case \eqref{eikonal negative}, the associated upwind Hamiltonian
is given by
$$
\cH_{\UU}(x_i,u_i,D^\pm u_i) =
  - |c(x_i)| \left( \sum_{k=1}^d \left[\min( D^-_k u_i, D^+_k u_i,0)\right]^2 \right)^{\frac{1}{2}} + 1.
$$
For terms involving the $\ell^1$-norm of the gradient $\pm \|\nabla u\|_1$, which corresponds to considering the control set $A$ (resp. $B$) being the $\ell^\infty$-unit sphere $A = \{ a\in \R^d\, : \ \|a\|_\infty = 1\}$, the upwind numerical approximation reads as 
\be  \label{eq:normu-1-hamiltonian}
   \| \nabla u\|_1 \approx \sum_{k=1}^d \max(p^-_k,p^+_k,0), \quad  \text{and} \quad
   -\| \nabla u\|_1 \approx \sum_{k=1}^d \min(p^-_k,p^+_k,0).
\ee

\medskip
\paragraph{Enforcing properness via the Kruzhkov transform.}

Consider the eikonal equation 
\be 
\label{eq:eikonal-0}
c(x) \|\nabla u(x)\|=1, \quad x\in \mO,\qquad
u=0 \ \mbox{on} \ \partial \mO,
\ee 
with $c(x)\geq 0$. Then, the properness condition {\bf (H3)} does not hold, both for Lax-Friedrichs \eqref{LxF scheme def} and Upwind \eqref{upwind scheme eikonal} discretizations.

However, the Kruzhkov transform 
$v:=1-e^{-u}$ yields
\be \label{eq:eikonal-1}
c(x) \|\nabla v(x)\| + v(x) - 1 = 0, \quad x\in \mO,\qquad
v=0 \ \mbox{on} \ \partial \mO.
\ee 
For this PDE, both the Lax-Friedrichs and Upwind schemes are proper with $\ml=1$.\\

More generally, consider
\be\label{eq:general-infinite-0}
\min_{a} \max_b 
\left(f(x,a,b) \cdot \nabla u(x) + r(x,a,b) u(x)  - \ell(x,a,b)\right) = 0,
\quad x\in \mO,
\qquad
u=0\ \text{on}\ \partial \mO,
\ee
with $r\geq 0$ and $\ell(x,a,b)\ge \underline \ell>0$,
and let us assume that $u\in [0,u_{\max}]$.
Then the transform $v:=\frac{1}{\mu} (1-e^{-\mu u})$, with $\mu>0$, leads to
\be\label{eq:general-infinite-1}
\hat H(x,v,\nabla v) := \min_{a} \max_b \hat H_{a,b}(x,v,\nabla v) = 0, \quad x\in \mO,
\qquad
v=0\ \text{on}\ \partial \mO,
\ee
where
$$
\hat H_{a,b}(x,v,\nabla v) :=
f(x,a,b) \cdot \nabla v 
+ \frac{r(x,a,b)}{\mu} (1-\mu v) \log \left(\frac{1}{1-\mu v}\right)
- (1 - \mu v)\ell(x,a,b).
$$ 
Hence,
provided that $\mu \geq \frac{1}{u_{\max}}$,
or $r\equiv 0$,  we have
$$
\frac{\partial \hat H_{a,b}}{\partial v} (x,v,\nabla v)
= \mu \ell(x,a,b) + r(x,a,b) \left(1- \log\left( \frac{1}{1-\mu v} \right)\right) \ \geq\  \mu \underline \ell \ > \  0
$$
(since $\mu u = \log \left( \frac{1}{1 - \mu v}\right)$).
In these cases, the associated numerical Hamiltonian $\cH$ satisfies \Ht\ with constant $\ml:=\mu \underline\ell$.

\section{Numerical algorithm}
\label{sec:numerical algos}

Our strategy to approximate the viscosity
solution of a Hamilton-Jacobi equation consists in training a neural network $\Phi(x;\theta)$  by minimizing the residual loss of the form
$$
  L(\Phi(\cdot;\theta)) = \sum_{j\in \cJ}
    \bigl(\cH(x_j, \Phi(x_j;\theta), \nabla_G \Phi(x_j;\theta))\bigr)^2
    + \mu_b \sum_{j\in \cB} \bigl(\Phi(x_j;\theta) - g(x_j)\bigr)^2,
$$
using SGD or a variant such as Adam. At each step, a minibatch of
collocation points $\hat\cJ_t \subset \cJ$ and $\hat\cB_t \subset \cB$ is
drawn, and the gradient with respect to $\theta$ is computed by automatic
differentiation. Note that, if we consider a Cartesian mesh of a $d$-dimensional domain, each
collocation point requires only $2d+1$ neural network evaluations.

Because $\Phi(\cdot;\theta)$ is defined on all of $\overline\Omega$, the
collocation points do not have to lie on a fixed grid. At each step, we can
sample them from $\Omega$ and its boundary directly, and treat each sample
$x$ as the central node of a local Cartesian stencil $x + \dx\, \Z^d$. To keep this stencil inside $\overline\Omega$, the interior collocation points are sampled at distance at least $\dx$ from the boundary $\partial\Omega$, while the boundary penalty is evaluated at points sampled on $\partial\Omega$ itself; this prevents the finite-difference stencil from reaching across the boundary.
The network is therefore trained against the residual of the scheme on a
family of admissible Cartesian graphs rather than on a single fixed one.
The results in Section \ref{sec: stability} suggest that minimizing a functional based on a numerical Hamiltonian is, in general, easier if the discretization step $\dx$ is larger. Note that Theorems \ref{thm:dim-robust error estimate proper case} and \ref{thm:  viscous HJ} ensure linear convergence of the gradient flow (with a dimension robust convergence rate) through a Polyak-\L{}ojasiewicz inequality, provided $\dx$ is sufficiently big. However, to obtain a more accurate approximation of the solution, we aim to minimize a finite-difference residual on a finer grid.

\subsection*{A multi-level training algorithm}
Let $\cH_{\alpha, \lambda, \dx} (\cdot)$ be a parametric family of monotone Hamiltonians, as the ones presented in Section \ref{sec:monotone-hamiltonians}, where $\dx$ is the discretization step, and $\lambda, \alpha$ are parameters related to the monotonicity properties \Ht\  and \Hq\ respectively.
We propose a multi-level strategy, optimizing a sequence of loss functionals
\begin{equation}
    L^{(k)}(\Phi(\cdot;\theta)) := L(\Phi(\cdot; \theta) ; \Delta x^{(k)}, \lambda^{(k)}, \alpha^{(k)}),~~~k=0,1,...
\end{equation}
with $\Delta x^{(k)}$, $\lambda^{(k)}$, and $\alpha^{(k)}$ all decreasing in $k$. The first level uses the coarsest grid and the largest damping
and artificial viscosity. The parameters there are chosen so that the
well-posedness theory of the previous sections applies: the numerical Hamiltonian is monotone, the residual operator is a homeomorphism, and the unique
critical point of the loss functional is the unique solution of the discrete problem. Moreover, in view of Theorems~\ref{thm:dim-robust error estimate proper case} and \ref{thm:  viscous HJ}, by choosing appropriate parameters in $\cH_{\alpha, \lambda,\dx} (\cdot)$, one can ensure linear convergence of gradient descent (through a Polyak-\L{}ojasiewicz inequality), with a convergence rate which does not suffer from the curse of dimensionality.
Subsequent levels refine the grid and reduce $\lambda, \, \dx$ and $\alpha$ toward
their target values. Once the iterate is close enough to the viscosity
solution, we are free to leave the monotonicity regime at the finer
levels: we may, for instance, switch to a higher-order scheme for which
the critical-point-uniqueness guarantee of
Theorem~\ref{thm:exist-uniqueness-scheme-gen} no longer applies. In
practice, each level also carries its own minibatch configuration. Finer
grids have more collocation points, so the minibatch size should grow
with the level to keep the variance of the stochastic gradient under
control.

In classical mesh-based settings, multi-level strategies require an
interpolation or prolongation operator to transfer a coarse-grid solution
to a finer grid. Neural networks sidestep this step, since the trained
network is defined on the entire domain and serves as its own
prolongation. This is especially useful in high-dimensional problems,
where classical interpolation on grids is already prohibitive.

\subsubsection*{The benefit of level-wise warm-starting}

The benefit of warm-starting each level from the preceding one is rooted
in how the conditioning of the discrete loss landscape changes with the
mesh. We analyze this first for full gradient descent on grid functions (see the numerical experiments in section \ref{subsec: full grid L2 flow}),
an idealized setting that ignores the cost of transferring between grids.
We then comment on the neural-network regime, which is what we actually
use.

\medskip\noindent\textbf{Grid-function regime.}
By Corollary~\ref{cor: L is coercive},
$L(u;\dx) = \|\cR_\dx(u)\|_2^2$ is coercive in the $u$ variable with a unique critical
point $u_\dx^\star$. Using the equivalence between $\ell^2$ and $\ell^\infty$-norms in the space of grid functions $\R^{V_\dx}$, or the finer $\ell^2$ estimates in Theorems \ref{thm:dim-robust error estimate proper case} and \ref{thm:  viscous HJ} when applicable, one can obtain a Polyak-Łojasiewicz inequality of the from
$$
L(u; \dx) = \| \cR_\dx (u)\|_2^2 \leq C(\dx) \| D\cR_\dx(u)^\top \cR_\dx(u) \|_2^2 = \frac{C (\dx)}{4} \| \nabla L(u; \dx) \|_2^2,
$$
which ensures linear convergence of gradient descent applied to $L(u)$ with rate depending on $C (\dx)$ (see \cite{karimi2016linear}).
More precisely, by choosing an appropriate step-size, the convergence rate is given by the condition number $\kappa (\dx) := K_L(\dx) \, C(\dx)$, where $K_L(\dx)$ is the Lipschitz constant of $\nabla L(u; \dx)$ with respect to $u$.
In view of the results in section \ref{sec: stability}, $C(\dx)$ is a decreasing function of $\dx$.

If we consider $\dx_1>0$, the initialization grid function $u_0 \in \R^{V_{\dx_1}}$, and apply $k_1$ iterations of gradient descent to $L(u_0; \dx_1)$ with learning rate $\eta_1 = \frac{1}{K_L(\dx_1)}$,  we obtain a function $u_1\in \R^{V_{\dx_1}}$ such that
\be\label{eq:k1-iterations}
L(u_1; \dx_1) \leq \left( 1 - \frac{1}{\kappa (\dx_1)} \right)^{k_1} \, L(u_0; \dx_1).
\ee
See for instance \cite[Theorem 1]{karimi2016linear}.

Let us now consider $\dx_2<\dx_1$ and a function $u: \bar \Omega \to \R$ such that $u(x_j) = u_{1,j}$, $\forall j\in V_{\dx_1}$.  Assuming that $u$ is Lipschitz and piecewise $C^2$, with a possibly singular region $\Sigma$ of at most co-dimension one, we can show that for any $0< \dx_2\leq \dx_1$,
\be\label{eq:interpol-estimate}
   |L(u;\dx_2) - L(u; \dx_1)| \leq C(u) \dx_1,
\ee  
where $C(u)>0$ depends on the Lipschitz constant of $u$ and the maximum of $D^2 u$ a.e. in $\Omega$.
Indeed, in the region where $u$ is regular then such a bound is easily obtained, in terms of the second order differential of $u$. In a region $\Sigma_{\dx_1}=\{x : \  d(x,\Sigma)\leq \dx_1\}$ where $u$ can 
be singular, the error between the residuals is bounded only by the Lipschitz contants of $u$, but the volume of the region is $O(\dx_1)$ compared to the volume of the whole region $\mO$. This leads to the desired estimate.
In the neural network setting, it is enough to assume that the neural network is a piecewise $C^2$ function with parameters (weights and biases) in a compact set.

With a slight abuse of notation we also use $u_1$ to denote the grid function obtained by evaluating $u(\cdot)$ in the finer mesh $\{ x_j\}_{j\in V_{\dx_2}}$, i.e. $u_1 = (u(x_j))_{j\in V_{\dx_2}}\in \R^{V_{\dx_2}}$.
Under \eqref{eq:interpol-estimate}, applying $k_2$ iterations of gradient descent to $L(u_1; \dx_2)$ starting from $u_1$ with learning rate $\eta_2 = \frac{1}{K_L (\dx_2)}$, combining \eqref{eq:k1-iterations} and \eqref{eq:interpol-estimate},
we obtain a grid function $u_2\in \R^{V_{\dx_2}}$ such that
\begin{eqnarray*}
    L(u_2; \dx_2) & \leq & \left( 1 - \frac{1}{\kappa(\dx_2)} \right)^{k_2} \, L(u_1; \dx_2) \\
    & \leq & \left( 1 - \frac{1}{\kappa (\dx_2)} \right)^{k_2} \left( L(u_1; \dx_1) + C(u_1) \dx_1 \right) \\
    & \leq & \Big( 1 - \frac{1}{\kappa (\dx_2)} \Big)^{k_2}  \left( \Big(1-\frac{1}{\kappa (\dx_1)}\Big)^{k_1} L(u_0,\dx_1) + C(u_1) \dx_1 \right).
\end{eqnarray*}
Note that $u_2$ is obtained after applying $k_1$ iterates of the gradient descent with mesh step $\dx_1$ (starting from $u_0$), and then $k_2$ iterates with mesh step $\dx_2$ (starting from $u_1$).

Since, as we proved, the condition number grows under mesh refinement (i.e. $\kappa(\dx_2) > \kappa(\dx_1)$ for $\dx_2 < \dx_1$), instead of directly doing gradient steps with mesh step $\dx_2$,
we can see that we have more interest in doing first as many as possible gradient steps with rough mesh step $\dx_1$
(meaning a loss reduction $(1-\frac{1}{\kappa(\dx_1)})$ at each step), untill the loss is sufficiently small 
or of the order of $C\dx_1$, before switching to a finer mesh step $\dx_2$, where the loss reduction per iteration is much smaller, in order to reach a similar loss reduction.

The same idea can be generalized to multi-level strategies, providing an efficient way to minimize residuals corresponding to fine discretizations. In the neural network framework, the interpolation step (which is the main bottleneck in high-dimensional settings) is sidestepped by using the neural network as interpolator.

\medskip\noindent\textbf{Neural-network regime.}
When $u$ is replaced by a neural network $\Phi(\cdot;\theta)$, gradient
flow on $\theta \mapsto L(\Phi(\cdot;\theta), \dx)$ reads
$$
  \dot\theta \;=\; -\,\partial_\theta \Phi^\top\,
    D\cR_\dx(\Phi (\cdot;\theta))^\top\, \cR_\dx(\Phi(\cdot;\theta)),
$$
where $\partial_\theta \Phi$ is the Jacobian of the Neural Network $\Phi (\cdot; \theta)$ with respect to the prameters
$\theta$ evaluated at the grid nodes. Differentiating the residual
$r_t := \cR_\dx(\Phi(\cdot;\theta_t))$ in time and substituting gives
$$
  \dot r_t \;=\; D\cR_\dx(\Phi(\cdot; \theta_t))\,\partial_\theta \Phi\,
    \dot\theta
  \;=\; -\, D\cR_\dx(\Phi (\cdot; \theta_t))\,\Theta_t\,
    D\cR_\dx(\Phi(\cdot;\theta_t))^\top\, r_t,
$$
where $\Theta_t := \partial_\theta \Phi\, \partial_\theta \Phi^\top$ is
the Gram matrix; it is the neural tangent kernel (NTK) evaluated along the
training trajectory on the grid nodes. Each mode of $r_t$ contracts at
the rate of the corresponding eigenvalue of this composite operator, not
at the rate set by $\kappa_\dx$ alone.

The NTK of standard architectures (fully connected networks with smooth
activations, Fourier features) typically has a sharp eigenvalue decay: a
small number of dominant eigenvalues carry most of the spectral mass, and
the tail decays rapidly. Empirically, these dominant eigenvalues are
aligned with low-frequency, slowly varying modes, while higher-frequency
modes lie in the tail; this is the spectral-bias phenomenon documented in
\cite{rahaman2019spectral,xu2019frequency, lu2026multilevel}. On a coarse grid with
spacing $2\dx$, the residual is resolved only up to frequency
$\sim 1/\dx$, so the modes it can represent are aligned with the dominant
portion of the NTK spectrum, where the eigenvalues of $\Theta_t$ are
large. The composite operator is then much better conditioned on the
coarse problem than $\kappa_{2\dx}$ alone would suggest. Training makes
rapid progress on the modes the network is best equipped to represent.

At the fine level, the additional residual content may still lie within
the dominant part of the NTK spectrum. Even so, the worst conditioning
$\kappa_\dx$ of $D\cR^\top D\cR$ means gradient descent needs more
iterations to converge than at the coarse level. The minibatch point
mentioned above reinforces this: the fine grid has more collocation
points, so variance control calls for larger minibatches, and the
per-iteration cost grows accordingly. Warm-starting each level from the
preceding one thus saves work on two fronts. It reduces both the number
of iterations needed and the cumulative cost per iteration.

This argument is heuristic. It omits subtleties such as the evolution
of $\Theta_t$ in finite-width regimes and the architecture dependence
of the dominant portion of the NTK spectrum. A rigorous treatment is
left to future work. The iteration counts reported in
Figures~\ref{fig: experiment 1 NNs summary} and
\ref{fig: experiment front prop dim 5 summary} support the picture
above. Coarse levels are trained in comparatively few iterations, and
each subsequent refinement benefits both from what has already been
resolved at the coarse level and from the smaller per-level training
cost.

\section{Numerical examples}
\label{sec:numerics}

We present several numerical experiments exemplifying various aspects of the framework.
Section~\ref{sec:1d_eikonal} uses the eikonal equation to illustrate the theoretical findings in a simple case. In subsection \ref{subsec: full grid L2 flow}, for the one-dimensional case, we compute numerically the $\ell^2$-gradient flow associated to finite-difference residuals in the space of grid functions. The latter experiment empirically verifies the condition number estimates of Theorem~\ref{thm:dim-robust error estimate proper case}, and show the benefits of the multi-level warm-start strategy. In subsection \ref{subsec: eikonal NNs}, we move to the NN framework, and show how residual minimization using neural networks can handle high-dimensional settings by solving the eikonal equation in the 4- and 8-dimensional cube.

In section~\ref{subsec:level set numerics}, we apply the same ideas to solve time-dependent Hamilton-Jacobi equations linked to level set methods. In subsection \ref{subsec:front prop}, we solve a front-propagation problem in a 5-dimensional domain, and in subsection \ref{subsec:obstacle}, we solve a non-linear advection problem with a rigid obstacle, in dimension $8$.

Section~\ref{sec:isaacs} treats a Hamilton--Jacobi--Isaacs equation with second-order diffusion arising from a stochastic differential game~\cite{ito-rei-zha-2021}, testing the framework on a non-convex Hamiltonian where the selection of the correct viscosity solution is non-trivial.

\paragraph{Common setup for the neural-network experiments.}
In the experiments below that use a neural network, the solution is approximated by a fully connected feed-forward network $\Phi_\theta$, trained with a stochastic optimizer (Adam, and AdamW in Section~\ref{sec:isaacs}) on a penalized least-squares loss
\begin{equation}\label{eq:generic NN loss}
  L(\theta) = \frac{1}{|\mathcal{X}|}\sum_{x\in\mathcal{X}} \cR(\Phi_\theta;x)^2
  + \frac{\mu_b}{|\mathcal{X}_b|}\sum_{x\in\mathcal{X}_b}\big(\Phi_\theta(x)-g(x)\big)^2 ,
\end{equation}
where $\cR$ is the residual of the relevant monotone scheme and the interior and boundary collocation sets $\mathcal{X}$ and $\mathcal{X}_b$ are drawn by Monte Carlo sampling and refreshed at each iteration, with the interior points kept at least one step $\dx$ from $\partial\Omega$ so that their stencils stay inside the domain. Following the multi-level strategy of Section~\ref{sec:numerical algos}, we minimize the residual over a decreasing sequence of discretization steps: at each level the step is held fixed and the network is trained until a stopping criterion on the loss is met, after which the step is reduced and training is warm-started from the previous level. Because the network is its own interpolant, no grid-transfer operator is needed between levels. When a non-monotone (for instance, second-order) discretization is used, we first pre-train with a monotone first-order scheme so that training begins near the viscosity solution. The architecture, learning rate, penalty weight $\mu_b$, step schedule, sample sizes, and error metrics are specified per experiment.

The numerical experiments were implemented in PyTorch 2.9.1 with CUDA 12.8 and executed on a NVIDIA Tesla V100-PCIE-16GB GPU.

\subsection{Impact of discretization on training efficiency}
\label{sec:1d_eikonal}

\subsubsection{Deterministic gradient descent over grid functions}
\label{subsec: full grid L2 flow}
We show how the discretization affects the convergence speed of gradient descent on the least-squares loss $L(u)$ in \eqref{Lp functional steady state eq}, and that controlling $\|\partial\cR(u)^{-1}\|$ enhances convergence.

For simplicity, we consider the one-dimensional eikonal equation $|\partial_x u| = 1$ in the unit interval $(0,1)$ 
with zero boundary condition, on a uniform grid $\{ j\dx \ : \ j\in \{0, 1, \ldots, N\}\}$, where $\dx = 1/N$.
To test the effect of the properness condition {\bf (H4)}, we also consider the modified equation $|\partial_x u| + u = 1$, which is equivalent to the eikonal equation through the Kruzkhov transform \eqref{eq:eikonal-1}.

To isolate the effect of the discretization from the NN parametrization, we first apply full gradient descent to $L(u; \dx): = \| \cR_\dx (u)\|_2^2$ over grid functions $\R^N$, producing a deterministic sequence that approximates the $\ell^2$-gradient flow of $\cR(u)$.

In the residual $\cR_\dx(u)$ in \eqref{residual function def}, we use the Lax-Friedrichs numerical Hamiltonian given by
\begin{equation}
\label{LxF Hamiltonian numerics 1d}
\cH_{\LF} \left( u_j, \dfrac{u_j - u_{j-1}}{\dx}, \dfrac{u_j - u_{j+1}}{\dx} \right) =
\left| \dfrac{u_{j+1} - u_{j-1}}{2 \dx} \right| - \alpha \dfrac{u_{j+1} + u_{j-1} - 2 u_j}{2\dx} + \lambda u_j - 1,
\end{equation}
for $j\in \{1, 2, \ldots, N-1\}$  with $\alpha = 2$, and null-boundary data $u_0 = u_N = 0$.
As discussed in Section~\ref{sec:monotone-hamiltonians}, $\cH_{\LF}$ satisfies \textbf{(H1)}, \textbf{(H2)} and {\bf (H4)} with $\subalpha = 1/2$. If $\lambda>0$, then {\bf (H3)} holds as well.
In view of Definition~\ref{def: graph distance constants}, the relevant constants are the stencil-sum bound $\overline{s} \leq 1$ and the graph distance to the boundary $d_\infty \leq \frac{N+1}{2}$, with $\dxmax = \dx_{min} = 1/N$.

By Lemma~\ref{lem:DR-invertible alpha}, $\|A^{-1}\|_\infty \leq \dx\, C_1^{1/\dx} + C_2$. On a fine grid (small $\dx$) the Jacobian $D\cR_\dx(u)$ may therefore be nearly singular at some $u$, where the $\ell^2$-gradient flow stalls. A coarser grid, or $\lambda>0$ (which by Lemma~\ref{lem:DR-invertible lambda} gives the $\dx$-independent bound $\|A^{-1}\|_\infty \leq 1/\lambda$), removes this.

\begin{figure}[htbp]
  \centering
  \begin{subfigure}{0.35\textwidth}
    \centering
    \includegraphics[width=\linewidth]{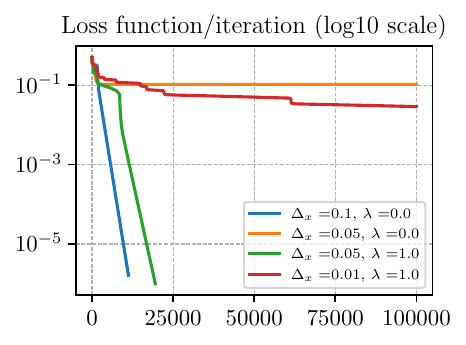}
  \end{subfigure}
  \hspace{1cm}
  \begin{subfigure}{0.35\textwidth}
    \centering
    \includegraphics[width=\linewidth]{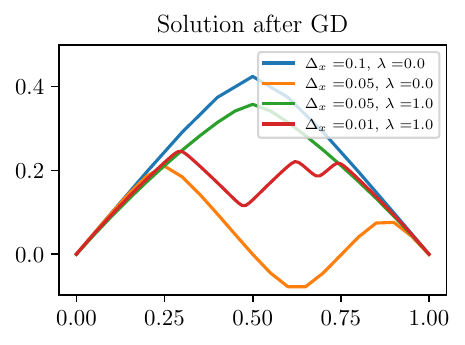}
  \end{subfigure}
  \\
  \begin{subfigure}{0.35\textwidth}
    \centering
    \includegraphics[width=\linewidth]{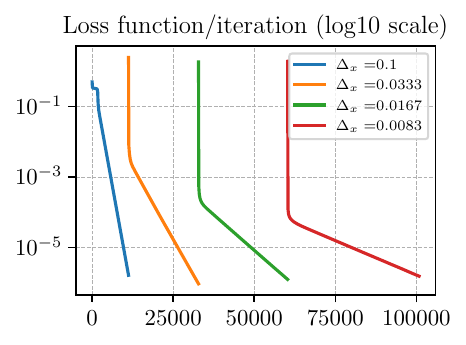}
  \end{subfigure}
    \hspace{1cm}
  \begin{subfigure}{0.35\textwidth}
    \centering
    \includegraphics[width=\linewidth]{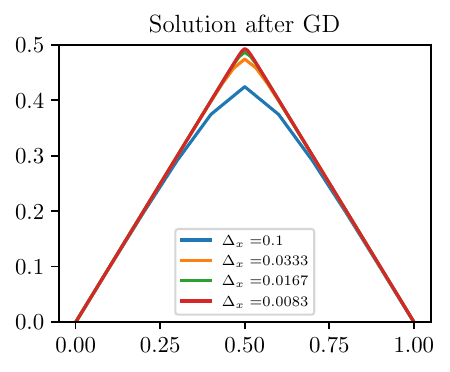}
  \end{subfigure}

  \caption{Full gradient descent on $L(u;\dx)$ in the grid-function space $\R^{N+1}$
($N=1/\dx$): loss evolution (left) and final solution (right). Top: four
independent runs from $u_0=0$. Bottom: a single run with $\lambda=0$ and
progressive grid refinement; the jumps in the loss (bottom-left) mark each
change of $\dx$ and the piecewise-constant interpolation between grids.
  }
  \label{fig:gradient descent}
\end{figure}

Figure~\ref{fig:gradient descent} (top) shows four independent runs from $u_0=0$, with the iteration $u_{n+1}=u_n-10^{-3}\nabla L(u_n;\dx)$ stopped once $\|\cR_\dx(u)\|_\infty<10^{-3}$. Convergence is fast for $(\lambda,\dx)=(0,0.1)$ and $(1,0.05)$ but slow for $(0,0.05)$ and $(1,0.01)$: at fixed $\lambda$, refining the grid degrades convergence, as our theoretical findings predict.

The bottom plots show the multi-level remedy in a single run with $\lambda=0$: starting from $u_0=0$, we minimize $L(u;\dx)$ over a decreasing sequence of $\dx$, warm-starting each level from the previous one (extended by piecewise-constant interpolation); the jumps in the loss mark the grid changes. Convergence is now fast even at small $\dx$, since the warm start keeps the iterate out of the flat regions of the landscape. As $\dx\to0$ this approximates the viscosity solution of $|\partial_x u|=1$.

\begin{figure}[htbp]
\centering
  \begin{subfigure}{0.35\textwidth}
    \centering
    \includegraphics[width=\linewidth]{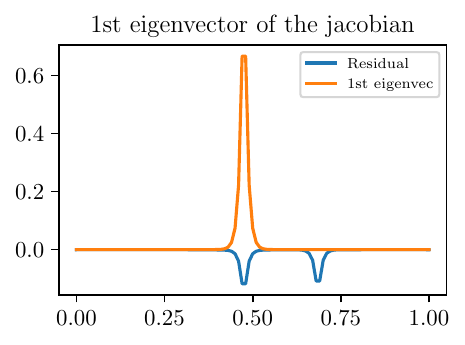}
  \end{subfigure}
  \hspace{1cm}
  \begin{subfigure}{0.35\textwidth}
    \centering
    \includegraphics[width=\linewidth]{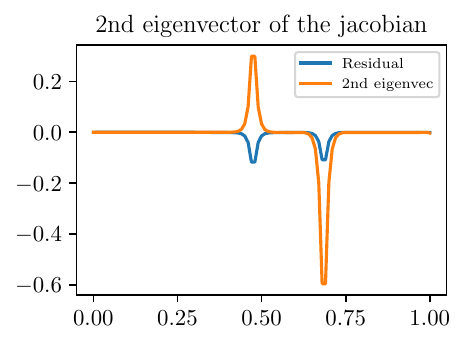}
  \end{subfigure}

  \caption{For the last run in Figure~\ref{fig:gradient descent} (top), with
$\dx=0.01$ and $\lambda=1$: (blue) the residual $\cR_\dx(u)$ after $10^5$
gradient-descent iterations; (orange) the eigenvectors of the Jacobian
$D\cR_\dx(u)$ for its smallest eigenvalue, $\approx 0.1$.}
  \label{fig:eigenvectors Algo 1}
\end{figure}

Figure~\ref{fig:eigenvectors Algo 1} examines the slowest run, $(\lambda,\dx)=(1,0.01)$. After $10^5$ iterations the iterate has not converged because $D\cR_\dx(u)$ is nearly singular: its smallest eigenvalue is $\approx 0.1$, and the residual (blue) is essentially a combination of the corresponding eigenvectors (orange). Their support, of width $\sim\dx$, concentrates where the computed profile has the wrong convexity---curving opposite to the viscosity solution and forming a spurious kink.

\subsubsection{Implementation with neural networks and uniform sampling}
\label{subsec: eikonal NNs}

Let us now replace the space of grid functions by a parametric family of NNs. Since we are no longer tied to a fixed grid, we can raise the dimension.

\paragraph{The equation.}
We consider the eikonal equation with zero boundary condition in the $d$-dimensional cube $\Omega = (-1,1)^d$ with $d= 4$ and $8$.

\paragraph{Monotone discretization.}
As for the numerical Hamiltonian, we consider the Upwind discretization given by
\begin{equation}
\label{Upwind scheme def}
\cH \left( u(x), D^\pm u(x) \right)  =  
\sqrt{\sum_{k=1}^d \Big(\max \left( D_k^+ u(x), D_k^- u(x) ,0 \right) \Big)^2}  - 1 + u(x),
\end{equation}
where $D_k^+ u(x)$ and $D_k^- u(x)$ are first- or second-order approximations of the derivatives in the directions $e_k$, where $\{e_k\}_{k=1,\ldots, d}$ is the canonical basis of $\R^d$.
In particular, we will compare the one-sided first-order approximations defined in \eqref{eq:one-sided-differences} 
and the second-order approximation given by
\begin{equation}
\label{2nd order grad eikonal numerics}
D_k^\pm u(x) := \dfrac{3u(x) - 4u( x \pm \dx\, e_k) + u(x\pm 2 \dx\, e_k)}{2\dx}, \qquad \text{for} \quad k = 1,\ldots, d.
\end{equation}
It is easy to verify (see Section~\ref{sec:monotone-hamiltonians}) that, with the first-order approximations \eqref{eq:one-sided-differences}, the numerical Hamiltonian $\cH$ satisfies {\bf (H1)}, {\bf (H2)}, and {\bf (H3)}. 
The numerical scheme $\cH$ in \eqref{Upwind scheme def} is consistent with the PDE $\| \nabla v\|_2 + v = 1$. But, since we are interested in solving the eikonal equation, we apply the inverse Kruzhkov transform \eqref{eq:eikonal-1} to the NN approximation $\Phi_\theta(x) \approx v(x)$, to obtain an approximation of the solution for the PDE $\|\nabla u\|_2 = 1$ as
$$
\hat u_\theta(x) = \log\left( \frac{1}{1 -  \Phi_\theta(x)} \right).
$$

\begin{figure}[htbp]
  \centering
  \begin{subfigure}{0.35\textwidth}
    \centering
    \includegraphics[width=\linewidth]{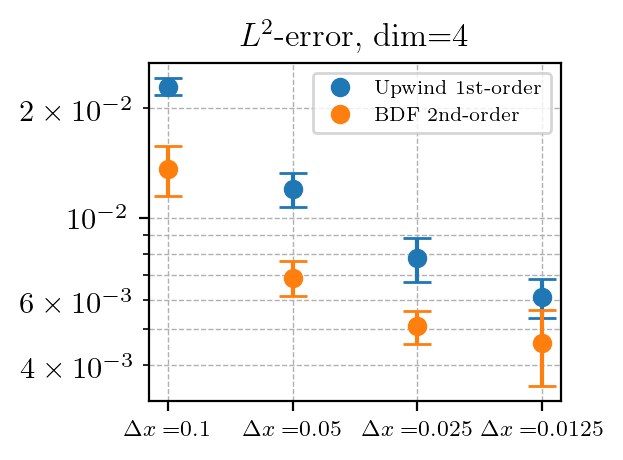}
    \includegraphics[width=\linewidth]{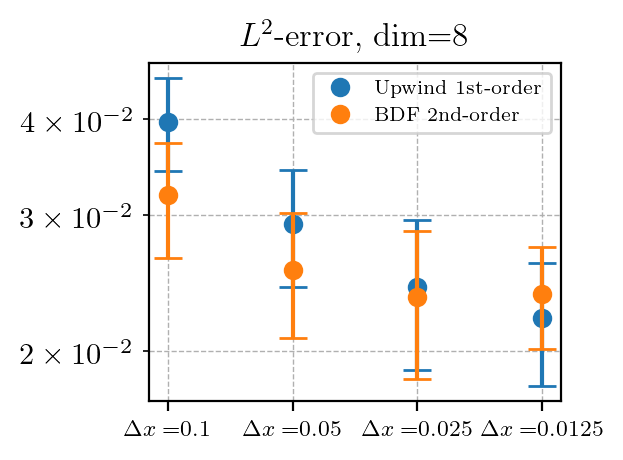}
  \end{subfigure}
  \begin{subfigure}{0.35\textwidth}
    \centering
    \includegraphics[width=\linewidth]{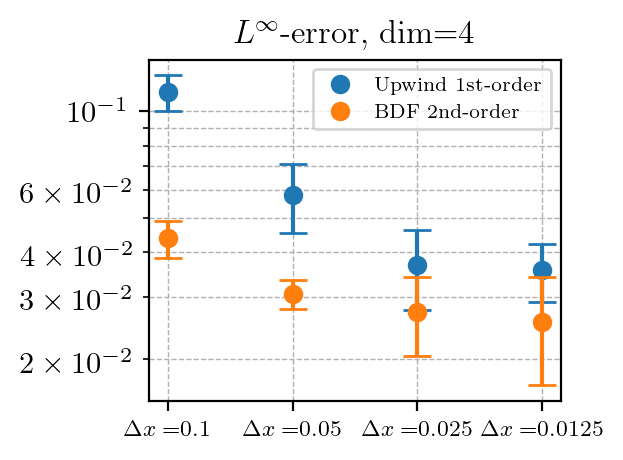}
    \includegraphics[width=\linewidth]{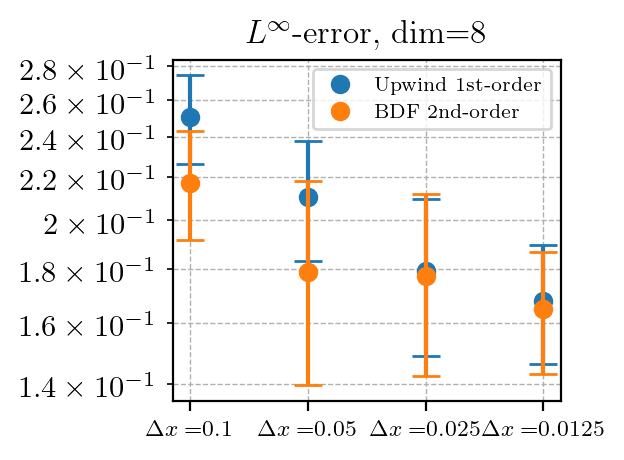}
  \end{subfigure}
  \caption{Errors for the eikonal equation in the $d$-dimensional cube
(Section~\ref{subsec: eikonal NNs}), $d=4$ (top) and $d=8$ (bottom), from
minimizing the upwind residual~\eqref{Upwind scheme def} with first-order
(blue) and second-order (orange) gradient approximations, each over $10$
independent runs. Both are pre-trained with the first-order scheme at
$\dx=0.1$ to ensure convergence to the viscosity solution.}
  \label{fig: experiment 1 NNs summary}
\end{figure}

\paragraph{Network, loss, and training.}
We minimize the penalized loss~\eqref{eq:generic NN loss} (with $\mu_b=10$) over the decreasing schedule
\begin{equation}
\label{training schedule exp 1}
(\dx_1, \, \dx_2, \, \dx_3, \, \dx_4) = (0.1,\,  0.05, \, 0.025, \, 0.0125).
\end{equation}
The network $\Phi_\theta$ has $4$ hidden layers of $120$ neurons ($p=58801$ parameters); at each level we run at most $10^4$ iterations, stopping once the loss stays below $5\cdot 10^{-4}$ for $20$ consecutive iterations, and draw $(N,N_b)=(12000,2000)$ interior/boundary points for $d=4$ and $(24000,6000)$ for $d=8$. Since the upwind scheme with the second-order gradient~\eqref{2nd order grad eikonal numerics} is not monotone, gradient descent need not converge to the viscosity solution; we therefore pre-train with the monotone first-order scheme~\eqref{eq:one-sided-differences} at $\dx=0.1$ before switching to the second-order one.
\paragraph{Reference solution and results.}
The exact solution is the distance function to the boundary.
Figure~\ref{fig: experiment 1 NNs summary} reports the $L^2$ and $L^\infty$ errors against this exact solution---means and standard deviations over $10$ runs, measured at the end of each refinement level in~\eqref{training schedule exp 1}---for both the first- and second-order gradient approximations.

\subsection{Level set method}
\label{subsec:level set numerics}

To demonstrate the capability of our framework to handle high-dimensional geometry without spatial meshes, we solve time-dependent front propagation problems in $\R^+\times \R^5$. First, we consider the propagation of a front that evolves along its outer normal. Then, we consider a non-linear advection problem with a rigid obstacle.

\subsubsection{Front propagation}
\label{subsec:front prop}

\begin{figure}[htbp]
  \centering
  \begin{subfigure}{0.25\textwidth}
    \centering
    \includegraphics[width=\linewidth]{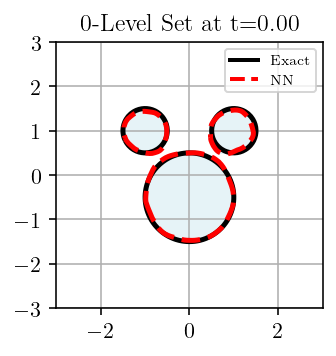}
  \end{subfigure}
  \hfill
  \begin{subfigure}{0.25\textwidth}
    \centering
    \includegraphics[width=\linewidth]{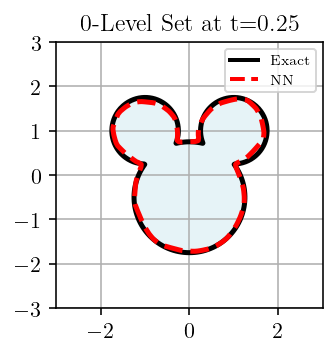}
  \end{subfigure}
  \hfill
  \begin{subfigure}{0.25\textwidth}
    \centering
    \includegraphics[width=\linewidth]{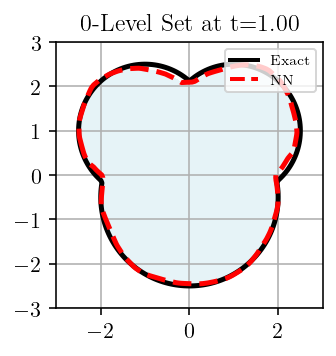}
  \end{subfigure}

  \caption{Zero-level set of the solution to \eqref{time-evol eikonal equation example} at different times. Each image corresponds to a two-dimensional slice of the space-domain $(-3,3)^5$. The NN approximation is trained using the first-order upwind discretization \eqref{implicit first-order upwind}, with decreasing values of $\dx$ and $\dt$ as shown in Figure~\ref{fig: experiment front prop dim 5 summary}.} 
  \label{fig: level-set plots 5d}
\end{figure}

\paragraph{The equation.}
We consider an initial set formed by the union of three disjoint balls
$S_0 := B(x_1, r_1)\cup B(x_2, r_2)\cup B(x_3, r_3)\subset \R^5$, with centers $x_1 = (-1, 1,0,0,0)$, $x_2 = (1,1,0,0,0)$ and $x_3 = (0,-\frac{1}{2},0,0,0)$, and radii $r_1=r_2 = \frac{1}{2}$ and $r_3 = 1$. The evolution of the boundary of $S_0$ in the normal outer direction to the boundary at speed 1 can be described, using the level-set method \cite{osher1988fronts}\cite{sethian1999level}\cite{osher2003level}, as the $0$-level set of the viscosity solution to the HJ equation
\begin{equation}
\label{time-evol eikonal equation example}
\begin{cases}
    \partial_t u + \|\nabla u\|_2 = 0 & (t,x)\in \R^+\times \R^5 \\
    u(0,x) = u_0(x) & x\in \R^5,
\end{cases}
\end{equation}
where the initial condition $u_0$ can be any Lipschitz function such that $S_0 = \{ u_0(x) < 0 \}$. For our experiments, we take
$$
u_0(x) := \min \left\{ \| x-x_1\| - r_1, \, \| x-x_2\| - r_2, \, \| x-x_3\| - r_3,\, \frac{1}{2}
\right\}.
$$
Since our framework considers only bounded domains, we restrict the solution to the hyper-cube $(0,T)\times \Omega:=(0,1)\times (-3,3)^5$. It is easy to prove that the solution in $(0,T)\times \Omega$ is uniquely determined by the initial condition, and no boundary condition needs to be prescribed. The reason is that the characteristic lines emanating from outside $\Omega$ never cross the boundary of this set. Figure~\ref{fig: level-set plots 5d} shows the computed zero-level set at three successive times.

\paragraph{Monotone discretization.}
The residual at the interior nodes uses the implicit first-order upwind scheme
\begin{equation}
\label{implicit first-order upwind}
\cR (u; t,x) := \frac{u(t,x) - u(t-\dt)}{\dt} + \sqrt{\sum_{k=1}^d \Big(\max \left( D_k^+ u(x), D_k^- u(x) ,0 \right) \Big)^2} ,
\end{equation}
where $D_k^+ u(x)$ and $D_k^- u(x)$ are the first-order one-sided differences \eqref{eq:one-sided-differences}; it satisfies {\bf (H1)}, {\bf (H2)} and {\bf (H3)} with $\ml = 1/\dt$.

\paragraph{Network, loss, and training.}
We use a fully connected network with $3$ hidden layers of $128$ neurons, and minimize the loss \eqref{eq:generic NN loss} with penalty weight $\mu_b=10$ and initial data $g=u_0$, drawing at each iteration $12000$ interior points uniformly in $(0,T)\times \Omega$ and $3000$ points uniformly in $\Omega$.

\begin{figure}[htbp]
  \centering
  \begin{subfigure}{0.32\textwidth}
    \centering
    \includegraphics[width=\linewidth]{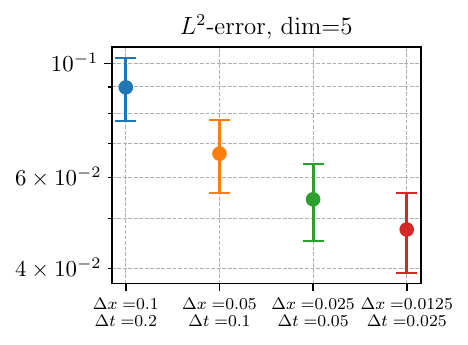}
  \end{subfigure}
  \hfill
  \begin{subfigure}{0.32\textwidth}
    \centering
    \includegraphics[width=\linewidth]{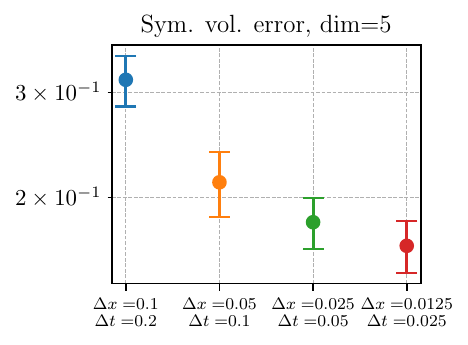}
  \end{subfigure}
  \hfill
  \begin{subfigure}{0.32\textwidth}
    \centering
    \includegraphics[width=\linewidth]{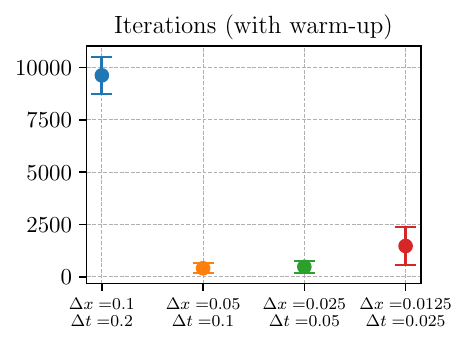}
  \end{subfigure}

  \caption{Front evolution in $\R^5$ (Section~\ref{subsec:front prop}). The NN
approximates~\eqref{time-evol eikonal equation example} by minimizing the
first-order upwind residual~\eqref{implicit first-order upwind}. For the
$L^2$-error~\eqref{L2 error narrow band} and symmetric volume
error~\eqref{sym volume error}, dots are means over $10$ runs and intervals
standard deviations; the four columns mark the end of each $(\dt,\dx)$
refinement round.}
  \label{fig: experiment front prop dim 5 summary}
\end{figure}

Concerning the discretization parameters, we use a 4-round training schedule with decreasing values of $\dx$ and $\dt$ as follows:
\begin{equation}
\label{training schedule front prop}
\begin{array}{rll}
(\dx_1, \dx_2, \dx_3, \dx_4) & = & (0.1,\ 0.05,\ 0.025,\ 0.0125),\\
(\dt_1, \dt_2, \dt_3, \dt_4) & = & (0.2,\ 0.10,\ 0.050,\ 0.0250).
\end{array}
\end{equation}
At each level we train until the loss stays below $10^{-4}$ for $50$ consecutive iterations, then refine to the next $(\dx, \dt)$.
\paragraph{Reference solution and results.}
The same experiment is repeated ten times to test the robustness of our method. In Figure~\ref{fig: experiment front prop dim 5 summary} (right), for each of the four training rounds, we see the number of iterations that it took (on average) to meet the stopping criterion. We observe that training the NN with larger values of $(\dx, \dt)$ improves convergence speed when we switch to smaller discretization steps.

To test the accuracy of our approximation, we use the exact solution as ground truth, which is given by the Hopf-Lax formula as
$$
u(t, x) = \min \left\{ u_0(y) \, :\,  y\in \overline{B(x, t)} \right\}.
$$
Since we are interested in approximating the 0-level set of the solution, we compute the error in an $\varepsilon$-neighborhood of this level set:
$$
\mathcal{N}_{\varepsilon} (t) := \{ x\in \R^5 \ : \quad | u(t,x)| < \varepsilon \}.
$$
In Figure~\ref{fig: experiment front prop dim 5 summary}, we see a summary of our experiment (which is repeated ten times). The error is computed at the end of each training round, right before the parameters $\dx$ and $\dt$ are decreased (following \eqref{training schedule front prop}). 

In Figure~\ref{fig: experiment front prop dim 5 summary} (left), we report the $L^2$-error:
\begin{equation}
\label{L2 error narrow band}
\| \hat{u} - u\|_\varepsilon :=\left(\int_0^T \int_{\mathcal{N}_\varepsilon (t)} | \hat{u} (t,x) - u(t,x)|^2 dx dt \right)^{1/2}.
\end{equation}
In Figure~\ref{fig: experiment front prop dim 5 summary} (middle), we report the symmetric volume error, which measures the volume of the region in which the ground truth $u(t,x)$ and the approximation $\hat{u}(t,x)$ have opposite sign:
\begin{equation}
\label{sym volume error}
\mathcal{E} (\hat{u}, u) := \int_0^T \dfrac{| \{x\in \mathcal{N}_\varepsilon (t) \ : \ \text{sign}(\hat{u}(t,x)) \neq \text{sign}(u(t,x)) \} |}{|\mathcal{N}_\varepsilon (t)|} dt .
\end{equation}
The $L^2$ error \eqref{L2 error narrow band} and the symmetric volume error \eqref{sym volume error} are approximated using the Monte Carlo method.

\subsubsection{Non-linear advection with an obstacle}
\label{subsec:obstacle}
\paragraph{The equation.}
Here, we solve a time-dependent non-linear advection problem in $\R^d$, with $d=8$, in the presence of a rigid obstacle. The problem consists in finding the viscosity solution to
\begin{equation}
\begin{cases}
    \min \left( \partial_t u
    + (f\cdot \nabla u)_+, \ u(t,x) - \Psi(x) \right) = 0, & (t,x) \in (0,T) \times \R^d, \\
    u(0, x) = g(x), & x\in \Omega,
    \label{eq:obstacle_problem}
\end{cases}
\end{equation}
where the drift is constant and given by $f:= (1,\dots,1) \in \R^d$, the rigid obstacle is a ball of radius $1/2$ centered at the origin, represented by $\{\Psi(x)\geq 0\}$, where
\begin{equation}
    \Psi(x) = \min\left( \frac{1}{2} - \| x\|, \
     \frac{1}{4}\right),
\end{equation}
and the initial set is the intersection of the complement of the obstacle with the unit ball centered at $a_0 = \left(-1,\ldots,-1\right)$, represented through the initial condition as $\{g(x) \leq 0\}$, where
$$
g(x):= \min \left(\max \left( \| x - a_0 \|_2-1  , \, \Psi(x) \right), \ \frac{1}{4}\right).
$$
Note that we have truncated both the obstacle $\Psi(x)$ and the initial condition $g(x)$ at the level $1/4$. While this does not change the zero-level set of the viscosity solution of \eqref{eq:obstacle_problem}, it is convenient to do so since it ensures that the solution is constant outside a compact set of radius depending on $T$.
In this way, it is sufficient to approximate the solution in this compact set, and one can omit the boundary condition term in the loss functional, which is already integrated in the definition of the initial condition and the obstacle.

\paragraph{Monotone discretization.}
The residual combines an implicit time discretization with an upwind discretization of the advection term,
$$
\cR (u(t,x)) = \min \left( \dfrac{u(t,x) - u(t-\dt, x)}{\dt} + \left( \sum_{k=1}^d D_k^- u(t,x) \right)_+, \  u(t,x) - \Psi (x) \right).
$$
Since the drift is constant ($f \equiv (1,\ldots,1)$), only the left-sided differences $D_k^- u$ \eqref{eq:one-sided-differences} enter the upwind term. We also use a second-order variant, replacing the first-order differences by the BDF2 approximation in time
$$
\partial_t u(t,x) \approx \dfrac{3u(t,x) - 4u(t-\dt, x) + u(t-2 \dt,x)}{2 \dt},
$$
(chosen for its stability properties~\cite{bok-pic-rei-21}),
and BDF2 in space by \eqref{2nd order grad eikonal numerics}.

\paragraph{Network, loss, and training.}
We use a fully connected feedforward neural network with $4$ hidden layers of $120$ neurons ($58801$ parameters), and minimize the loss \eqref{eq:generic NN loss} with penalty weight $\mu_b=10$ and the initial data $g$ above. At each iteration the interior batch $\mathcal{X}\subset (0,T)\times \R^d$ has $2\cdot 10^5$ points---the time sampled uniformly in $(0,T)$, $T=2$, and the space variables from $\mathcal{N}(0,\sigma^2 I_d)$ with $\sigma=1$---while the initial batch $\mathcal{X}_0\subset\{0\}\times\R^d$ is sampled the same way. For each scheme we pre-train with the coarse first-order discretization ($\dt=0.2$, $\dx=0.1$) for $10^4$ iterations, then run $10^5$ iterations at the finer steps ($\dt=0.0125$, $\dx=0.00625$).

\begin{figure}[!htbp]
  \centering
  \hspace*{-5ex}
  \begin{tabular}{@{} r @{\hspace{0.5em}} ccc @{}}
    & $t=0$ & $t=1$ & $t=2$ \\
    \rotatebox{90}{\makebox[2.5em]{\hspace{20ex} 1st order}} &
    \includegraphics[width=0.26\textwidth]{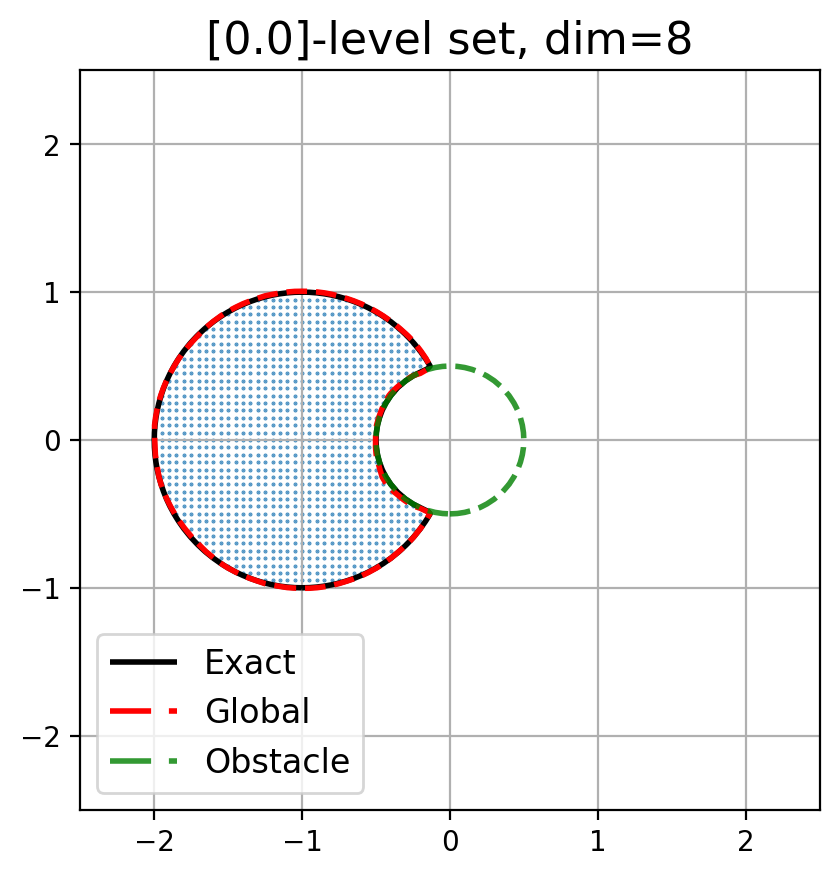} &
    \includegraphics[width=0.26\textwidth]{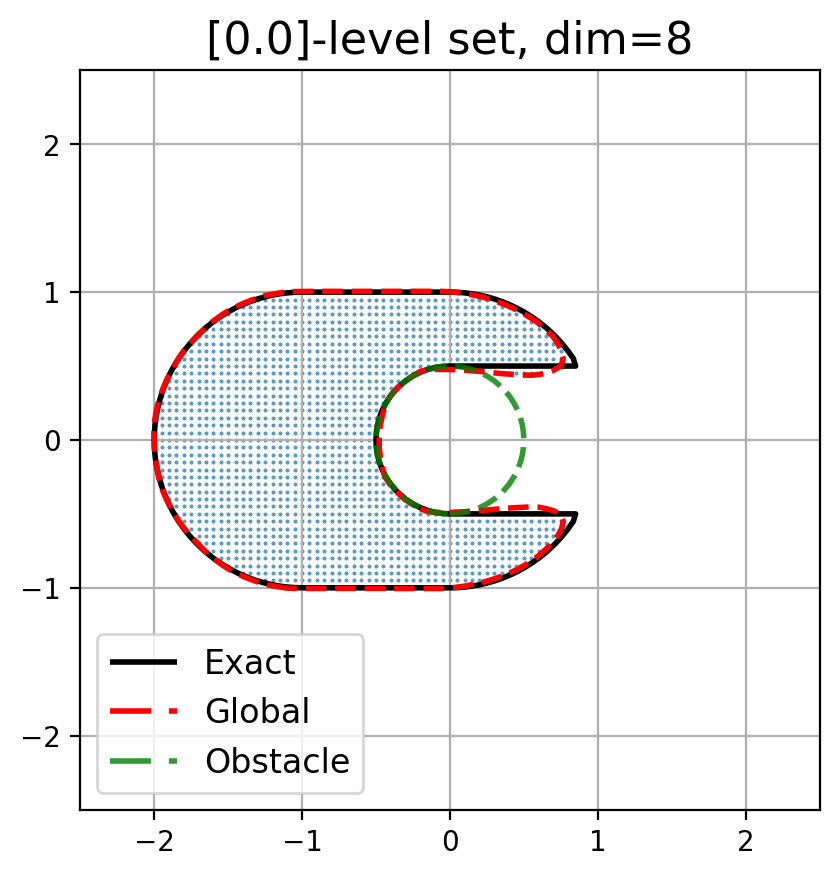} &
    \includegraphics[width=0.26\textwidth]{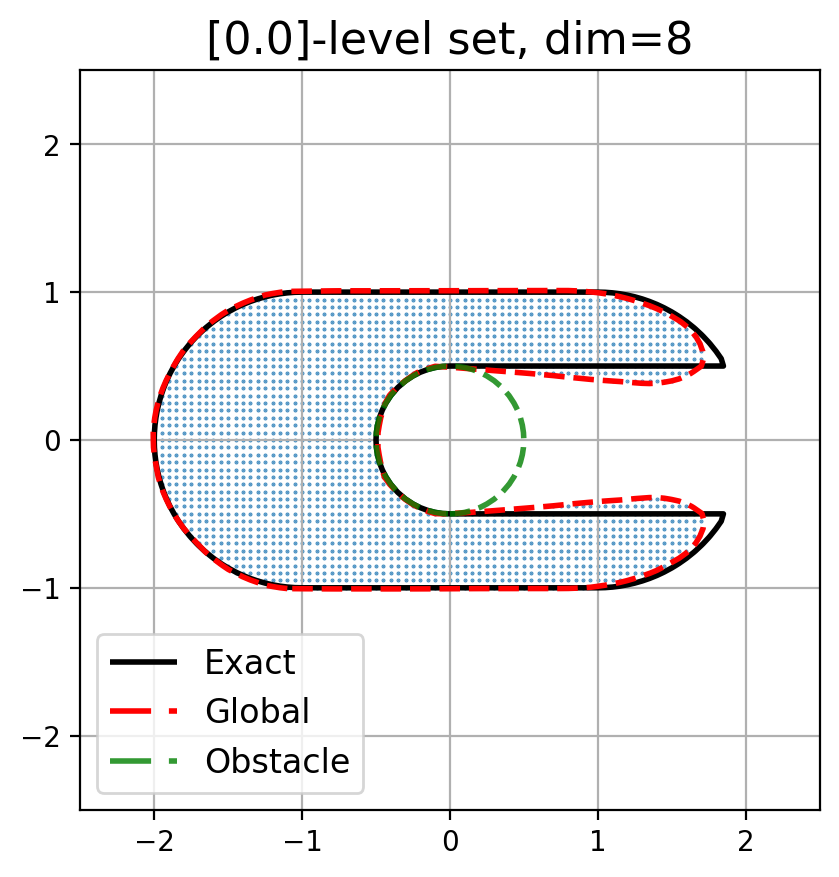} \\[1ex]
    \rotatebox{90}{\makebox[2.5em]{\hspace{20ex} 2nd order}} &
    \includegraphics[width=0.26\textwidth]{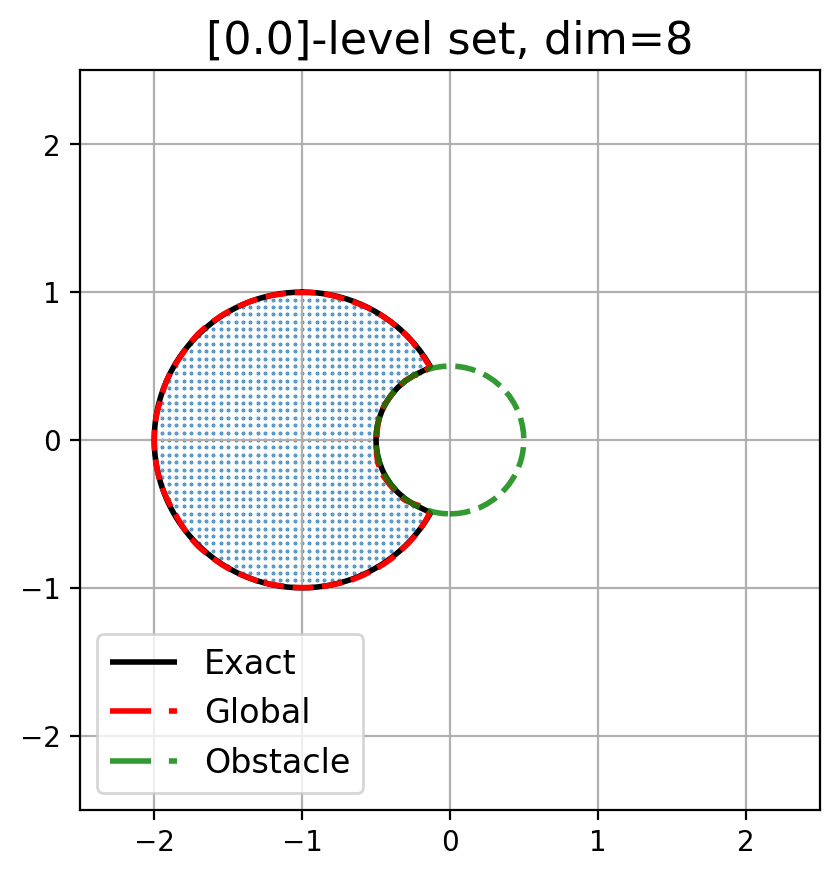} &
    \includegraphics[width=0.26\textwidth]{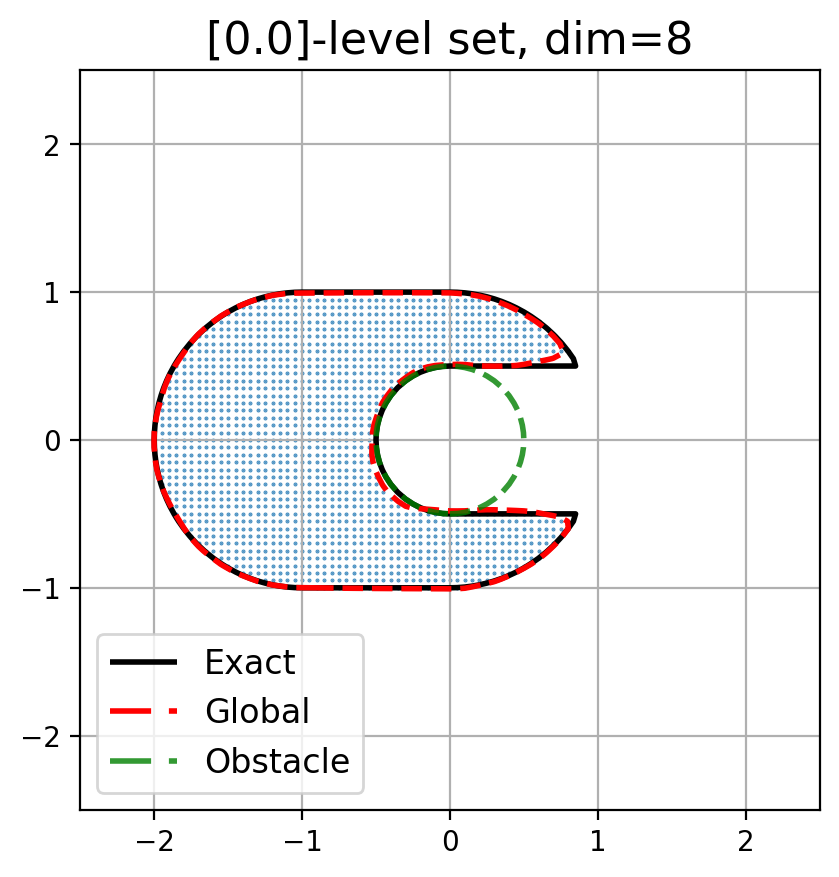} &
    \includegraphics[width=0.26\textwidth]{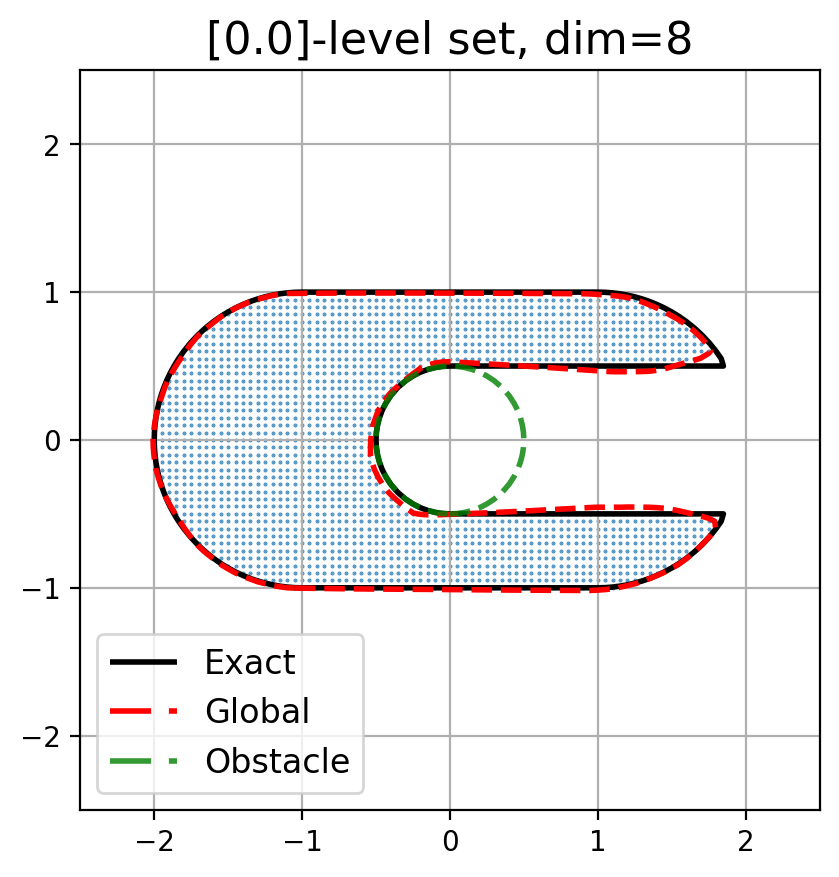}
  \end{tabular}
  \caption{\label{fig:obstacle_8d} 
  Zero-level set at $t=0,1,2$ (central 2D slice of the 8D domain): first-order (top) vs.\ second-order (bottom) discretization, both pre-trained first-order at larger $\dx$ for faster convergence.}
\end{figure}

\paragraph{Reference solution and results.}
Figure~\ref{fig:obstacle_8d} compares the trained network with the exact solution: the zero-level set at $t=0,1,2$ on a central cross section of the $8$-dimensional domain, rotated so that the drift $f=(1,\ldots,1)$ is aligned with the $x$-axis. The top row uses the first-order scheme, the bottom row the second-order BDF2 variant.

\subsection{A viscous Hamilton--Jacobi--Isaacs equation}
\label{sec:isaacs}

Our final example is a second-order Hamilton--Jacobi--Isaacs equation
arising from a stochastic differential game: the $d$-dimensional extension
of the Zermelo navigation problem of Ito, Reisinger and
Zhang~\cite{ito-rei-zha-2021}, originally posed in dimension two. Its
Hamiltonian is neither convex nor concave, so the selection of the correct
viscosity solution is nontrivial, while the diffusion term makes the
discrete scheme uniformly elliptic in the sense of {\bf (H4)}.

\paragraph{The equation.}
Writing $x=(x_1,x')\in\R^d$ with $x'=(x_2,\dots,x_d)$, the value function
$u$ solves
\begin{equation}
\label{eq:HJI-dimd}
  -\tfrac12\sigma_1^2\,\partial_{x_1x_1}u
  -\tfrac12\sigma_2^2\sum_{i=2}^{d}\partial_{x_ix_i}u
  - v_c(x)\,\partial_{x_1}u
  + v_s\,\|\nabla u\|_2
  - \kappa\,\|\Sigma\nabla u\|_1^{\ast}
  = 1,
  \qquad x\in\Omega,
\end{equation}
where the anisotropic, non-smooth first-order term $\|\Sigma\nabla u\|_1^{\ast}$ is given by
\begin{equation}
\label{eq:HJI-l1star}
  \|\Sigma\nabla u\|_1^{\ast}
  := \sigma_1\,|\partial_{x_1}u|
   + \sigma_2\,\big\|(\partial_{x_2}u,\dots,\partial_{x_d}u)\big\|_2 .
\end{equation}
The domain is the spherical shell
$\Omega := \{\,x\in\R^d : r_0<\|x\|_2<r_1\,\}$ with $r_0=\tfrac12$,
$r_1=\sqrt2$, and the boundary conditions are of Dirichlet type: $u=0$ on
$\{\|x\|_2=r_0\}$ (the target to be reached) and $u=1$ on
$\{\|x\|_2=r_1\}$ (the set to be avoided). The drift acts only along the
$x_1$-axis through the wind field
\[
  v_c(x) = 1 - a\,\sin\!\Big(\pi\,\frac{\|x\|_2^2-r_0^2}{r_1^2-r_0^2}\Big),
  \qquad x\in\Omega .
\]
The equation couples an anisotropic diffusion ($\sigma_1$ along $x_1$,
$\sigma_2$ in the remaining directions), the convex term $v_s\|\nabla u\|_2$,
and the concave term $-\kappa\|\Sigma\nabla u\|_1^{\ast}$. We take
$\sigma_1=0.5$, $\sigma_2=0.2$, $v_s=1.0$, $\kappa=0.1$, $a=0.2$, no
zeroth-order term, and $d=8$.

\paragraph{Reference solution.}
The problem is invariant under rotations of the last $d-1$ coordinates: all
data depend on $x$ only through $(x_1,\rho)$ with $\rho:=\|x'\|_2$, so
$u(x)=U(x_1,\rho)$ and~\eqref{eq:HJI-dimd} reduces to a two-dimensional
equation in $(x_1,\rho)$, the diffusion in the last $d-1$ variables becoming
$\tfrac12\sigma_2^2\big(\partial_{\rho\rho}U+\tfrac{d-2}{\rho}\,\partial_\rho U\big)$;
the dimension thus enters only through the radial coefficient
$\tfrac{d-2}{\rho}$. We obtain a reference solution by solving this reduced
problem with a second order finite-difference scheme
by a semi-smooth Newton's method---on a $800\times800$ grid of $[-2,2]^2$. 

\paragraph{Monotone discretization.}
We discretize~\eqref{eq:HJI-dimd} with the first-order monotone numerical
Hamiltonian of Section~\ref{sec:monotone-hamiltonians}. With the one-sided
differences $D_k^{\pm}u_i=(u_i-u_{i\pm e_k})/\dx_k$, the second derivatives
use the central stencil
$\delta_k^2u_i=(u_{i+e_k}-2u_i+u_{i-e_k})/\dx_k^2=-(D_k^{+}u_i+D_k^{-}u_i)/\dx_k$,
while the first-order terms are upwinded:
\begin{equation}
\label{eq:HJI-Hnum}
\begin{aligned}
  \cH(x_i,u_i,D^{\pm}u_i)
  ={}& -\tfrac12\sigma_1^2\,\delta_1^2u_i
       -\tfrac12\sigma_2^2\sum_{k=2}^{d}\delta_k^2u_i
       \;+\; v_c(x_i)\,D_1^{+}u_i \\[2pt]
     & +\, v_s\sqrt{\sum_{k=1}^{d}\big[\max(D_k^{+}u_i,D_k^{-}u_i,0)\big]^2}
       \;+\;\kappa\,H_1^{\ast} - 1,
\end{aligned}
\end{equation}
with
\[
  H_1^{\ast}
  = \sigma_1\,\min\!\big(D_1^{-}u_i,D_1^{+}u_i\big)
  - \sigma_2\,\Big\|\big(\min(D_k^{-}u_i,D_k^{+}u_i)\big)_{k=2}^{d}\Big\|_2 .
\]
Since $v_c\ge0$, the term $v_c\,D_1^{+}u_i$ is the monotone upwind
approximation of $-v_c\,\partial_{x_1}u$; the eikonal term is the standard monotone approximation of $\|\nabla u\|_2$; and $H_1^{\ast}$ approximates
$-\|\Sigma\nabla u\|_1^{\ast}$ through
$-|\partial_{x_k}u|=\min(\partial_{x_k}u,-\partial_{x_k}u)\approx
\min(D_k^{-}u_i,D_k^{+}u_i)$. The diffusion supplies the uniform ellipticity
that makes the scheme satisfy {\bf (H4)} once $\dx$ is small enough, so the
theory of Sections~\ref{sec: main results stationary case}--\ref{sec: stability}
applies.

\begin{figure}[htbp]
  \centering
  \begin{tabular}{ccc}
    \includegraphics[width=0.35\linewidth,height=4.3cm]{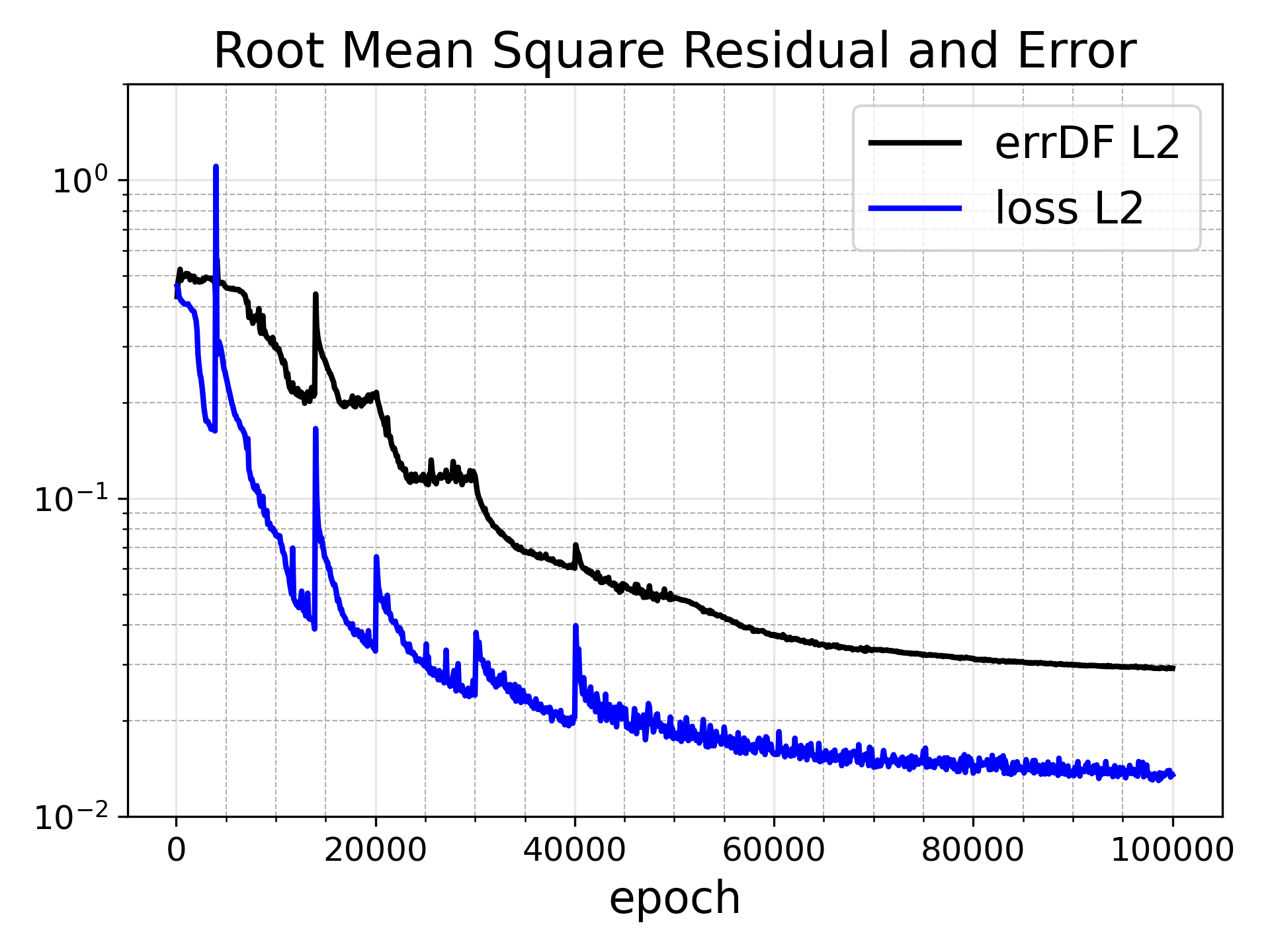} &
    \includegraphics[width=0.30\linewidth,height=4.5cm]{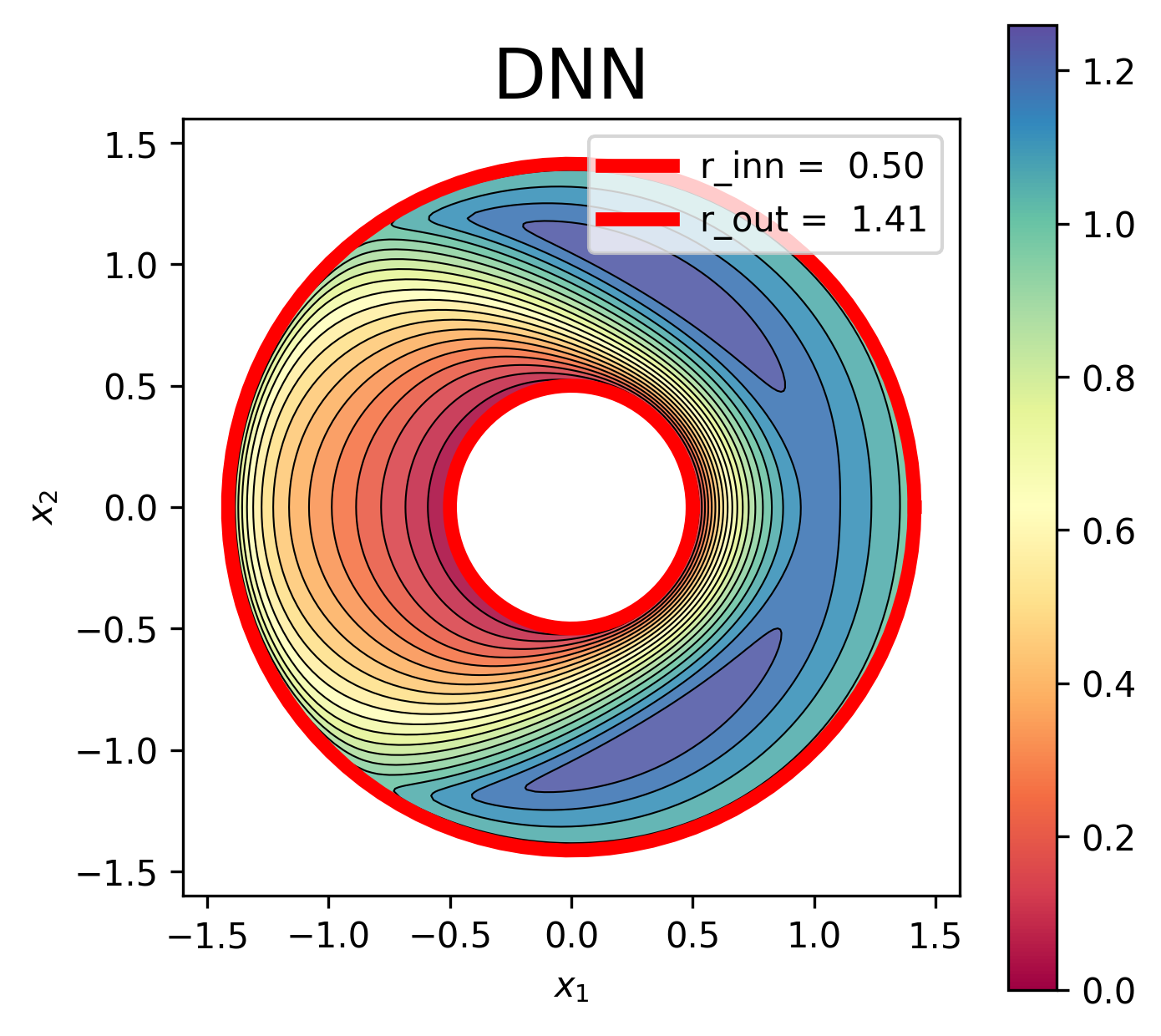} &
    \includegraphics[width=0.30\linewidth,height=4.5cm]{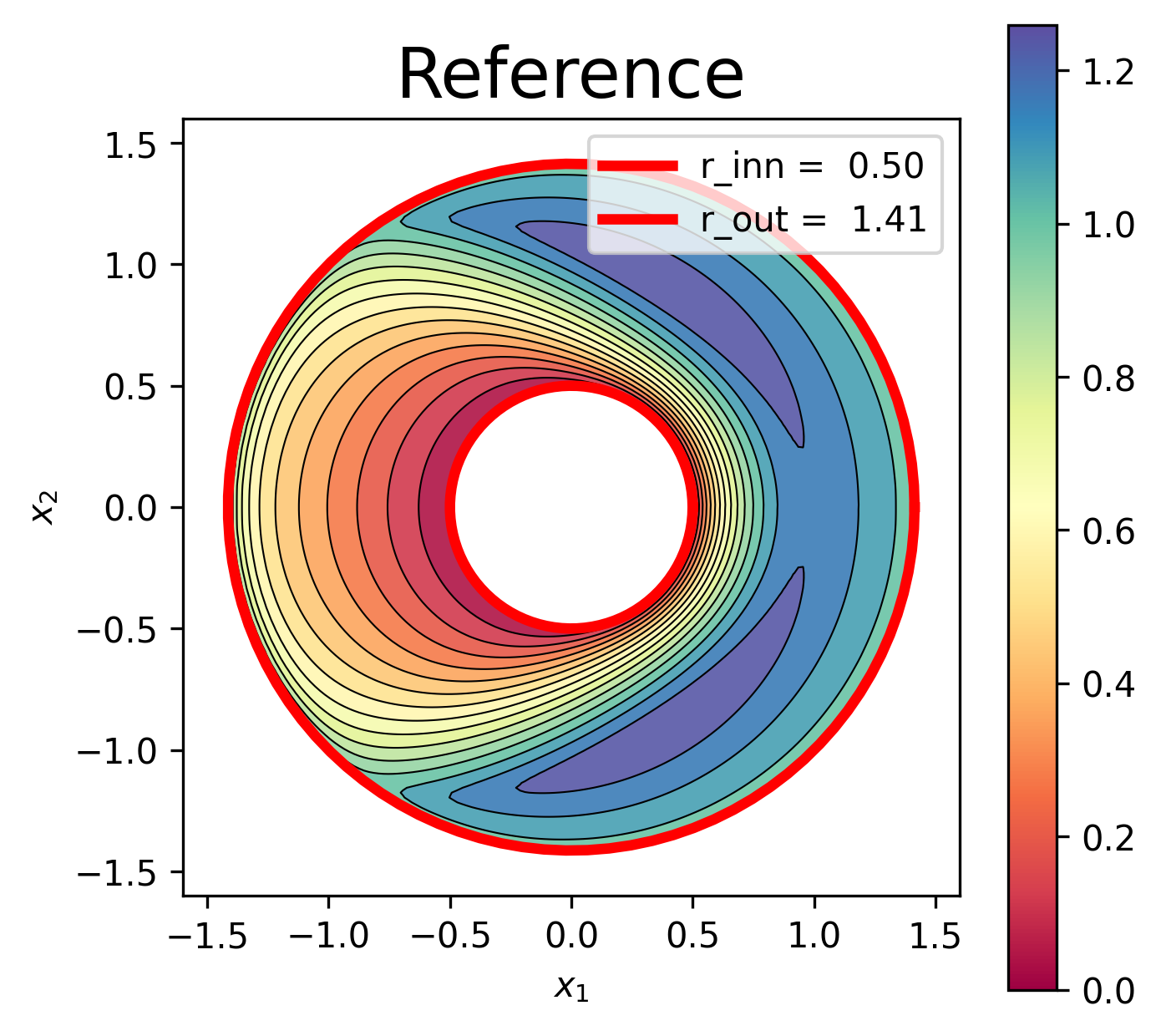}
  \end{tabular}
  \caption{\label{fig:HJI-conv-dim8}
  Convergence history for the $d=8$ HJI experiment: (left) residual loss - in blue - and
  $\ell^2$ error against the reference - in black - versus Adam iterations; the
  jumps mark the successive grid refinements of the multi-level
  schedule (see~\eqref{eq:HJI-schedule}); level sets of the DNN solution (center) and reference (right).}
\end{figure}

\paragraph{Network, loss, and training.}
We approximate $u$ by a fully connected network $\Phi_\theta$ with $5$
hidden layers of $80$ neurons and $\tanh$ activation, trained with the
AdamW optimizer on the penalized loss~\eqref{eq:generic NN loss} with $\mu_b=50$ and residual $\cR(\Phi_\theta;x)=\cH\big(x,\Phi_\theta(x),D^{\pm}\Phi_\theta(x)\big)$.
At each iteration the interior set $\mathcal{X}$ consists of $5\times10^4$
points drawn uniformly from the shrunken shell $r_0+\dx<\|x\|_2<r_1-\dx$, so that the stencil of every interior point stays inside $\Omega$, and boundary conditions are enforced on
$\mathcal{X}_b$ of $5\times 10^4$ points split evenly between the two spheres of $\partial\Omega$.
Following the multi-level strategy of Section~\ref{sec:numerical algos}, we
minimize the residual over a sequence of decreasing discretization steps
\begin{equation}
\label{eq:HJI-schedule}
  (\dx_1,\dx_2,\dx_3,\dx_4) = (0.2,\ 0.1,\ 0.05,\ {0.02}),
\end{equation}
each level warm-started from the network trained at the preceding one. 
The learning rate, initially $10^{-3}$, is reduced as the grid is refined, 
and the number of Adam iterations grows from about $10^4$ at the coarsest 
level to $4\times10^4$ at the finest.
{For the last step ($\dx=0.02$), we switched to centered  differences 
$D_ku_i:=(u_{i+e_k}-u_{i-e_k})/(2\dx_k)$ for the first-order derivative since, for this value of $\dx$, the overall scheme is now 
monotone due to the diffusion terms. }
\paragraph{Results.}
The trained network matches the reference closely (Figure~\ref{fig:HJI-conv-dim8} center and right), and
both the residual loss and the $\ell^2$ error against the reference decrease
by about an order of magnitude over training
(Figure~\ref{fig:HJI-conv-dim8}, left).
{Here the $\ell^2$ mean error is computed by using the reference value and considering a 2D slice in the first two variables $(x_1,x_2)\in[-2,2]^2$ and $x_3=\cdots=x_8=0$, on a $200^2$ cartesian grid.}

\begin{figure}[htbp]
\end{figure}

\section{Conclusion}

We have developed and analyzed a method for solving Hamilton--Jacobi equations in higher dimensions by minimizing the squared residual of monotone finite-difference discretizations, with the minimization carried out by gradient-based training of a neural network. The framework covers first-order equations as well as fully nonlinear second-order equations, and it applies to the standard monotone schemes, including Lax--Friedrichs and upwind.

The well-posedness of this approach rests on two results. First, under either a properness condition or a uniform ellipticity condition, every critical point of the residual loss is its global minimizer and solves the monotone scheme together with the prescribed Dirichlet conditions, so gradient-based optimization cannot stall at a spurious solution. Second, the error of any approximation is controlled by its residual through computable \emph{a~posteriori} bounds, so the trained network carries a checkable error estimate. Under the properness condition this control is dimension-robust in the $\ell^2$ norm and takes the form of a Polyak--\L{}ojasiewicz inequality, from which the gradient flow converges linearly, at a rate set by the local stencil rather than by the number of grid points (overcoming the curse of dimensionality). The same theory extends to time-dependent problems through a single space-time residual, where the implicit and explicit Euler steppers are both well-posed and, within the CFL regime required for dimension-robust training, equivalent; our analysis therefore favors neither.

Because the conditioning of the optimization deteriorates as the grid is refined, we train in several levels, fitting the network on coarse grids and warm-starting it on finer ones; since the network is its own interpolant, no grid-transfer operator is required between levels. The resulting method scales to high dimensions, as we demonstrate on eikonal equations in dimensions four and eight, front propagation in five space dimensions, nonlinear advection past an obstacle in eight dimensions, and a Hamilton--Jacobi--Isaacs equation in dimension height, whose Hamiltonian is neither convex nor concave, arising from a stochastic differential game.

Several questions remain open. Since monotone schemes are only first-order accurate, the multi-level structure invites the use of a higher-order scheme, such as ENO or WENO \cite{osher1991high,jiang2000weighted}, at the finest level, which would recover higher-order accuracy while preserving the well-posedness of the monotone scheme on the coarser levels.
It also remains to quantify whether, and to what extent, the multi-level method mitigates the curse of dimensionality, as the scaling of cost with dimension is not yet understood. Finally, the efficient selection of collocation points in high dimensions remains challenging, since concentration of measure makes uniform sampling increasingly wasteful as the dimension grows; first steps in this direction were taken in \cite{esteve2025finite}.

\section*{Acknowledgments}

Richard Tsai is supported partially by National Science Foundation grant DMS-2513857.
Part of the research is also supported by the Swedish Research Council under grant no. 2021-06594, while the author was in residence at Institut Mittag-Leffler in Djursholm, Sweden during the Fall 2025 semester.
Carlos Esteve-Yagüe is supported by Agencia Estatal de Investigación under the Ram\'on y Cajal 2022 grant RYC2022-035966-I (Spain).
Olivier Bokanowski and Carlos Esteve-Yagüe benefited from the support of the FMJH Program Gaspard Monge for optimization and operations research and their interactions with data science.

\bibliographystyle{abbrv}
\bibliography{references.bib}

\end{document}